\pdfoutput=1
\RequirePackage{ifpdf}
\ifpdf % We are running pdfTeX in pdf mode
\documentclass[pdftex]{sigma}
\else
\documentclass{sigma}
\fi

\numberwithin{equation}{section}

\newtheorem{Theorem}{Theorem}[section]
\newtheorem{Lemma}[Theorem]{Lemma}
\newtheorem{Proposition}[Theorem]{Proposition}
{\theoremstyle{definition}
\newtheorem{Definition}[Theorem]{Definition}
}

\usepackage{sgmlcmpt}
\usepackage{multirow}

\usepackage{tikz}
\usetikzlibrary{decorations.markings}
\usetikzlibrary{decorations.pathreplacing}
\usetikzlibrary{arrows,shapes,positioning}
\tikzstyle directed=[postaction={decorate,decoration={markings,
 mark=at position #1 with {\arrow{>}}}}]
\tikzstyle rdirected=[postaction={decorate,decoration={markings,
 mark=at position #1 with {\arrow{<}}}}]

\tikzset{anchorbase/.style={baseline={([yshift=-0.5ex]current bounding box.center)}},
anchorzero/.style={baseline={([yshift=-0.5ex]0,0)}},
arrowinthemiddle/.style={postaction=decorate,decoration={markings,mark=at position 0.5 with {\arrow{>}}}},
arrowinthemiddlerev/.style={postaction=decorate,decoration={markings,mark=at position 0.5 with {\arrow{<}}}},
cross line/.style={preaction={draw=white,line width=4pt,-}},
int/.style={thick},
zero/.style={thin,dotted},
uno/.style={thin}}

\usepackage[all]{xy}
\SelectTips{cm}{}

\newcommand{\onenn}[1]{{\mathbf 1}_{#1}}
\newcommand{\onel}{{\mathbf 1}_{\lambda}}
\newcommand{\onell}[1]{{\mathbf 1}_{#1}}
\newcommand{\onelp}{{\mathbf 1}_{\lambda'}}
\def\cal#1{\mathcal{#1}}%
\newcommand\Ts{{T''_{i}}}
\newcommand{\Tsw}[1]{{T''_{#1}}}
\newcommand\Tp{{T'_{i}}}
\newcommand{\Tpw}[1]{{T'_{#1}}}
\newcommand{\T}{T_{i}}
\newcommand{\Tw}[1]{T_{#1}}

\newcommand\E{{\sf{E}}}
\newcommand\F{{\sf{F}}}

\newcommand{\U}{\dot{{\bf U}}}
\newcommand{\Ucat}{\cal{U}}
\newcommand{\sln}{\mathfrak{sl}_n}
\newcommand{\slm}{\mathfrak{sl}_m}
\newcommand{\slnn}[1]{\mathfrak{sl}_{#1}}
\newcommand{\gln}{\mathfrak{gl}_n}
\newcommand{\glm}{\mathfrak{gl}_m}

\newcommand{\hsln}{\hat{\sln}}
\newcommand{\hslm}{\hat{\slm}}

\newcommand{\Bfoam}[3][N]{#2\cat{BFoam}_{#3}(N)}
\newcommand{\cat}[1]{\ensuremath{\mbox{\bfseries {\upshape {#1}}}}}

\let\hat=\widehat
\let\tilde=\widetilde

\let\phi=\varphi
\let\epsilon=\varepsilon

\usepackage{bbm}
\def\C{{\mathbb C}}
\def\R{{\mathbb R}}
\def\Z{{\mathbb Z}}
\def\Q{{\mathbb Q}}

\def\1{\mathbbm{1}}%

\newcommand{\bbpef}[1]{\xybox{%
 (-6,0)*{};
 (6,0)*{};
 (-4,0)*{}="t1";
 (4,0)*{}="t2";
 "t1";"t2" **\crv{(-4,-6) & (4,-6)}; ?(.15)*\dir{>} ?(.9)*\dir{>}
 ?(.5)*\dir{}+(0,-2)*{\scriptstyle{#1}};
}}
\newcommand{\bbpfe}[1]{\xybox{%
 (-6,0)*{};
 (6,0)*{};
 (-4,0)*{}="t1";
 (4,0)*{}="t2";
 "t2";"t1" **\crv{(4,-6) & (-4,-6)}; ?(.15)*\dir{>} ?(.9)*\dir{>} ?(.5)*\dir{}+(0,-2)*{\scriptstyle{#1}};
}}

\newcommand{\bbcfe}[1]{\xybox{%
 (-6,0)*{};
 (6,0)*{};
 (-4,0)*{}="t1";
 (4,0)*{}="t2";
 "t1";"t2" **\crv{(-4,6) & (4,6)}; ?(.15)*\dir{>} ?(.9)*\dir{>}
 ?(.5)*\dir{}+(0,2)*{\scriptstyle{#1}};
}}

\newcommand{\bbcef}[1]{\xybox{%
 (-6,0)*{};
 (6,0)*{};
 (-4,0)*{}="t1";
 (4,0)*{}="t2";
 "t2";"t1" **\crv{(4,6) & (-4,6)}; ?(.15)*\dir{>}
 ?(.9)*\dir{>} ?(.5)*\dir{}+(0,2)*{\scriptstyle{#1}};
}}

\newcommand{\Efoam}[1][1]{
\xy
(0,0)*{
\begin{tikzpicture} [fill opacity=0.2, decoration={markings,
 mark=at position 0.6 with {\arrow{>}}; }, scale=#1]
 \filldraw [fill=red] (-2,0) rectangle (2,3);
 \draw[very thick] (-2,3) -- (2,3);
 \draw[very thick] (-1,-1) -- (3,-1);
 \draw[ultra thick, red] (0,0) -- (0,3);
 \filldraw [fill=blue] (1,-1) -- (1,2) -- (0,3) -- (0,0) -- cycle;
 \draw[very thick, postaction={decorate}] (1,-1) -- (0,0);
 \draw[very thick, postaction={decorate}] (1,2) -- (0,3);
 \filldraw [fill=red] (-1,-1) rectangle (3,2);
 \draw[very thick] (-1,2) -- (3,2);
 \draw[very thick] (-2,0) -- (2,0);
 \draw[ultra thick, red] (1,-1) -- (1,2);
\end{tikzpicture}};
\endxy
}

\newcommand{\dotEfoam}[1][1]{
\xy
(0,0)*{
\begin{tikzpicture} [fill opacity=0.2, decoration={markings,
 mark=at position 0.6 with {\arrow{>}}; }, scale=#1]
 \filldraw [fill=red] (-2,0) rectangle (2,3);
 \draw[very thick] (-2,3) -- (2,3);
 \draw[very thick] (-1,-1) -- (3,-1);
 \draw[ultra thick, red] (0,0) -- (0,3);
 \node [opacity=1] at (0.5,1) {$\bullet$};
 \filldraw [fill=blue] (1,-1) -- (1,2) -- (0,3) -- (0,0) -- cycle;
 \draw[very thick, postaction={decorate}] (1,-1) -- (0,0);
 \draw[very thick, postaction={decorate}] (1,2) -- (0,3);
 \filldraw [fill=red] (-1,-1) rectangle (3,2);
 \draw[very thick] (-1,2) -- (3,2);
 \draw[very thick] (-2,0) -- (2,0);
 \draw[ultra thick, red] (1,-1) -- (1,2);
\end{tikzpicture}};
\endxy
}

\newcommand{\crossingEEfoam}[1][1]{
\xy
(0,0)*{
\begin{tikzpicture} [fill opacity=0.2, decoration={markings,
 mark=at position 0.6 with {\arrow{>}}; }, scale=#1]
 \filldraw [fill=red] (2,2) rectangle (-2,-1);
 \draw[very thick](-2,2)--(2,2);
 \draw[ultra thick, red] (-.75,2) arc (180:360:.75);
 \draw[very thick] (-2,-1)--(2,-1);
 \draw[ultra thick, red] (-.75,-1) arc (180:0:.75);
 \path [fill=blue] (1.75,-2) arc (0:60:0.75) -- (0.375,-0.35) arc (60:0:0.75);
 \path [fill=blue] (0.375,-0.35) arc (60:180:0.75) -- (0.25,-2) arc (180:60:0.75);
 \draw[very thick, postaction={decorate}] (.25,-2) -- (-.75,-1);
 \draw[very thick, postaction={decorate}] (1.75,-2) --(.75,-1);
 \filldraw[ultra thick, red] [fill=blue] (1,-1.25) -- (0,-0.25) -- (0,1.25) -- (1,0.25) -- cycle;
 \path [fill=blue] (.25,1) arc (180:240:0.75) -- (-0.375,1.35) arc (240:180:0.75);
 \path [fill=blue] (0.625,0.35) arc (240:360:0.75) -- (.75,2) arc (0:-120:0.75);
 \filldraw [fill=red] (-1,-2) rectangle (3,1);
 \draw[very thick] (-1,1)--(3,1);
 \draw[very thick] (-1,-2)--(3,-2);
 \draw[ultra thick, red] (.25,-2) arc (180:0:.75);
 \draw[ultra thick, red] (.25,1) arc (180:360:.75);
 \draw[very thick, postaction={decorate}] (.25,1) -- (-.75,2);
 \draw[very thick, postaction={decorate}] (1.75,1) --(.75,2);
\end{tikzpicture}};
\endxy
}

\newcommand{\crossingEoneEtwofoam}[1][1]{
\xy
(0,0)*{
\begin{tikzpicture} [fill opacity=0.2, decoration={markings, mark=at position 0.6 with {\arrow{>}};}, scale=#1]
 \filldraw [fill=red] (-3,1) rectangle (1,4);
 \draw[very thick] (-3,1) -- (1,1);
 \draw[very thick] (-3,4) -- (1,4);
 \draw [very thick, red] (-1.75,4) .. controls (-1.75,3) and (-.25,2) .. (-.25,1);
 \draw[very thick, postaction={decorate}] (-.75,3) -- (-1.75,4);
 \draw[very thick, postaction={decorate}] (.75,0) -- (-.25,1);
 \path [fill=blue] (-1.75,4) .. controls (-1.75,3) and (-.25,2) .. (-.25,1) --
 (.75,0) .. controls (.75,1) and (-.75,2) .. (-.75,3);
 \filldraw [fill=red] (-2,0) rectangle (2,3);
 \draw[very thick] (-2,0) -- (2,0);
 \draw[very thick] (-2,3) -- (2,3);
 \draw [very thick, red] (-.75,3) .. controls (-.75,2) and (.75,1) .. (.75,0);
 \draw [very thick, red] (.75,3) .. controls (.75,2) and (-.75,1) .. (-.75,0);
 \draw[very thick, postaction={decorate}] (1.75,2) -- (.75,3);
 \draw[very thick, postaction={decorate}] (.25,-1) -- (-.75,0);
 \path [fill=blue] (.75,3) .. controls (.75,2) and (-.75,1) .. (-.75,0) --
 (.25,-1) .. controls (.25,0) and (1.75,1) .. (1.75,2);
 \filldraw [fill=red] (-1,-1) rectangle (3,2);
 \draw[very thick] (-1,-1) -- (3,-1);
 \draw[very thick] (-1,2) -- (3,2);
 \draw[very thick, red] (.25,-1) .. controls (.25,0) and (1.75,1) .. (1.75,2);
\end{tikzpicture}}
\endxy
}

\newcommand{\crossingEtwoEonefoam}[1][1]{
\xy
(0,0)*{
\begin{tikzpicture} [fill opacity=0.2, decoration={markings, mark=at position 0.6 with {\arrow{>}};}, scale=#1]
 \filldraw [fill=red] (-3,1) rectangle (1,4);
 \draw[very thick] (-3,1) -- (1,1);
 \draw[very thick] (-3,4) -- (1,4);
 \draw [very thick, red] (-.25,4) .. controls (-.25,3) and (-1.75,2) .. (-1.75,1);
 \draw[very thick, postaction={decorate}] (.75,3) -- (-.25,4);
 \draw[very thick, postaction={decorate}] (-.75,0) -- (-1.75,1);
 \path [fill=blue] (-.25,4) .. controls (-.25,3) and (-1.75,2) .. (-1.75,1) --
 (-.75,0) .. controls (-.75,1) and (.75,2) .. (.75,3);
 \filldraw [fill=red] (-2,0) rectangle (2,3);
 \draw[very thick] (-2,0) -- (2,0);
 \draw[very thick] (-2,3) -- (2,3);
 \draw [very thick, red] (-.75,3) .. controls (-.75,2) and (.75,1) .. (.75,0);
 \draw [very thick, red] (.75,3) .. controls (.75,2) and (-.75,1) .. (-.75,0);
 \draw[very thick, postaction={decorate}] (.25,2) -- (-.75,3);
 \draw[very thick, postaction={decorate}] (1.75,-1) -- (.75,0);
 \path [fill=blue] (-.75,3) .. controls (-.75,2) and (.75,1) .. (.75,0) --
 (1.75,-1) .. controls (1.75,0) and (.25,1) .. (.25,2);
 \filldraw [fill=red] (-1,-1) rectangle (3,2);
 \draw[very thick] (-1,-1) -- (3,-1);
 \draw[very thick] (-1,2) -- (3,2);
 \draw[very thick, red] (1.75,-1) .. controls (1.75,0) and (.25,1) .. (.25,2);
\end{tikzpicture}}
\endxy
}

\newcommand{\crossingEthreeEonefoam}[1][1]{
\xy
(0,0)*{
\begin{tikzpicture} [fill opacity=0.2,decoration={markings, mark=at position 0.6 with {\arrow{>}}; }, scale=#1]
 \filldraw [fill=red] (-4,2) rectangle (0,6);
 \draw[very thick] (-4,2) -- (0,2);
 \draw[very thick] (-4,6) -- (0,6);
 \draw[very thick, red] (-1.25,6) .. controls (-1.25,5) and (-2.75,3) .. (-2.75,2);
 \draw[very thick, postaction={decorate}] (-1.75,1) -- (-2.75,2);
 \draw[very thick, postaction={decorate}] (-.25,5) -- (-1.25,6);
 \path [fill=blue] (-1.25,6) .. controls (-1.25,5) and (-2.75,3) .. (-2.75,2) --
 (-1.75,1) .. controls (-1.75,2) and (-.25,4) .. (-.25,5);
 \filldraw [fill=red] (-3,1) rectangle (1,5);
 \draw[very thick] (-3,1) -- (1,1);
 \draw[very thick] (-3,5) -- (1,5);
 \draw[very thick, red] (-1.75,1) .. controls (-1.75,2) and (-.25,4) .. (-.25,5);
 \filldraw [fill=red] (-2,0) rectangle (2,4);
 \draw[very thick] (-2,0) -- (2,0);
 \draw[very thick] (-2,4) -- (2,4);
 \draw[very thick, red] (-.75,4) .. controls (-.75,3) and (.75,1) .. (.75,0);
 \draw[very thick, postaction={decorate}] (.25,3) -- (-.75,4);
 \draw[very thick, postaction={decorate}] (1.75,-1) -- (.75,0);
 \path [fill=blue] (-.75,4) .. controls (-.75,3) and (.75,1) .. (.75,0) --
 (1.75,-1) .. controls (1.75,0) and (.25,2) .. (.25,3);
 \filldraw [fill=red] (-1,-1) rectangle (3,3);
 \draw[very thick] (-1,-1) -- (3,-1);
 \draw[very thick] (-1,3) -- (3,3);
 \draw[very thick, red] (1.75,-1) .. controls (1.75,0) and (.25,2) .. (.25,3);
\end{tikzpicture}}
\endxy
}

\newcommand{\crossingEoneEthreefoam}[1][1]{
\xy
(0,0)*{
\begin{tikzpicture} [fill opacity=0.2,decoration={markings, mark=at position 0.6 with {\arrow{>}}; }, scale=#1]
 \filldraw [fill=red] (-4,2) rectangle (0,6);
 \draw[very thick] (-4,2) -- (0,2);
 \draw[very thick] (-4,6) -- (0,6);
 \draw[very thick, red] (-2.75,6) .. controls (-2.75,5) and (-1.25,3) .. (-1.25,2);
 \draw[very thick, postaction={decorate}] (-.25,1) -- (-1.25,2);
 \draw[very thick, postaction={decorate}] (-1.75,5) -- (-2.75,6);
 \path [fill=blue] (-2.75,6) .. controls (-2.75,5) and (-1.25,3) .. (-1.25,2) --
 (-.25,1) .. controls (-.25,2) and (-1.75,4) .. (-1.75,5);
 \filldraw [fill=red] (-3,1) rectangle (1,5);
 \draw[very thick] (-3,1) -- (1,1);
 \draw[very thick] (-3,5) -- (1,5);
 \draw[very thick, red] (-.25,1) .. controls (-.25,2) and (-1.75,4) .. (-1.75,5);
 \filldraw [fill=red] (-2,0) rectangle (2,4);
 \draw[very thick] (-2,0) -- (2,0);
 \draw[very thick] (-2,4) -- (2,4);
 \draw[very thick, red] (.75,4) .. controls (.75,3) and (-.75,1) .. (-.75,0);
 \draw[very thick, postaction={decorate}] (1.75,3) -- (.75,4);
 \draw[very thick, postaction={decorate}] (.25,-1) -- (-.75,0);
 \path [fill=blue] (.75,4) .. controls (.75,3) and (-.75,1) .. (-.75,0) --
 (.25,-1) .. controls (.25,0) and (1.75,2) .. (1.75,3);
 \filldraw [fill=red] (-1,-1) rectangle (3,3);
 \draw[very thick] (-1,-1) -- (3,-1);
 \draw[very thick] (-1,3) -- (3,3);
 \draw[very thick, red] (.25,-1) .. controls (.25,0) and (1.75,2) .. (1.75,3);
\end{tikzpicture}}
\endxy
}

\newcommand{\cupFEfoam}[1][1]{
\xy
(0,0)*{
\begin{tikzpicture} [fill opacity=0.2, decoration={markings,
 mark=at position 0.6 with {\arrow{>}}; }, scale=#1]
 \filldraw [fill=red] (-2,0) rectangle (2,3);
 \draw[very thick](-2,3)--(2,3);
 \draw[very thick] (-2,0)--(2,0);
 \draw[ultra thick, red] (-.75,3) arc (180:360:.75);
 \path [fill=blue] (.25,2) arc (180:240:0.75) -- (-0.375,2.35) arc (240:180:0.75);
 \path [fill=blue] (0.625,1.35) arc (240:360:0.75) -- (.75,3) arc (0:-120:0.75);
 \filldraw [fill=red] (-1,-1) rectangle (3,2);
 \draw[very thick] (-1,-1)--(3,-1);
 \draw[very thick] (-1,2)--(3,2);
 \draw[very thick, postaction={decorate}] (-.75,3) -- (.25,2);
 \draw[very thick, postaction={decorate}] (1.75,2) --(.75,3);
 \draw[ultra thick, red] (.25,2) arc (180:360:.75);
\end{tikzpicture}};
\endxy
}

\newcommand{\cupEFfoam}[1][1]{
\xy
(0,0)*{
\begin{tikzpicture} [fill opacity=0.2, decoration={markings,
 mark=at position 0.6 with {\arrow{>}}; }, scale=#1]
 \filldraw [fill=red] (-2,0) rectangle (2,3);
 \draw[very thick](-2,3)--(2,3);
 \draw[very thick] (-2,0)--(2,0);
 \draw[ultra thick, red] (-.75,3) arc (180:360:.75);
 \path [fill=blue] (.25,2) arc (180:240:0.75) -- (-0.375,2.35) arc (240:180:0.75);
 \path [fill=blue] (0.625,1.35) arc (240:360:0.75) -- (.75,3) arc (0:-120:0.75);
 \filldraw [fill=red] (-1,-1) rectangle (3,2);
 \draw[very thick] (-1,-1)--(3,-1);
 \draw[very thick] (-1,2)--(3,2);
 \draw[very thick, postaction={decorate}] (.25,2) -- (-.75,3) ;
 \draw[very thick, postaction={decorate}] (.75,3) -- (1.75,2);
 \draw[ultra thick, red] (.25,2) arc (180:360:.75);
\end{tikzpicture}};
\endxy
}

\newcommand{\capFEfoam}[1][1]{
\xy
(0,0)*{
\begin{tikzpicture} [fill opacity=0.2, decoration={markings,
 mark=at position 0.6 with {\arrow{>}}; }, scale=#1]
 \filldraw [fill=red] (-2,0) rectangle (2,3);
 \draw[very thick](-2,3)--(2,3);
 \draw[very thick] (-2,0)--(2,0);
 \draw[ultra thick, red] (-.75,0) arc (180:0:.75);
 \path [fill=blue] (1.75,-1) arc (0:60:0.75) -- (0.375,0.65) arc (60:0:0.75);
 \path [fill=blue] (0.375,0.65) arc (60:180:0.75) -- (0.25,-1) arc (180:60:0.75);
 \filldraw [fill=red] (-1,-1) rectangle (3,2);
 \draw[very thick] (-1,-1)--(3,-1);
 \draw[very thick] (-1,2)--(3,2);
 \draw[ultra thick, red] (.25,-1) arc (180:0:.75);
 \draw[very thick, postaction={decorate}] (-.75,0) -- (.25,-1);
 \draw[very thick, postaction={decorate}] (1.75,-1) --(.75,0);
 \end{tikzpicture}};
\endxy
}

\newcommand{\capEFfoam}[1][1]{
\xy
(0,0)*{
 \begin{tikzpicture} [fill opacity=0.2, decoration={markings,
 mark=at position 0.6 with {\arrow{>}}; }, scale=#1]
 \filldraw [fill=red] (-2,0) rectangle (2,3);
 \draw[very thick](-2,3)--(2,3);
 \draw[very thick] (-2,0)--(2,0);
 \draw[ultra thick, red] (-.75,0) arc (180:0:.75);
 \path [fill=blue] (1.75,-1) arc (0:60:0.75) -- (0.375,0.65) arc (60:0:0.75);
 \path [fill=blue] (0.375,0.65) arc (60:180:0.75) -- (0.25,-1) arc (180:60:0.75);
 \filldraw [fill=red] (-1,-1) rectangle (3,2);
 \draw[very thick] (-1,-1)--(3,-1);
 \draw[very thick] (-1,2)--(3,2);
 \draw[ultra thick, red] (.25,-1) arc (180:0:.75);
 \draw[very thick, postaction={decorate}] (.25,-1) -- (-.75,0) ;
 \draw[very thick, postaction={decorate}] (.75,0) -- (1.75,-1);
\end{tikzpicture}};
\endxy
}

\DeclareMathOperator{\End}{End}
\begin{document}

\newcommand{\arXivNumber}{1309.4985}

\allowdisplaybreaks

\renewcommand{\thefootnote}{$\star$}

\renewcommand{\PaperNumber}{030}

\FirstPageHeading

\ShortArticleName{Skein Modules from Skew Howe Duality and Affine Extensions}

\ArticleName{Skein Modules from Skew Howe Duality\\
and Affine Extensions\footnote{This paper is a~contribution to the Special Issue on New Directions in Lie
Theory. The full collection is available
at \href{http://www.emis.de/journals/SIGMA/LieTheory2014.html}{http://www.emis.de/journals/SIGMA/LieTheory2014.html}}}

\Author{Hoel QUEFFELEC}

\AuthorNameForHeading{H.~Queffelec}

\Address{Mathematical Sciences Institute, The Australian National University,\\
J.D.~27 Union Lane, Acton ACT 2601, Australia}
\Email{\href{mailto:hoel.queffelec@anu.edu.au}{hoel.queffelec@anu.edu.au}}
\URLaddress{\url{http://maths-people.anu.edu.au/~queffelech/}}

\ArticleDates{Received July 22, 2014, in f\/inal form March 30, 2015; Published online April 15, 2015}

\Abstract{We show that we can release the rigidity of the skew Howe duality process for~$\sln$ knot invariants
by rescaling the quantum Weyl group action,
and recover skein modules for web-tangles.
This skew Howe duality phenomenon can be extended to the af\/f\/ine $\slm$ case,
corresponding to looking at tangles embedded in a~solid torus.
We investigate the relations between the invariants constructed by evaluation representations (and af\/f\/inization of them) and usual skein modules,
and give tools for interpretations of annular skein modules as sub-algebras of intertwiners for particular $U_q(\sln)$ representations.
The categorif\/ication proposed in a~joint work with A.~Lauda and D.~Rose also admits a direct extension in the af\/f\/ine case.}

\Keywords{skein modules; quantum groups; annulus; af\/f\/ine quantum groups}

\Classification{81R50; 17B37; 17B67; 57M25; 57M27}

\renewcommand{\thefootnote}{\arabic{footnote}}
\setcounter{footnote}{0}

\section{Introduction}

\subsection{Webs and skew-Howe duality}

Cautis,
Kamnitzer and Licata~\cite{Cautis, CKL} recently introduced a~reformulation of the $\sln$ Reshetikhin--Turaev invariants
for knots and links based on the \emph{quantum skew Howe duality}.
This duality phenomenon involves two commuting actions of $U_q(\mathfrak{sl}_m)$ and $U_q(\mathfrak{sl}_n)$
on the quantum exterior algebra ${\bigwedge}_q(\C^n\otimes \C^m)$,
where $n$ corresponds to the $\sln$-invariants we look at,
and $m$ governs the braiding of $m$-fold tensor products of $\sln$-representations.
In this framework,
braidings arise from the so-called \emph{quantum Weyl group action}~\cite{KamnTin, Lus4} on $U_q(\slm)$.

This new process is naturally related to the concept of \emph{webs},
which emerge from the study of~$\sln$ knot invariants and describe intertwiners of $\sln$-representations
(see~\cite{Kim,Kup,Mor} for detailed studies of the spider categories they form).
For each $n$,
$\sln$ webs are trivalent oriented graphs with edges labeled with integers in $\{1,\dots,n\}$.
At each vertex,
the sum of the indices of the incoming edges equals the sum of the indices of the outgoing edges:
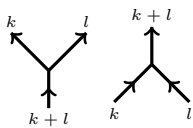
\begin{gather*}
\xy
(0,0)*{
\begin{tikzpicture} [scale=.5]
\draw [very thick, directed=.55] (1,0) -- (1,1);
\draw [very thick, directed=1] (1,1) -- (0,2);
\draw [very thick, directed=1] (1,1) -- (2,2);
\node at (1,-.3) {\tiny $k+l$};
\node at (0,2.3) {\tiny $k$};
\node at (2,2.3) {\tiny $l$};
\end{tikzpicture}};
\endxy
\xy
(0,0)*{
\begin{tikzpicture} [scale=.5]
\draw [very thick, directed=1] (1,1) -- (1,2);
\draw [very thick, directed=.55] (0,0) -- (1,1);
\draw [very thick, directed=.55] (2,0) -- (1,1);
\node at (1,2.3) {\tiny $k+l$};
\node at (0,-.3) {\tiny $k$};
\node at (2,-.3) {\tiny $l$};
\end{tikzpicture}};
\endxy
\end{gather*}

Here is an example of a~web:
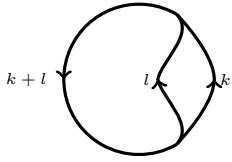
\begin{gather*}
\xy
(0,0)*{
\begin{tikzpicture} [scale=.5,decoration={markings,
mark=at position .5 with {\arrow{>}}; }]
\draw[postaction={decorate}, very thick] (1,1.73) arc (60:300:2);
\draw[postaction={decorate}, very thick] (1,-1.73) .. controls (1.5,-1.2) and (.6,-.5) .. (.5,0) -- (.5,0) .. controls (.6,.5) and (1.5,1.2) .. (1,1.73);
\draw[postaction={decorate}, very thick] (1,-1.73) .. controls (1.5,-1.2) and (2,-.5) .. (2,0) -- (2,0) .. controls (2,.5) and (1.5,1.2) .. (1,1.73);
\node at (-3,0) {\tiny $k+l$};
\node at (2.3,0) {\tiny $k$};
\node at (.2,0) {\tiny $l$};
\end{tikzpicture}};
\endxy
\end{gather*}

These diagrams are to be understood up to some local relations (see~\cite{CKM} for example),
which are a~diagrammatic analogue of relations between morphisms of $U_q(\sln)$-representations that are at the origin of the def\/inition of webs.
Note that there are more ref\/ined notions of $\sln$-webs,
in particular concerning what to do with $n$-labeled strands.
Indeed,
a $k$-labeled strand corresponds to the $k$-th exterior power $\bigwedge_q^k(\C^n)$ of the standard $U_q(\sln)$ representation $\C^n$.
The $0$-th power is the trivial representation,
and it appears natural not to depict it in webs.
Similarly,
the maximal exterior power $\bigwedge_q^n(\C^n)$ is just the trivial representation,
and it is usually forgotten as well (which comes with a~correspondence between an edge labeled by $k$ and the same edge
with opposite orientation labeled by $n-k$,
see for example~\cite{MS2}).
However,
it appears that this maximal exterior power plays a~non-trivial role in some places,
in particular when looking at categorif\/ication questions.
An heuristic interpretation of this could be the fact that this representation corresponds to the determinant representation,
which is indeed a~trivial $\sln$-representation,
but is not a~trivial $\gln$-representation.
This non-triviality has been encoded by \emph{tags} in some places~\cite{Mor,CKM},
and applications to categorif\/ied knot invariants using these tags in the $\slnn{2}$ case can be found in a~work by Clark,
Morrison and Walker~\cite{CMW}.
One can also choose to keep all the $n$-edges (which we will then depict doubled).
Although the dif\/ference is at f\/irst sight minimal on the level of webs,
it seems to play an important role at the categorif\/ied level,
as suggested by Blanchet's work~\cite{Blan} and developed in~\cite{LQR1}.
We sometimes refer to these webs as \emph{enhanced}.

For example,
the enhanced web:
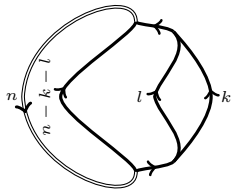
\begin{gather*}
\xy
(0,0)*{
\begin{tikzpicture} [scale=.5,decoration={markings,
mark=at position .55 with {\arrow{>}}; }]
\draw[postaction={decorate},very thick] (1,1.73) .. controls (.8,1.9) and (0,1.9) .. (0,2);
\draw[postaction={decorate}, very thick] (0,-2) .. controls (0,-1.9) and (.8,-1.9) .. (1,-1.73);
\draw[postaction={decorate}, very thick] (0,-2) .. controls (0,-1.7) and (-2,-.5) .. (-2,0) -- (-2,0) .. controls (-2,.5) and (0,1.7) ..(0,2);
\draw [postaction={decorate}, double] (0,2) .. controls (0,3) and (-2.9,2) .. (-3,0) -- (-3,0) .. controls (-2.9,-2) and (0,-3) .. (0,-2);
\draw[postaction={decorate}, very thick] (1,-1.73) .. controls (1.5,-1.2) and (.6,-.5) .. (.5,0) -- (.5,0) .. controls (.6,.5)
 and (1.5,1.2) .. (1,1.73);
\draw[postaction={decorate}, very thick] (1,-1.73) .. controls (1.5,-1.2) and (2,-.5) .. (2,0) -- (2,0) .. controls (2,.5) and (1.5,1.2) .. (1,1.73);
\node at (-3.3,0) {\tiny $n$};
\node [rotate=90] at (-2.4,0) {\tiny $n-k-l$};
\node at (2.4,0) {\tiny $k$};
\node at (.1,0) {\tiny $l$};
\end{tikzpicture}};
\endxy
\end{gather*}
is represented in the tagged version of webs as
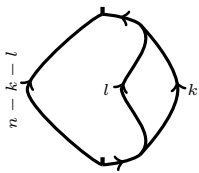
\begin{gather*}
\xy
(0,0)*{
\begin{tikzpicture} [scale=.5,decoration={markings,
mark=at position .55 with {\arrow{>}}; }]
\draw[postaction={decorate},very thick] (1,1.73) .. controls (.8,1.9) and (.2,2) .. (0,2);
\draw[postaction={decorate}, very thick] (0,-2) .. controls (.2,-2) and (.8,-1.9) .. (1,-1.73);
\draw[postaction={decorate}, very thick] (0,-2) .. controls (-.2,-2) and (-2,-.5) .. (-2,0) -- (-2,0) .. controls (-2,.5) and (-.2,2) ..(0,2);
\draw [ultra thick] (0,2) -- (0,2.2);
\draw [ultra thick] (0,-2) -- (0,-1.8);
\draw[postaction={decorate}, very thick] (1,-1.73) .. controls (1.5,-1.2) and (.6,-.5) .. (.5,0) -- (.5,0) .. controls (.6,.5)
 and (1.5,1.2) .. (1,1.73);
\draw[postaction={decorate}, very thick] (1,-1.73) .. controls (1.5,-1.2) and (2,-.5) .. (2,0) -- (2,0) .. controls (2,.5) and (1.5,1.2) .. (1,1.73);
\node [rotate=90] at (-2.4,0) {\tiny $n-k-l$};
\node at (2.4,0) {\tiny $k$};
\node at (.1,0) {\tiny $l$};
\end{tikzpicture}};
\endxy
\end{gather*}

The only dif\/ference between the two pictures above lies in the way to deal with the strands decorated with the maximum exterior power
of the fundamental representation: while we keep them completely in the f\/irst case,
they only appear locally in the second case.

The Jones polynomial and its $\sln$ analogues naturally take place in the spider categories.
Their reformulation in terms of quantum skew Howe duality proved to be a~very powerful tool for understanding these categories,
and led Cautis,
Kamnitzer and Morrison~\cite{CKM} to solve conjectures on generators and relations for categories of representations ${\rm Rep}_q(\sln)$.
Furthermore,
this process admits a~very natural categorif\/ication~\cite{Cautis},
linking~\cite{LQR1} topological categorif\/ications based on skein theory~\cite{BN2, Kh1,Kh2,Kh5} and categorif\/ied
quantum groups~\cite{KL,KL3,KL2,KLMS}.

However,
the skew Howe duality process is quite rigid,
allowing to deal only with \emph{ladder webs},
which are a~particular class of webs with only upward oriented edges.
This is a~generalization to the web case of the notion of upward-oriented tangles,
with the additional requirement that webs are presented in a~rigid structure where strands are either vertical
(the uprights of the ladder) or elementary horizontal pieces (the rungs of the ladder).
For example,
a ladder version of the previous web would be
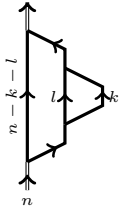
\begin{gather*}
\xy
(0,0)*{
\begin{tikzpicture} [scale=.5]
\draw [double, directed=.55] (0,0) -- (0,.75);
\draw [very thick, directed=.55] (0,.75) -- (0,4.25);
\draw [double, directed=1] (0,4.25) -- (0,5);
\draw [very thick, directed=.55] (0,.75) -- (1,1.25) -- (1,1.75);
\draw [very thick, directed=.55] (1,1.75) -- (1,3.25);
\draw [very thick, directed=.55] (1,1.75) -- (2,2.25) -- (2,2.75) -- (1,3.25);
\draw [very thick, directed=.55] (1,3.25) -- (1,3.75) -- (0,4.25);
\node at (0,-.3) {\tiny $n$};
\node at (.7,2.5) {\tiny $l$};
\node at (2.3,2.5) {\tiny $k$};
\node [rotate=90] at (-.4,2.5) {\tiny $n-k-l$};
\end{tikzpicture}};
\endxy
\end{gather*}

Furthermore,
the relation established by Cautis--Kamnitzer--Licata between the braiding (or $R$-matrix) and the quantum Weyl group action
does not allow to ignore crossings involving $0$-labeled strands.
For example,
the def\/inition of the braidings as it appears in~\cite{Lus4} gives for such crossings a~smoothing map $\Psi$ as follows
\begin{gather*}
\Psi_{\Tp}\left(
\xy
(0,0)*{
\begin{tikzpicture} [scale=.5]
\draw [very thick, directed=1] (0,0) -- (2,2);
\draw [dotted] (2,0) -- (1.2,.8);
\draw [dotted, directed=1] (.8,1.2) -- (0,2);
\node at (0,-.5) {\tiny $k$};
\node at (2,-.5) {\tiny 0};
\end{tikzpicture}};
\endxy
\right)
=
(-1)^kq^{k}
\xy
(0,0)*{
\begin{tikzpicture} [scale=.5]
\draw [very thick, directed=1,directed=.55] (0,0) -- (0,.5) -- (2,1.5) -- (2,2);
\draw [dotted, directed=.55] (2,0) -- (2,1.5);
\draw [dotted, directed=1] (0,.5) -- (0,2);
\node at (0,-.5) {\tiny $k$};
\node at (2,-.5) {\tiny 0};
\end{tikzpicture}};
\endxy
\end{gather*}
while we would like this crossing to be smoothed without creating any coef\/f\/icient,
since we do not want to consider $0$-labeled strands in the skein context.
Similarly,
def\/initions provided in~\cite{Cautis} wouldn't produce any coef\/f\/icients,
but the use of tags in some of the Reidemeister-like web moves produces dif\/f\/iculties.

In this paper,
we f\/ind an appropriate rescaling of the Weyl group action that removes the rigidity in the diagrammatic formulation
of link (and knotted web) invariants: the goal is to f\/ind a~skew-Howe based process from which we can extract smoothing rules
that yield a~skein module for web-tangles.
The next paragraph gives more details about these ideas.

\subsection{Obtaining a~skein module}

We give in this paper a~detailed explanation of the skew-Howe duality process,
focusing on obtaining $\sln$ skein modules from this rather rigid context,
for any value of $n$.
One of the problems that usually appears when looking at a~local crossing in a~skein context is that it can be understood in dif\/ferent ways.
For example,
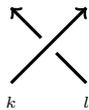
\begin{gather*}
\xy
(0,0)*{
\begin{tikzpicture} [scale=.5]
\draw [very thick, directed=1] (0,0) -- (2,2);
\draw [very thick] (2,0) -- (1.2,.8);
\draw [very thick, directed=1] (.8,1.2) -- (0,2);
\node at (0,-.5) {\tiny $k$};
\node at (2,-.5) {\tiny $l$};
\end{tikzpicture}};
\endxy
\end{gather*}
can be translated as a~positive $(k,l)$ crossing,
or (if we look at it from left to right) as a~nega\-ti\-ve~$(l,k)$ crossing with the $l$ strand in reverse direction.
These two crossings would give rise to dif\/ferent smoothings in their ladder transcriptions.

The ref\/inement we introduce in this paper is based on both a~convenient rescaling of Lusztig's def\/inition of the braidings~\cite{Lus4}
with a~$\glm$-information,
and keeping the whole \emph{enhanced} information of webs,
following ideas of Blanchet~\cite{Blan}.
It is interesting to note that the original construction of Murakami--Ohtsuki--Yamada~\cite{MOY} was actually also using $n$-labeled
strands and is consistent with this presentation.
In the $\mathfrak{sl}_2$ case,
this leads to a~rather unusual presentation of the skein module,
since $2$-labeled crossings produce when smoothed some non-trivial coef\/f\/icients,
so that the smoothing map $\Psi$ would behave as follows
\begin{gather*}
\Psi\left(
\xy
(0,0)*{
\begin{tikzpicture} [scale=.5]
\draw [double, directed=1] (0,0) -- (2,2);
\draw [double] (2,0) -- (1.2,.8);
\draw [double, directed=1] (.8,1.2) -- (0,2);
\end{tikzpicture}};
\endxy
\right)
=
q^{-2} \xy
(0,0)*{
\begin{tikzpicture} [scale=.5]
\draw [double, directed=1] (0,0) .. controls (1,1) .. (0,2);
\draw [double, directed=1] (2,0) .. controls (1,1) .. (2,2);
\end{tikzpicture}};
\endxy
,
\qquad
\Psi\left(
\xy
(0,0)*{
\begin{tikzpicture} [scale=.5]
\draw [double, directed=1] (2,0) -- (0,2);
\draw [double] (0,0) -- (.8,.8);
\draw [double, directed=1] (1.2,1.2) -- (2,2);
\end{tikzpicture}};
\endxy
\right)
=
q^{2}
\xy
(0,0)*{
\begin{tikzpicture} [scale=.5]
\draw [double, directed=1] (2,0) .. controls (1,1) .. (2,2);
\draw [double, directed=1] (0,0) .. controls (1,1) .. (0,2);
\end{tikzpicture}};
\endxy.
\end{gather*}

In the $\mathfrak{sl}_3$ case,
we similarly keep $1$-,
$2$- and $3$-labeled strands,
which produce again dif\/ferent coef\/f\/icients in the smoothings.

The main result is then that there exists a~version of the skew-Howe duality process from which the def\/inition of the braiding
can be used locally to def\/ine an invariant of framed web-tangles.
A~good understanding of the behavior of the braidings back in the representation-theory world is of great help in the proof
of the invariance under Kauf\/fman's web-moves and considerably simplify them,
and also clarif\/ies the categorif\/ication of these results.

\subsection{Af\/f\/ine extensions}

The skew-Howe duality process is based on two commuting actions of $U_q(\sln)$ and $U_q(\slm)$ on the module ${\bigwedge}_q^N(\C^n\otimes \C^m)$.
$U_q(\sln)$ corresponds to the quantum invariant we are looking at,
and we want to keep it unchanged,
but $U_q(\slm)$ appears more as a~parameter,
and we may want to consider extensions of it.
A~f\/irst step is to replace $U_q(\slm)$ by its af\/f\/ine version $U_q(\hslm)$.

Classical representation-theoretic tools tell us that we can extend the action of $U_q(\slm)$ to a~$U_q(\hslm)$ one,
keeping by construction the commutation property with the $U_q(\sln)$ action.
These extensions can be achieved by the process of \emph{evaluation representations}~\cite{ChariPressley}.
This naturally provides knotted-web invariants for the cylinder,
and the only question is to relate these invariants to the usual skein module associated to the thickened surface.
We show that the evaluation representations with a~particular choice of the parameter give the skein module of the f\/illed cylinder,
that can be ref\/ined by passing to the af\/f\/inization of the representations.
These constructions therefore provide a~very natural extension of Jones' construction in the case of web-tangles drawn on the cylinder.

We also investigate better descriptions of the annular skein module in terms of $U_q(\slm)$,
and we take a~f\/irst step towards a~realization as a~sub-algebra of intertwiners for an expli\-cit~$U_q(\sln)$ representation,
which would give to it the same kind of representation-theory f\/lavored interpretation as we have in the linear case.

Many proofs use the fact that relations for $U_q(\slm)$ and $U_q(\hslm)$ locally have the same form.
Thus,
just as at the uncategorif\/ied level,
the categorif\/ication of the skew Howe process provided in~\cite{LQR1} admits a~direct extension to the af\/f\/ine case.

Note that a~recent paper by Mackaay and Thiel~\cite{McKTh} presents a~categorif\/ication of af\/f\/ine $q$-Schur algebras.
Although their paper does not directly deal with annular knots,
it would be interesting to understand its implications in terms of categorif\/ied invariants of annular web-tangles and the links
with categorif\/ied quantum skew-Howe duality.

\section{Skew Howe duality and skein modules}

\subsection{Skew Howe duality}

\subsubsection{Context}

We f\/irst give a~short description of the skew Howe duality phenomenon for usual $\sln$ as explained in~\cite{CKL} and~\cite{CKM}.

We look at the quantum group $U_q(\slm)$ as the $\C[q,q^{-1}]$-algebra generated by the Chevalley elements $E_i$,
$F_i$,
$K_i^{\pm 1}$,
for $1\leq i\leq n-1$,
subject to the relations:
\begin{gather*}
K_iK_i^{-1} = K_i^{-1}K_i = 1,
\qquad
K_iK_j = K_jK_i,
\\
K_iE_jK_i^{-1} = q^{a_{ij}} E_j,
\qquad
K_iF_jK_i^{-1} = q^{-a_{ij}} F_j,
\qquad
E_iF_j - F_jE_i = \delta_{ij} \frac{K_i-K_{i}^{-1}}{q-q^{-1}},
\\
E_i^2E_j-\big(q+q^{-1}\big)E_iE_jE_i+E_jE_i^2 =0
\qquad
\text{if} \ \ j=i\pm 1,
\\
F_i^2F_j-\big(q+q^{-1}\big)F_iF_jF_i+F_jF_i^2=0
\qquad
\text{if} \ \ j=i\pm 1,
\\
E_iE_j=E_jE_i,\qquad F_iF_j = F_jF_i
\qquad
\text{if} \ \ |i-j|>1.
\end{gather*}

The idempotented version of $U_q(\slm)$ will be denoted $\U_q(\slm)$.
Generators are $\onel$,
$E_i\onel$ and~$F_i\onel$,
for all weights $\lambda$.
The unit is then replaced by a~collection of orthogonal idempotents~$\onel$ indexed by the weight lattice of~$\mathfrak{sl}_m$,
\begin{gather*}
\onel \onelp = \delta_{\lambda\lambda'} \onel,
\end{gather*}
such that if $\lambda = (\lambda_1,\lambda_2,
\dots,
\lambda_{m-1})$,
then
\begin{gather*}
K_i\onel =\onel K_i= q^{\lambda_i} \onel,
\qquad
E_i^{}\onel = \onell{\lambda+\alpha_i}E_i,
\qquad F_i\onel = \onell{\lambda-\alpha_i}F_i,
\end{gather*}
where
\begin{gather*}
\lambda +\alpha_i =
\begin{cases}
(\lambda_1+2,
\lambda_2-1,\lambda_3,\dots,\lambda_{m-2},
\lambda_{m-1})
& \text{if $i=1$},
\\
(\lambda_1,
\lambda_2,\dots,\lambda_{m-3},\lambda_{m-2}-1,
\lambda_{m-1}+2)
& \text{if $i=m-1$},
\\
(\lambda_1,
\dots,
\lambda_{i-1}-1,
\lambda_i+2,
\lambda_{i+1}-1,
\dots,
\lambda_{m-1})
& \text{otherwise.}
\end{cases}
\end{gather*}

$U_q(\slm)$ can be endowed with the structure of a~Hopf algebra,
with coproduct $\Delta\colon U_q(\slm)\mapsto U_q(\slm)\otimes U_q(\slm)$ given on Chevalley generators by
\begin{gather*}
\Delta(E_i) = 1\otimes E_i + E_i\otimes K_i,
\qquad
\Delta(F_i) = K_i^{-1}\otimes F_i + F_i\otimes 1,
\qquad
\Delta(K_i^{\pm 1})= K_i^{\pm1} \otimes K_i^{\pm 1}.
\end{gather*}

Def\/ine ${\bigwedge}_q(\C^r)$ as the algebra generated by $r$ variables
\begin{gather*}
{\bigwedge}_q(\C^r)=\C\big[q,q^{-1}\big]\langle X_1,\dots,X_r\rangle /\big(X_i^2,\ X_iX_j+q^{-1}X_jX_i\ \text{for} \ i<j\big).
\end{gather*}
This algebra can be given a~$U_q(\mathfrak{sl}_r)$ action,
extending the natural representation\footnote{Actually,
we choose here a~non-standard form (dual) for the natural representation in order to obtain the same conventions as in~\cite{LQR1}.}.
More precisely:
\begin{gather*}
E_iX_i=X_{i+1},
\qquad
E_iX_j=0
\qquad
\text{if}\quad j\neq i,
\\
F_i X_{i+1}=X_i,\qquad F_iX_j=0\quad \text{if}\qquad j\neq i+1,
\\
K_iX_i=q^{-1}X_i, \qquad K_iX_{i+1}=qX_{i+1},
\qquad K_iX_j=X_j\quad \text{otherwise}.
\end{gather*}

We now consider ${\bigwedge}_q(\C^n\otimes \C^m)$,
where,
following~\cite{CKM},
the generating variables can be deno\-ted~$z_{ij}$ with $1\leq i\leq n$,
$1\leq j\leq m$,
subject to skew-commutation relations.

There are two isomorphisms:
\begin{gather*}
{\bigwedge}_q(\C^n)^{\otimes m}\leftarrow {\bigwedge}_q(\C^n\otimes \C^m) \rightarrow {\bigwedge}_q(\C^m)^{\otimes n}.
\end{gather*}

We can thus endow this module with actions of $U_q(\sln)$ and $U_q(\slm)$,
which Cautis,
Kamnitzer and Licata have proved to commute,
calling this \emph{quantum skew Howe duality}.
Furthermore,
$U_q(\sln)$ and $U_q(\slm)$ form a~Howe pair,
which is a~key argument in~\cite{CKM}.\footnote{Note that proving the commutation is an easy computation,
while it is much harder to prove that both algebras are each other commutant.} The actions of the two quantum groups can be deduced,
for $U_q(\slm)$ for example,
from the expressions on the variables~$z_{ij}$:
\begin{gather*}
E_j z_{ij}=z_{i,j+1},
\qquad
E_kz_{ij}=0\quad \text{if}\quad k\neq j,
\\
F_jz_{i,j+1}=z_{i,j},
\qquad
F_kz_{ij}=0
\quad
\text{if}\quad k\neq j+1,
\\
K_jz_{ij}=q^{-1}z_{ij},
\qquad
K_jz_{i,j+1}=qz_{i,j+1},
\qquad
K_kz_{ij}=z_{ij}
\quad
\text{otherwise}.
\end{gather*}

We can assign degree one to each generating variable.
Given an integer $N$,
the subspace of degree $N$ decomposes as an $\sln$-representation as follows:
\begin{gather*}
{\bigwedge}_q^{N}(\C^n\otimes \C^m)=\bigoplus_{a_1+\cdots +a_m=N}{\bigwedge}_q^{a_1}(\C^n)\otimes \cdots \otimes{\bigwedge}_q^{a_m}(\C^n).
\end{gather*}
Each direct summand is an $m$-fold tensor product of minuscule $\sln$-representations,
but is not stable under the action of $U_q(\slm)$.
However,
it appears (tracking it from the explicit def\/inition of the actions)
that each subspace $\bigwedge_q^{a_1}(\C^n)\otimes \cdots \otimes \bigwedge_q^{a_m}(\C^n)$ is a~$U_q(\slm)$ weight space
of weight $(a_2-a_1,a_3-a_2,\dots,a_m-a_{m-1})$
In particular,
if $m=2$,
the subspaces are of the form $\bigwedge_q^k(\C^n)\otimes \bigwedge_q^l(\C^n)$,
which are both $U_q(\sln)$-modules
and $U_q(\slnn{2})$ weight spaces of weight $l-k$.
The action of $U_q(\slm)$ can be explicitly tracked (see Table (\ref{sl2Action}) for the case where $m=n=N=2$),
and we see that
$E_i\colon\bigwedge_q^{a_1}(\C^n)\otimes\cdots\otimes\bigwedge_q^{a_i}(\C^n)\otimes\bigwedge_q^{a_{i+1}}(\C^n)\otimes \cdots
\otimes \bigwedge_q^{a_m}(\C^n) \mapsto \bigwedge_q^{a_1}(\C^n)\otimes \cdots \otimes \bigwedge_q^{a_i-1}(\C^n)
\otimes \bigwedge_q^{a_{i+1}+1}(\C^n)\otimes \cdots \otimes \bigwedge_q^{a_m}(\C^n)$.

This $U_q(\slm)$ action can be depicted by some particular diagrams called \emph{ladders}.
To a~direct summand ${\bigwedge}_q^{a_1}(\C^n)\otimes \cdots \otimes {\bigwedge}_q^{a_m}(\C^n)$ we assign a~sequence $(a_1,
\dots,a_m)$ depicted as weighted upward strands
\begin{gather*}
(a_1,\dots,a_m) \mapsto
\xy
(0,0)*{
\begin{tikzpicture} [scale=.75]
\draw [very thick, directed=1] (0,0) -- (0,1);
\draw [very thick, directed=1] (2,0) -- (2,1);
\draw [very thick, directed=1] (3,0) -- (3,1);
\draw [very thick, directed=1] (5,0) -- (5,1);
\node at (0,-.3) {\tiny $a_1$};
\node at (1,.5) {$\cdots$};
\node at (2,-.3) {\tiny $a_{i}$};
\node at (3,-.3) {\tiny $a_{i+1} $};
\node at (4,.5) {$\cdots$};
\node at (5,-.3) {\tiny $a_{m}$};
\node at (0,1.3) {};
\end{tikzpicture}};
\endxy
\end{gather*}

Strands labeled by zero will be erased,
and we will sometimes depict the $n$-labeled strands doubled.

We represent the action of $E_i$ and $F_i$ as follows
\begin{gather}
E_i \mapsto
\xy
(0,0)*{
\begin{tikzpicture} [scale=.75]
\draw [very thick, directed=1] (-2,0) -- (-2,3);
\draw [very thick, directed=.55] (0,0) -- (0,1);
\draw[very thick, directed=1] (0,1) -- (0,3);
\draw [very thick, directed=.55] (1,0) -- (1,2);
\draw[very thick, directed=1] (1,2) -- (1,3);
\draw [very thick, directed=.55] (0,1) -- (1,2);
\draw [very thick, directed=1] (3,0) -- (3,3);
\node at (-2,-.3) {\tiny $a_1$};
\node at (-1,1.5) {$\cdots$};
\node at (0,-.3) {\tiny $a_i$};
\node at (1,-.3) {\tiny $a_{i+1}$};
\node at (-.3,3.3) {\tiny $a_{i} - 1$};
\node at (1.3,3.3) {\tiny $a_{i+1}+1$};
\node at (2,1.5) {$\cdots$};
\node at (3,-.3) {\tiny $a_m$};
\end{tikzpicture}};
\endxy
\qquad \text{and}
\qquad
F_i \mapsto
\xy
(0,0)*{
\begin{tikzpicture} [scale=.75]
\draw [very thick, directed=1] (-2,0) -- (-2,3);
\draw [very thick, directed=.55] (0,0) -- (0,2);
\draw[very thick, directed=1] (0,2) -- (0,3);
\draw [very thick, directed=.55] (1,0) -- (1,1);
\draw[very thick, directed=1] (1,1) -- (1,3);
\draw [very thick, directed=.55] (1,1) -- (0,2);
\draw [very thick, directed=1] (3,0) -- (3,3);
\node at (-2,-.3) {\tiny $a_1$};
\node at (-1,1.5) {$\cdots$};
\node at (0,-.3) {\tiny $a_i$};
\node at (1,-.3) {\tiny $a_{i+1}$};
\node at (-.3,3.3) {\tiny $a_{i} + 1$};
\node at (1.3,3.3) {\tiny $a_{i+1} - 1$};
\node at (2,1.5) {$\cdots$};
\node at (3,-.3) {\tiny $a_m$};
\end{tikzpicture}};
\endxy
\label{Phi_Def}
\end{gather}
where we only depicted the strands $1$,
$i$,
$i+1$ and $m$: straight strands with indices $a_j$,
$j\neq 1$,
$i$,
$i+1$,
$m$ have to be added in place of the dots.
The diagrams are to be read from bottom to top.

We can def\/ine the notion of ladder as any morphism obtained by composition of identities and elementary morphisms
given by the images of $E_i$ and $F_i$.

The above diagrams have an interpretation in terms of webs.
Recall that $\sln$ webs are trivalent oriented graphs with edges indexed by integers $1,\dots,n$,
so that at each vertex,
the sum of outgoing labels equals the sum of ingoing labels (in the literature,
the $n$-strands are usually erased or only kept as \emph{tags} on the other strands,
which we will not do here).
These graphs are considered modulo some local relations (see for example~\cite{CKM},
or Def\/inition~\ref{def:slnWebs},
for a~more precise description).
The webs may be understood as $\sln$ analogues of the skein module in the $\slnn{2}$ case (see also~\cite{MOY} or~\cite{MS2} for more details).

An interesting fact is that all of the web relations in the spider category can be recovered from the relations in $U_q(\slm)$
via its action on webs\footnote{This statement holds in the case where $n$-strands are only kept as tags.
We will therefore add some relations on the $n$-strands later.}.
We refer for this and for a~complete description of the spider category to~\cite{CKM}.

We give below an example of the translation process,
which gives a~ladder whose closure is the web depicted in the introduction (with $k=l=1$) with $n=3$,
$m=3$ and $N=3$:
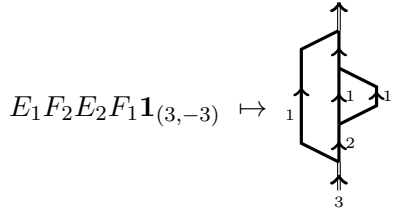
\begin{gather*}
E_1F_2E_2F_1\onell{(3,-3)} \  \mapsto \xy
(0,0)*{
\begin{tikzpicture} [scale=.5]
\draw[directed=.55,
double] (1,0) -- (1,.75);
\draw [directed=.55,
very thick] (1,.75) -- (0,1.25) -- (0,3.75) -- (1,4.25);
\draw [directed=.55,
very thick] (1,.75) -- (1,1.75);
\draw [directed=.55,
very thick] (1,1.75) -- (1,3.25);
\draw [directed=.55,
very thick] (1,1.75) -- (2,2.25)-- (2,2.75) -- (1,3.25);
\draw [directed=.55,
very thick] (1,3.25) -- (1,4.25);
\draw[directed=1,
double] (1,4.25) -- (1,5);
\node at (1,-.3) {\tiny $3$};
\node at (-.3,2) {\tiny $1$};
\node at (1.3,1.25) {\tiny $2$};
\node at (2.3,2.5) {\tiny $1$};
\node at (1.3,2.5) {\tiny $1$};
\end{tikzpicture}};
\endxy
\end{gather*}

In the case where $m=n=N=2$,
the $U_q(\slm)$ action can be explicitly given:
\begin{gather}\label{sl2Action}
\begin{tabular}{|c|c|c|c|}
\hline
summand & generator & image under $E$ & image under $F$
\\
\hline \hline
${\bigwedge}_q^2\otimes {\bigwedge}_q^0\vphantom{{\bigwedge_{q_q}}^{0^0}}$
        & $z_{11}\otimes z_{21}$ & $z_{11}\otimes z_{22}+q^{-1}z_{12}\otimes z_{21}$ & 0
\\
\hline \hline
\multirow{4}{*}{${\bigwedge}_q^1\otimes {\bigwedge}_q^1$} & $z_{11}\otimes z_{12}$ & $0$ & $0$
\\
\cline{2-4}
& $z_{21}\otimes z_{22}$ & $0$ & $0$
\\
\cline{2-4}
& $z_{11}\otimes z_{22}$ & $qz_{12}\otimes z_{22}$ & $qz_{11}\otimes z_{21} $
\\
\cline{2-4}
& $z_{12}\otimes z_{21}$ & $z_{12}\otimes z_{22}$ & $z_{11}\otimes z_{21}$
\\
\hline \hline
${\bigwedge}_q^0\otimes {\bigwedge}_q^2\vphantom{{\bigwedge_{q_q}}^{0^0}}$
        & $z_{12}\otimes z_{22}$ & $0$ & $q^{-1}z_{12}\otimes z_{21}+ z_{11}\otimes z_{22}$
\\
\hline
\end{tabular}
\end{gather}

The previous table corresponds in the diagrammatic world to the next situation:
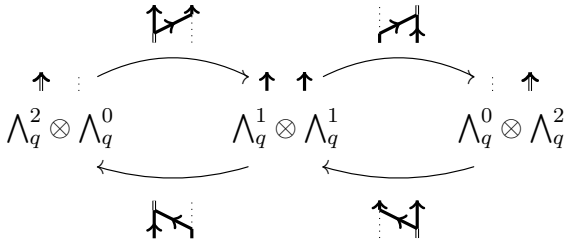
\begin{gather*}
\xy
(0,0)*{
\begin{tikzpicture} [scale=.5,decoration={markings,
mark=at position 1 with {\arrow{>}}; }]
\draw [double,
postaction={decorate}] (0,0) -- (0,.5);
\draw [dotted] (1,0) -- (1,.5);
\node at (.5,-1) {${\bigwedge}_q^2\otimes {\bigwedge}_q^0$};
\end{tikzpicture}};
(30,0)*{
\begin{tikzpicture} [scale=.5,decoration={markings,
mark=at position 1 with {\arrow{>}}; }]
\draw [very thick,
postaction={decorate}] (0,0) -- (0,.5);
\draw [very thick,
postaction={decorate}] (1,0) -- (1,.5);
\node at (.5,-1) {${\bigwedge}_q^1\otimes {\bigwedge}_q^1$};
\end{tikzpicture}};
(60,0)*{
\begin{tikzpicture} [scale=.5,decoration={markings,
mark=at position 1 with {\arrow{>}}; }]
\draw [dotted] (0,0) -- (0,.5);
\draw [double,
postaction={decorate}] (1,0) -- (1,.5);
\node at (.5,-1) {${\bigwedge}_q^0\otimes {\bigwedge}_q^2$};
\end{tikzpicture}};
(5,5); (25,5) **\crv{(15,10)} *\dir{>};
(35,5); (55,5) **\crv{(45,10)} *\dir{>};
(55,-7); (35,-7) **\crv{(45,-10)} *\dir{>};
(25,-7); (5,-7) **\crv{(15,-10)} *\dir{>};
(15,12)*{
\begin{tikzpicture} [scale=.5,]
\draw [double] (0,0) -- (0,.25);
\draw [dotted] (1,0) -- (1,.75);
\draw [very thick, directed=.55] (0,.25) -- (1,.75);
\draw [very thick, directed=1] (0,.25) -- (0,1);
\draw [very thick, directed=1] (1,.75) -- (1,1);
\end{tikzpicture}};
(45,12)*{
\begin{tikzpicture} [scale=.5]
\draw [very thick] (0,0) -- (0,.25);
\draw [very thick, directed=.55] (1,0) -- (1,.75);
\draw [very thick, directed=.55] (0,.25) -- (1,.75);
\draw [dotted] (0,.25) -- (0,1);
\draw [double] (1,.75) -- (1,1);
\end{tikzpicture}};
(45,-14)*{
\begin{tikzpicture} [scale=.5]
\draw [dotted] (0,0) -- (0,.75);
\draw [double] (1,0) -- (1,.25);
\draw [very thick, directed=.55] (1,.25) -- (0,.75);
\draw [very thick, directed=1] (0,.75) -- (0,1);
\draw [very thick, directed=1] (1,.25) -- (1,1);
\end{tikzpicture}};
(15,-14)*{
\begin{tikzpicture} [scale=.5]
\draw [very thick, directed=.55] (0,0) -- (0,.75);
\draw [very thick] (1,0) -- (1,.25);
\draw [very thick, directed=.55] (1,.25) -- (0,.75);
\draw [double] (0,.75) -- (0,1);
\draw [dotted] (1,.25) -- (1,1);
\end{tikzpicture}};
\endxy
\end{gather*}
where in the above pictures,
$0$-strands are depicted dotted and $2$-strands are doubled.

\subsubsection{Quantum Weyl group action}

The action of the Weyl group $\mathfrak{S}_m$ of $\slm$ on the weights $q$-deforms to give rise to a~braid group action
on representations of $U_q(\slm)$.
This phenomenon is referred to as the \emph{quantum Weyl group action} (see~\cite{CKL,KamnTin,Lus4}).

Generators of the braid group action are elements of the completion $\tilde{U_q(\slm)}$ of $U_q(\slm)$.
This ring is def\/ined (see~\cite{KamnTin} for example) as a~quotient of the ring of series $\sum\limits_{k=1}^{\infty} X_k$ of elements of $U_q(\slm)$,
acting on each irreducible representation $V(\lambda)$ of highest weight $\lambda$ by zero but for f\/initely many terms $X_k$.
We then consider the quotient of this ring by the two-sided ideal of elements acting by zero on all~$V(\lambda)$.

Following~\cite{Lus4},
to $s_i$ the elementary transposition corresponding to the root $\alpha_i$,
we associate the map $\Ts \in \tilde{U_q(\slm)}$:
\begin{gather*}%\label{qWeylaction_1}
\Ts\onel := \sum\limits_{a-b+c=-\lambda_i} (-1)^b q^{-ac+b} E_{i}^{(a)}F_{i}^{(b)}E_i^{(c)}\onel.
\end{gather*}

With this def\/inition,
$\Ts$ gives an endomorphism of any f\/inite-dimensional representation.
Note that if $v$ is a~weight vector of weight $\lambda$,
$T_i(v)$ is a~weight vector of weight $s_i(\lambda)$.

Taking $m=2$ for simplicity,
we have $T^{''\pm}\in \tilde{U_q(\slnn{2})}$,
acting on ${\bigwedge}_q^N(\C^n\otimes \C^2)$.
This stabilizes the whole representation,
and gives a~morphism of $U_q(\sln)$ representations,
from ${\bigwedge}_q^k(\C^n)\otimes {\bigwedge}_q^l(\C^n)$ to ${\bigwedge}_q^l(\C^n)\otimes {\bigwedge}_q^k(\C^n)$.
It is shown in~\cite{CKL} that this $U_q(\sln)$ endomorphism recovers the braiding.
This is the starting point of a~reinterpretation of Reshetikhin--Turaev invariants in terms of skew-Howe duality~\cite{Cautis},
which admits natural categorif\/ications~\cite{Cautis,
LQR1}.

The name \emph{quantum Weyl group} is used by dif\/ferent authors with slightly dif\/ferent signif\/ications.
The f\/irst one,
where we use the notation $\Ts$,
consists in considering morphisms of representations,
acting on the category of f\/inite-dimensional modules.
We can also use it to build automorphisms of the quantum group itself,
by conjugation.
Following~\cite{KamnTin},
we denote the latter $C_{\Ts}\;\colon X\mapsto \Ts X\Ts^{-1}$.
We will use both versions in this paper.
We will need some results concerning the behavior of these elements for later use.

For $w=s_{i_1}\cdots s_{i_n}$ element of the Weyl group written in reduced form,
where $s_i$ are simple ref\/lections,
we def\/ine $\Tsw{w}=\Tsw{i_1}\cdots \Tsw{i_k}$.

\begin{Proposition}[\protect{\cite[Theorem~8.1.2]{ChariPressley},  \cite[Section~37.1.3]{Lus4}, \cite{KamnTin}}] \label{propCTs}
\begin{gather*}
C_{\Ts}(E_i\onel)=-q^{-\lambda_i}F_i\onenn{s_i(\lambda)},
\qquad
C_{\Ts}(\onel F_i)=-q^{\lambda_i}\onenn{s_i(\lambda)}E_i.
\end{gather*}
For $w\in W$ such that $w(\alpha_i)=\alpha_j$,
$C_{\Tsw{w}}(E_i)=E_j$.
\end{Proposition}

Other intertwiners,
def\/ined in~\cite{Lus4},
may also be of interest:
\begin{gather*}%\label{qWeylaction_2}
\Tp\onel := \sum\limits_{a-b+c=\lambda_i} (-1)^b q^{-ac+b} F_{i}^{(a)}E_{i}^{(b)}F_i^{(c)}\onel.
\end{gather*}

We have an analogue of Proposition~\ref{propCTs}:
\begin{Proposition} \label{propCTp}
\begin{gather*}
C_{\Tp}(\onel E_i)=-q^{-\lambda_i}\onenn{s_i(\lambda)}F_i,
\qquad
C_{\Tp}(F_i\onel)=-q^{\lambda_i}E_i\onenn{s_i(\lambda)}.
\end{gather*}
For $w\in W$ such that $w(\alpha_i)=\alpha_j$,
$C_{\Tpw{w}}(E_i)=E_j$.
\end{Proposition}

The relation between the actions of both def\/initions is given by:
\begin{Proposition}[\protect{\cite[Sections 5.2.3,
37.1.2]{Lus4}}] \label{PropTsTp}
$\Ts$ and $\overline{\Tp}$,
where the bar corresponds to chan\-ging~$q$ to~$q^{-1}$,
are inverse of each other.
$\Ts\onel$ and $(-1)^{\lambda_i}q^{\lambda_i}\Tp\onel$ act the same way on any integrable module.
\end{Proposition}

\subsubsection{Skew Howe duality and quantum invariants for knots}

As we have seen,
the skew-Howe duality process gives us dif\/ferent pieces of the Jones (or Reshetikhin--Turaev) invariants:
\begin{itemize}\itemsep=0pt
\item minuscule representations $\bigwedge_q^k(\C^n)$ of $U_q(\sln)$.
This means that we are looking at knot invariants where we decorate the strands with minuscule representations.
In particular,
this does not deal with the colored Jones polynomial or its $\sln$ generalizations,
where the strands of the knot can be decorated with any f\/inite-dimensional representations.
Paths using Jones-Wenzl projectors,
and their categorif\/ications in the categorif\/ied case,
can be given to relate the general invariants to the ones we study here~\cite{CoopKrush,FSS,Rose,SS}.
\item elementary morphisms between tensor products of these representations,
given as images of $E_i$ and $F_i\in U_q(\slm)$.
These morphisms involve minuscule representations,
but do not directly deal with duals,
which in the language of knots means that we are looking at upward tangles (or their generalization for webs).
The bridge with general knots or links is established in our case in~\cite{CKM} (see also~\cite{MS2}).
\item braiding between minuscule representations,
understood in terms of the quantum Weyl group action of $U_q(\slm)$.
Again,
this is given in the framework of ladders,
and relaxing this structure will be one of the goals of the next section.
\end{itemize}

\subsection{Skein modules}

\subsubsection{Braidings for skein modules}

Let us now turn toward knots,
or ladder analogues of them.
The previous diagrammatic process gives us an algebraic interpretation of ladder webs,
as well as a~def\/inition of the braiding for the tensor product of two minuscule representations.
This braiding corresponds in the diagrammatic world to a~crossing between two adjacent strands in a~ladder,
the explicit formulas for $\Tp$ or $\Ts$ giving a~way to smooth it and replace it by a~sum of ladders without crossing.

We start by def\/ining more precisely the notion of skein module,
before relating it to the previous analysis.
By skein module,
we usually refer here both to the module itself and to the Kauf\/fman bracket def\/ining a~map from web-tangles to the module.
The def\/inition below is adapted from~\cite{CKM}.

\begin{Definition} \label{def:slnWebs} Let $n{\bf Web}$,
the $\sln$ web skein module,
be the $\Z[q,q^{-1}]$-module generated by webs (planar oriented trivalent graphs with preserved f\/low),
possibly with boundary,
up to isotopy and the following web relations:
\begin{gather}\label{webrel1}
\xy
(0,0)*{
\begin{tikzpicture} [scale=.5]
\draw [very thick, directed=.55] (0,0) -- (0,1);
\draw [very thick, directed=.55] (0,1) .. controls (-.5,2) .. (0,3);
\draw [very thick,directed=.55] (0,1) .. controls (.5,2) .. (0,3);
\draw [very thick,directed=1] (0,3) -- (0,4);
\node at (-1,.3) {\tiny $k+l$};
\node at (-.9,2) {\tiny $k$};
\node at (.9,2) {\tiny $l$};
\node at (-1,3.3) {\tiny $k+l$};
\end{tikzpicture}};
\endxy
=
{k+l \brack l}_q \xy
(0,0)*{
\begin{tikzpicture} [scale=.5]
\draw [very thick, directed=1] (0,0) -- (0,4);
\node at (-1,2) {\tiny $k+l$};
\end{tikzpicture}};
\endxy
\;,\qquad
\xy
(0,0)*{
\begin{tikzpicture} [scale=.5]
\draw [very thick, directed=.55] (0,0) -- (0,1);
\draw [very thick,
rdirected=.45] (0,1) .. controls (-.5,2) .. (0,3);
\draw [very thick,directed=.55] (0,1) .. controls (.5,2) .. (0,3);
\draw [very thick,directed=1] (0,3) -- (0,4);
\node at (-1,.3) {\tiny $k$};
\node at (-.9,2) {\tiny $l$};
\node at (1.6,2) {\tiny $k+l$};
\node at (-1,3.3) {\tiny $k$};
\end{tikzpicture}};
\endxy
=
{n-k \brack l}_q
\xy
(0,0)*{
\begin{tikzpicture} [scale=.5]
\draw [very thick, directed=1] (0,0) -- (0,4);
\node at (-.4,2) {\tiny $k$};
\end{tikzpicture}};
\endxy
\\
\label{webrel2}
\xy
(0,0)*{
\begin{tikzpicture} [scale=.5]
\draw [very thick, directed=.55] (0,0) -- (.5,1);
\draw [very thick, directed=.55] (1,0) -- (.5,1);
\draw [very thick, directed=.55] (.5,1) -- (1,2);
\draw [very thick, directed=.55] (2,0) -- (1,2);
\draw [very thick, directed=1] (1,2) -- (1,3);
\node at (0,-.3) {\tiny $k$};
\node at (1,-.3) {\tiny $l$};
\node at (2,-.3) {\tiny $m$};
\node at (1,3.4) {\tiny $k+l+m$};
\end{tikzpicture}
};
\endxy
=
\xy
(0,0)*{
\begin{tikzpicture} [scale=.5]
\draw [very thick, directed=.55] (0,0) -- (1,2);
\draw [very thick, directed=.55] (1,0) -- (1.5,1);
\draw [very thick, directed=.55] (2,0) -- (1.5,1);
\draw [very thick, directed=.55] (1.5,1) -- (1,2);
\draw [very thick, directed=1] (1,2) -- (1,3);
\node at (0,-.3) {\tiny $k$};
\node at (1,-.3) {\tiny $l$};
\node at (2,-.3) {\tiny $m$};
\node at (1,3.4) {\tiny $k+l+m$};
\end{tikzpicture}
};
\endxy
,\qquad
\xy
(0,0)*{
\begin{tikzpicture} [scale=.5]
\draw [very thick,
directed =.55] (0,0) -- (0,.75);
\draw [very thick,
directed =.55] (0,.75) -- (0,1.75);
\draw [very thick,
directed =1] (0,1.75) -- (0,3);
\draw [very thick,
directed =.55] (1,0) -- (1,1.25);
\draw [very thick,
directed =.55] (1,1.25) -- (1,2.25);
\draw [very thick,
directed =1] (1,2.25) -- (1,3);
\draw [very thick, directed=.55] (0,.75) -- (1,1.25);
\draw [very thick, directed=.55] (0,1.75) -- (1,2.25);
\node at (0,-.3) {\tiny $k$};
\node at (1,-.3) {\tiny $l$};
\node at (.5,1.3) {\tiny $s$};
\node at (.5,2.3) {\tiny $r$};
\node at (0,
3.3) {};
\node at (1,3.3) {};
\end{tikzpicture}
};
\endxy
=
{r+s \brack r}_q
\xy
(0,0)*{
\begin{tikzpicture} [scale=.5]
\draw [very thick,
directed =.55] (0,0) -- (0,1.25);
\draw [very thick,
directed =1] (0,1.25) -- (0,3);
\draw [very thick,
directed =.55] (1,0) -- (1,1.75);
\draw [very thick,
directed =1] (1,1.75) -- (1,3);
\draw [very thick, directed=.55] (0,1.25) -- (1,1.75);
\node at (0,-.3) {\tiny $k$};
\node at (1,-.3) {\tiny $l$};
\node [rotate=90] at (.5,2.2) {\tiny $r+s$};
\node at (0,
3.3) {};
\node at (1,3.3) {};
\end{tikzpicture}
};
\endxy
\\
\label{webrel3}
\xy
(0,0)*{
\begin{tikzpicture} [scale=.5]
\draw [very thick,
directed =.55] (0,0) -- (0,.75);
\draw [very thick,
directed =.55] (0,.75) -- (0,2.25);
\draw [very thick,
directed =1] (0,2.25) -- (0,3);
\draw [very thick,
directed =.55] (1,0) -- (1,1.25);
\draw [very thick,
directed =.55] (1,1.25) -- (1,1.75);
\draw [very thick,
directed =1] (1,1.75) -- (1,3);
\draw [very thick, directed=.55] (0,.75) -- (1,1.25);
\draw [very thick, directed=.55] (1,1.75) -- (0,2.25);
\node at (0,-.3) {\tiny $k$};
\node at (1,-.3) {\tiny $l$};
\node at (.5,1.3) {\tiny $s$};
\node at (.5,2.3) {\tiny $r$};
\node at (0,
3.3) {};
\node at (1,3.3) {};
\end{tikzpicture}
};
\endxy
=
\sum\limits_t {k-l+r-s \brack t}_q
\xy
(0,0)*{
\begin{tikzpicture} [scale=.5]
\draw [very thick,
directed =.55] (0,0) -- (0,1.25);
\draw [very thick,
directed =.55] (0,1.25) -- (0,1.75);
\draw [very thick,
directed =1] (0,1.75) -- (0,3);
\draw [very thick,
directed =.55] (1,0) -- (1,.75);
\draw [very thick,
directed =.55] (1,.75) -- (1,2.25);
\draw [very thick,
directed =1] (1,2.25) -- (1,3);
\draw [very thick, directed=.55] (1,.75) -- (0,1.25);
\draw [very thick, directed=.55] (0,1.75) -- (1,2.25);
\node at (0,-.3) {\tiny $k$};
\node at (1,-.3) {\tiny $l$};
\node [rotate=90] at (.5,.2) {\tiny $s-t$};
\node [rotate=90] at (.5,2.8) {\tiny $r-t$};
\node at (0,
3.3) {};
\node at (1,3.3) {};
\end{tikzpicture}
};
\endxy
\end{gather}

All equations come with the ones obtained by mirror image and arrow reversion.
Recall that $[k]= \frac{q^k-q^{-k}}{q-q^{-1}}$,
$[k]!=[k][k-1]\cdots[1]$ and ${p \brack k}_q=\frac{[p]!}{[k]![p-k]!}$.

\end{Definition}

The skein module described above can be given the structure of a~category,
with objects given by oriented points with labels on a~horizontal line (the boundary of strands),
and morphisms the webs joining the dots on two such parallel lines.
A~subcategory is of particular interest for the skew-Howe interpretation,
and turns out to essentially represent all the information we need.
Assuming a~value of $N$ is f\/ixed,
let us call $\Phi \colon \U_q(\slm) \mapsto n{\bf Web}$ the map described in equation~\eqref{Phi_Def}.

\begin{Definition}
Def\/ine $n{\bf Web}^+_m$ to be the image category $\Phi(\U_q(\slm))$,
with objects,
sequences $(a_1,\dots,a_m)$ ($0\leq a_i\leq n$) labeling points on a~horizontal line,
together with a~zero object,
and morphisms,
$\sln$ webs between such sequences (in the sense of Def\/inition~\ref{def:slnWebs}),
that are composition of the images of~$E_i$ and~$F_i$,
as depicted in equation~\eqref{Phi_Def}.
\end{Definition}

We sometimes refer to such webs as upward webs,
or ladders.
These ladder webs have their boundary split in two parts,
with all strands oriented inside for the bottom part,
and outside for the upper part.
Although general webs are more general than this particular situation,
it is shown in~\cite{CKM} that they can be related to the particular class of webs obtained from ladders using the tool of pivotal categories.

The use of tags makes the situation somewhat simpler (but harder to f\/it in a~skein module formulation!),
but the next relations (and the ones obtained by symmetry on the next ones) are particular realizations of the ones given in~\cite{CKM} in the case
where we keep the $n$-labeled strands (and are special cases of Def\/inition~\ref{def:slnWebs}):
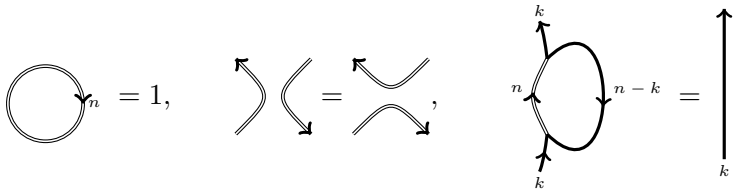
\begin{gather*}
\xy
(0,0)*{
\begin{tikzpicture} [scale=.5,decoration={markings,
mark=at position 0.5 with {\arrow{>}}; }]
\draw [double,
postaction={decorate}] (-1,0) arc (180:-180:1);
\node at (1.3,0) {\tiny $n$};
\end{tikzpicture}};
\endxy
= 1
,\qquad
\xy
(0,0)*{
\begin{tikzpicture} [scale=.5]
\draw [double, directed=1] (0,0) .. controls (1,1) .. (0,2);
\draw [double, directed=1] (2,2) .. controls (1,1) .. (2,0);
\end{tikzpicture}};
\endxy
=
\xy
(0,0)*{
\begin{tikzpicture} [scale=.5]
\draw [double, directed=1] (0,0) .. controls (1,1) .. (2,0);
\draw [double, directed=1] (2,2) .. controls (1,1) .. (0,2);
\end{tikzpicture}};
\endxy
,
\qquad
\xy
(0,0)*{
\begin{tikzpicture} [scale=.5]
\node at (-.2,-.3) {\tiny $k$};
\draw [very thick, directed=.55] (-.2,0) .. controls (-.1,.2) and (0,.9) .. (0,1);
\draw [double, directed=.55] (0,1) .. controls (-.5,2) .. (0,3);
\node at (-.8,2.2) {\tiny $n$};
\draw [very thick, directed=.55] (0,3) .. controls (2,5) and (2,-1) .. (0,1);
\node at (2.4,2.2) {\tiny $n-k$};
\draw [very thick, directed=1] (0,3)..
controls (0,3.1) and (-.1,3.8) .. (-.2,4);
\node at (-.2,4.3) {\tiny $k$};
\end{tikzpicture}};
\endxy
=\xy
(0,0)*{
\begin{tikzpicture} [scale=.5]
\draw [very thick, directed=1] (0,0) -- (0,4);
\node at (0,-.3) {\tiny $k$};
\end{tikzpicture}};
\endxy
\end{gather*}

In the above pictures,
the $n$-th strands are depicted doubled to emphasize their particular role.

The skein module $n{\bf Web}$ is a~natural target for maps from (equivalence classes of) diagrams of knotted webs.
We call knotted webs,
or web-tangles,
the natural generalization of knots and tangles to webs.
For example,
knotted webs are isotopy classes of embeddings into $\R^3$ of closed webs,
and produce diagrams of knotted webs as generic projections onto a~plane.
Web-tangles are the natural generalization allowing boundaries.

Recall from~\cite{Kauf} (see also~\cite[Theorem 2]{CarterFoams}) the relations that generalize Reidemeister moves (in a~framed version,
where a~numbered circle on a~strand stands for twists): any two diagrams representing the same web-tangle are related by a~sequence
of moves of the following kind:
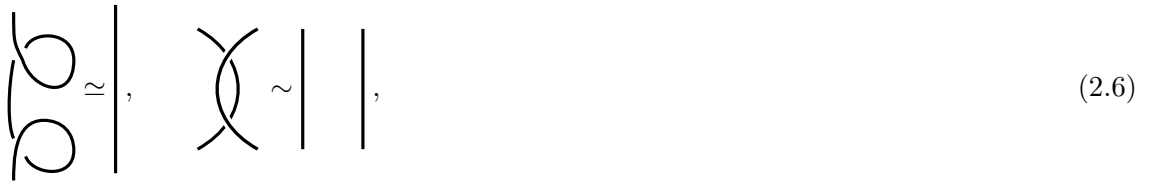
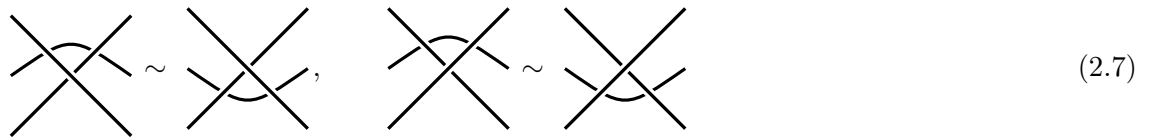
\begin{gather}
\xy
(0,0)*{
\begin{tikzpicture} [scale=.8,decoration={markings,
mark=at position 0.8 with {\arrow{>}}; }]
\draw [very thick] (0,0) .. controls (0,.3) and (0,1) .. (.5,1) -- (.5,1) .. controls (.8,1) and (1,.8) .. (1,.5) -- (1,.5) .. controls (1,0)
 and (.3,.1) .. (.2,.4);
\draw [very thick] (0,.7) .. controls (-.1,.9) and (-.1,1.5) .. (0,2);
\draw [very thick] (.2,2.2) .. controls (.3,2.5) and (1,2.5) .. (1,2) -- (1,2) .. controls (1,1.3)
 and (.3,1.5) .. (.15,2) -- (.15,2) .. controls (0,2.3) .. (0,2.8);
\end{tikzpicture}};
\endxy
\simeq
\xy
(0,0)*{
\begin{tikzpicture} [scale=.8,decoration={markings,
mark=at position 0.8 with {\arrow{>}}; }]
\draw [very thick] (0,0) -- (0,2.8);
\end{tikzpicture}};
\endxy
\;,
\qquad
\xy
(0,0)*{
\begin{tikzpicture} [scale=.8,decoration={markings,
mark=at position 0.8 with {\arrow{>}}; }]
\draw [very thick] (0,0) .. controls (.9,.5) and (.9,1.5) .. (0,2);
\draw [draw =white,
double=black,
very thick,
double distance=1.25pt] (1,0) .. controls (.1,.5) and (.1,1.5) .. (1,2);
\end{tikzpicture}};
\endxy
\sim
\xy
(0,0)*{
\begin{tikzpicture} [scale=.8,decoration={markings,
mark=at position 0.8 with {\arrow{>}}; }]
\draw [very thick] (0,0) -- (0,2);
\draw [very thick] (1,0) -- (1,2);
\end{tikzpicture}};
\endxy
\;,
\label{ReidUsual12}
\\
\xy
(0,0)*{
\begin{tikzpicture} [scale=.8,decoration={markings,
mark=at position 0.8 with {\arrow{>}}; }]
\draw [very thick] (0,1) .. controls (1,1.7) .. (2,1);
\draw [very thick,
draw=white,
double=black,
double distance=1.25pt] (0,0) -- (2,2);
\draw [very thick,
draw=white,
double=black,
double distance=1.25pt] (2,0) -- (0,2);
\end{tikzpicture}};
\endxy
\sim \;
\xy
(0,0)*{
\begin{tikzpicture} [scale=.8,decoration={markings,
mark=at position 0.8 with {\arrow{>}}; }]
\draw [very thick] (0,1) .. controls (1,.3) .. (2,1);
\draw [very thick,
draw=white,
double=black,
double distance=1.25pt] (0,0) -- (2,2);
\draw [very thick,
draw=white,
double=black,
double distance=1.25pt] (2,0) -- (0,2);
\end{tikzpicture}};
\endxy
,\qquad\xy
(0,0)*{
\begin{tikzpicture} [scale=.8,decoration={markings,
mark=at position 0.8 with {\arrow{>}}; }]
\draw [very thick] (0,1) .. controls (1,1.7) .. (2,1);
\draw [very thick,
draw=white,
double=black,
double distance=1.25pt] (2,0) -- (0,2);
\draw [very thick,
draw=white,
double=black,
double distance=1.25pt] (0,0) -- (2,2);
\end{tikzpicture}};
\endxy
\sim \;
\xy
(0,0)*{
\begin{tikzpicture} [scale=.8,decoration={markings,
mark=at position 0.8 with {\arrow{>}}; }]
\draw [very thick] (0,1) .. controls (1,.3) .. (2,1);
\draw [very thick,
draw=white,
double=black,
double distance=1.25pt] (2,0) -- (0,2);
\draw [very thick,
draw=white,
double=black,
double distance=1.25pt] (0,0) -- (2,2);
\end{tikzpicture}};
\endxy \label{ReidUsual2}
\\
\xy
(0,0)*{
\begin{tikzpicture} [scale=.8,decoration={markings,
mark=at position 0.8 with {\arrow{>}}; }]
\draw [very thick] (0,1) .. controls (1,1.7) .. (2,1);
\draw [very thick,
draw=white,
double=black,
double distance=1.25pt] (1,0) -- (1,1);
\draw [very thick,
draw=white,
double=black,
double distance=1.25pt] (1,1) -- (2,2);
\draw [very thick,
draw=white,
double=black,
double distance=1.25pt] (1,1) -- (0,2);
\end{tikzpicture}};
\endxy
\sim \;
\xy
(0,0)*{
\begin{tikzpicture} [scale=.8,decoration={markings,
mark=at position 0.8 with {\arrow{>}}; }]
\draw [very thick] (0,1) .. controls (1,.3) .. (2,1);
\draw [very thick,
draw=white,
double=black,
double distance=1.25pt] (1,0) -- (1,1);
\draw [very thick,
draw=white,
double=black,
double distance=1.25pt] (1,1) -- (2,2);
\draw [very thick,
draw=white,
double=black,
double distance=1.25pt] (1,1) -- (0,2);
\end{tikzpicture}};
\endxy
,\qquad\xy
(0,0)*{
\begin{tikzpicture} [scale=.8,decoration={markings,
mark=at position 0.8 with {\arrow{>}}; }]
\draw [very thick,
draw=white,
double=black,
double distance=1.25pt] (1,0) -- (1,1);
\draw [very thick,
draw=white,
double=black,
double distance=1.25pt] (1,1) -- (2,2);
\draw [very thick,
draw=white,
double=black,
double distance=1.25pt] (1,1) -- (0,2);
\draw [very thick,
draw=white,
double=black,
double distance=1.25pt] (0,1) .. controls (1,1.7) .. (2,1);
\end{tikzpicture}};
\endxy
\sim \;
\xy
(0,0)*{
\begin{tikzpicture} [scale=.8,decoration={markings,
mark=at position 0.8 with {\arrow{>}}; }]
\draw [very thick,
draw=white,
double=black,
double distance=1.25pt] (1,0) -- (1,1);
\draw [very thick,
draw=white,
double=black,
double distance=1.25pt] (1,1) -- (2,2);
\draw [very thick,
draw=white,
double=black,
double distance=1.25pt] (1,1) -- (0,2);
\draw [very thick,
draw=white,
double=black,
double distance=1.25pt] (0,1) .. controls (1,.3) .. (2,1);
\end{tikzpicture}};
\endxy
,\label{ReidWebs1}
\\
\xy
(0,0)*{
\begin{tikzpicture} [scale=.8,decoration={markings,
mark=at position 0.8 with {\arrow{>}}; }]
\draw [very thick] (0,0) -- (0,1);
\draw [very thick] (0,1) -- (-.5,2) -- (.5,3);
\draw [very thick] (0,1) -- (.5,2) -- (.2,
2.3);
\draw [very thick] (-.2,2.7) -- (-.5,3);
\draw [very thick] (-.5,3) -- (-.5,4);
\draw [very thick] (.5,3) -- (.5,4);
\end{tikzpicture}};
\endxy
\quad \sim \quad
\xy
(0,0)*{
\begin{tikzpicture} [scale=.8,decoration={markings,
mark=at position 1 with {\arrow{>}}; }]
\draw [very thick] (0,0) -- (0,1);
\draw [very thick] (0,1) -- (-.5,2);
\draw [very thick] (0,1) -- (.5,2);
\draw [very thick] (-.5,2) -- (-.5,3.5);
\draw [very thick] (.5,2) -- (.5,3.5);
\node at (0,.7) {$\circ$};
\node at (-.5,2.7) {$\circ$};
\node at (.5,2.7) {$\circ$};
\node at (.4,.7) {\tiny$\frac{1}{2}$};
\node at (-.1,2.7) {\tiny$\frac{-1}{2}$};
\node at (1,2.7) {\tiny$\frac{-1}{2}$};
\end{tikzpicture}};
\endxy
\label{ReidWebs2}
\end{gather}

\begin{Definition}
A~Kauf\/fman bracket for $\sln$ webs is a~map $\Psi$ from diagrams of $\sln$ web-tangles to $n{\bf Web}$,
def\/ined locally by replacing a~crossing by a~linear combination of smoothings,
and subject to relations \eqref{ReidUsual12}--\eqref{ReidWebs2}.
\end{Definition}

Let us now restrict to the knotted analogue of $n{\bf Web}^+_m$,
and def\/ine a~knotted ladder (or web-tangle in ladder position) to be a~vertical composition of images of $E_i\onel\in \U_q(\slm)$,
$F_i\onel \in \U_q(\slm)$ and crossings between two adjacent uprights in the ladder.
Interpreting a~crossing between the $i$-th and $(i+1)$-th strands as the quantum Weyl group action given by $\Ts$,
one obtains a~smoothing process for crossings:
\begin{gather*}
\Psi_{\Ts}\left(
\xy
(0,0)*{
\begin{tikzpicture} [scale=.5,decoration={markings,
mark=at position .7 with {\arrow{>}}; }]
\draw [very thick,
postaction={decorate}] (0,0) -- (0,.75);
\draw [very thick,
postaction={decorate}] (0,.75) -- (0,4.25);
\draw [very thick,
postaction={decorate}] (0,4.25) -- (0,5);
\draw [very thick,
postaction={decorate}] (3,0) -- (3,.75);
\draw [very thick,
postaction={decorate}] (3,.75) -- (3,4.25);
\draw [very thick,
postaction={decorate}] (3,4.25) -- (3,5);
\draw [very thick,
postaction={decorate}] (0,.75) -- (1,1.25) -- (1,2) -- (1.5,2.5);
\draw [very thick,
postaction={decorate}] (1.5,2.5) -- (2,3) -- (2,3.75) -- (3,4.25);
\draw [very thick,
postaction={decorate}] (3,.75) -- (2,1.25) -- (2,2) -- (1.6,2.4);
\draw [very thick,
postaction={decorate}] (1.4,2.6) -- (1,3) -- (1,3.75) -- (0,4.25);
\node at (0,-.3) {\tiny $2$};
\node at (3,-.3) {\tiny $2$};
\node at (-.3,2.5) {\tiny $1$};
\node at (3.3,2.5) {\tiny $1$};
\end{tikzpicture}};
\endxy
\right)
=
\xy
(0,0)*{
\begin{tikzpicture} [scale=.5,decoration={markings,
mark=at position .7 with {\arrow{>}}; }]
\draw [very thick,
postaction={decorate}] (0,0) -- (0,.75);
\draw [very thick,
postaction={decorate}] (0,.75) -- (0,4.25);
\draw [very thick,
postaction={decorate}] (0,4.25) -- (0,5);
\draw [very thick,
postaction={decorate}] (3,0) -- (3,.75);
\draw [very thick,
postaction={decorate}] (3,.75) -- (3,4.25);
\draw [very thick,
postaction={decorate}] (3,4.25) -- (3,5);
\draw [very thick,
postaction={decorate}] (0,.75) -- (1,1.25) -- (1,2) -- (1,3) -- (1,3.75) -- (0,4.25);
\draw [very thick,
postaction={decorate}] (3,.75) -- (2,1.25) -- (2,2) -- (2,3) -- (2,3.75) -- (3,4.25);
\node at (0,-.3) {\tiny $2$};
\node at (3,-.3) {\tiny $2$};
\node at (-.3,2.5) {\tiny $1$};
\node at (3.3,2.5) {\tiny $1$};
\end{tikzpicture}};
\endxy
-q
\xy
(0,0)*{
\begin{tikzpicture} [scale=.5,decoration={markings,
mark=at position .7 with {\arrow{>}}; }]
\draw [very thick,
postaction={decorate}] (0,0) -- (0,.75);
\draw [very thick,
postaction={decorate}] (0,.75) -- (0,4.25);
\draw [very thick,
postaction={decorate}] (0,4.25) -- (0,5);
\draw [very thick,
postaction={decorate}] (3,0) -- (3,.75);
\draw [very thick,
postaction={decorate}] (3,.75) -- (3,4.25);
\draw [very thick,
postaction={decorate}] (3,4.25) -- (3,5);
\draw [very thick,
postaction={decorate}] (0,.75) -- (1,1.25) -- (1,2.25);
\draw [very thick,
postaction={decorate}] (1,2.25) -- (1,2.75);
\draw [very thick,
postaction={decorate}] (1,2.75) -- (1,3.75) -- (0,4.25);
\draw [very thick,
postaction={decorate}] (3,.75) -- (2,1.25) -- (2,1.75) -- (1,2.25);
\draw [very thick,
postaction={decorate}] (1,2.75) -- (2,3.25) -- (2,3.75) -- (3,4.25);
\node at (0,-.3) {\tiny $2$};
\node at (3,-.3) {\tiny $2$};
\node at (-.3,2.5) {\tiny $1$};
\node at (3.3,2.5) {\tiny $1$};
\node at (.7,2.5) {\tiny $2$};
\end{tikzpicture}};
\endxy
\end{gather*}

Thus,
smoothing all crossings in a~ladder web-tangle,
one obtains a~formal sum of non-knotted ladders that one can see as an element of a~skein module.

To obtain a~more powerful skein module allowing less rigidity,
we want to forget the $0$-labeled strands.
Indeed,
in ladder position,
even if the $0$-labeled strands are not depicted,
one knows where they are.
If we want to start from any diagram and use the same smoothing rules as in the ladder case,
we cannot know where $0$-strands should be and we want to make sure that crossings involving $0$-labeled strands do not play any role.

The goal of this section is to obtain a~skew-Howe duality process with a~conveniently rescaled braiding,
so that the smoothing rules derived from this braiding induce a~Kauf\/fman bracket for general $\sln$ web-tangles.

Using~\cite[Proposition 5.2.2]{Lus4},
we can show that the use of the smoothing rules provided by~$\Ts$ gives the expected non-rescaled diagonal strand for the ``trivial''
positive $(k,0)$ crossing.
Similarly,
using $\Tp^{-1}$ gives the expected result for the negative $(0,k)$ crossing.

\begin{Proposition}%\label{0smoothings}
\begin{gather*}
\Psi_{\Ts}\left(
\xy
(0,0)*{
\begin{tikzpicture} [scale=.5]
\draw [very thick, directed=1] (0,0) -- (2,2);
\draw [dotted] (2,0) -- (1.2,.8);
\draw [dotted, directed=1] (.8,1.2) -- (0,2);
\node at (0,-.5) {\tiny $k$};
\node at (2,-.5) {\tiny $0$};
\end{tikzpicture}};
\endxy
\right)
=
\xy
(0,0)*{
\begin{tikzpicture} [scale=.5]
\draw [very thick,
postaction={decorate}, directed=.55, directed=1] (0,0) -- (0,.5) -- (2,1.5) -- (2,2);
\draw [dotted, directed=.55] (2,0) -- (2,1.5);
\draw [dotted, directed=1] (0,.5) -- (0,2);
\node at (0,-.5) {\tiny $k$};
\node at (2,-.5) {\tiny $0$};
\end{tikzpicture}};
\endxy
,
\qquad
\Psi_{\Tp^{-1}}\left(
\xy
(0,0)*{
\begin{tikzpicture} [scale=.5]
\draw [very thick, directed=1] (2,0) -- (0,2);
\draw [dotted] (0,0) -- (.8,.8);
\draw [dotted, directed=1] (1.2,1.2) -- (2,2);
\node at (2,-.5) {\tiny $k$};
\node at (0,-.5) {\tiny $0$};
\end{tikzpicture}};
\endxy
\right)
=
\xy
(0,0)*{
\begin{tikzpicture} [scale=.5]
\draw [very thick, directed=.55, directed=1] (2,0) -- (2,.5) -- (0,1.5) -- (0,2);
\draw [dotted] (0,0) -- (0,1.5);
\draw [dotted, directed=1] (2,.5) -- (2,2);
\node at (2,-.5) {\tiny $k$};
\node at (0,-.5) {\tiny $0$};
\end{tikzpicture}};
\endxy
\end{gather*}
\end{Proposition}

Note that if we use $\Ts$ on a~$(0,k)$ crossing and $\Tp$ on a~$(k,0)$ crossing,
the situation is dif\/ferent.
Using Proposition~\ref{PropTsTp},
we have:
\begin{gather*}
\Psi_{\Tp}\left(
\xy
(0,0)*{
\begin{tikzpicture} [scale=.5]
\draw [very thick, directed=1] (0,0) -- (2,2);
\draw [dotted] (2,0) -- (1.2,.8);
\draw [dotted, directed=1] (.8,1.2) -- (0,2);
\node at (0,-.5) {\tiny $k$};
\node at (2,-.5) {\tiny $0$};
\end{tikzpicture}};
\endxy
\right)
=
(-1)^kq^{k}\xy
(0,0)*{
\begin{tikzpicture} [scale=.5]
\draw [very thick,
postaction={decorate}, directed=.55, directed=1] (0,0) -- (0,.5) -- (2,1.5) -- (2,2);
\draw [dotted, directed=.55] (2,0) -- (2,1.5);
\draw [dotted, directed=1] (0,.5) -- (0,2);
\node at (0,-.5) {\tiny $k$};
\node at (2,-.5) {\tiny $0$};
\end{tikzpicture}};
\endxy
,
\qquad
\Psi_{\Ts^{-1}}
\left(
\xy
(0,0)*{
\begin{tikzpicture} [scale=.5]
\draw [very thick, directed=1] (2,0) -- (0,2);
\draw [dotted] (0,0) -- (.8,.8);
\draw [dotted, directed=1] (1.2,1.2) -- (2,2);
\node at (2,-.5) {\tiny $k$};
\node at (0,-.5) {\tiny $0$};
\end{tikzpicture}};
\endxy
\right)
=
(-1)^kq^{-k}\xy
(0,0)*{
\begin{tikzpicture} [scale=.5]
\draw [very thick, directed=.55, directed=1] (2,0) -- (2,.5) -- (0,1.5) -- (0,2);
\draw [dotted] (0,0) -- (0,1.5);
\draw [dotted, directed=1] (2,.5) -- (2,2);
\node at (2,-.5) {\tiny $k$};
\node at (0,-.5) {\tiny $0$};
\end{tikzpicture}};
\endxy
\end{gather*}

It appears that we cannot choose one of the two solutions and apply it in all cases,
since there would always be a~situation where a~trivial crossing would lead to a~non-trivially rescaled piece of strand in the associated skein module.
A~natural idea would be to use a~braiding mixing both def\/initions,
which may produce some gaps if we still want to have some instances of Propositions~\ref{propCTs} and~\ref{propCTp} (which will prove useful later).

In order to avoid these distortions,
we introduce additional rescalings that utilize $U_q(\glm)$ data that is naturally encoded in the representation we are looking at (namely,
the sequence $(a_1,\dots,a_m)$,
which is determined by the $\slm$ weight and the choice of an integer $N$).
Note that the following def\/initions are rather symmetric in the $\Tp$'s and $\Ts$'s
\begin{gather}
\T\onel=(-1)^{-a_{i+1}}q^{-a_{i+1}}\Ts\onel = (-1)^{-a_i}q^{-a_i}\Tp\onel,
\nonumber\\
\T^{-1}\onel=(-1)^{a_i}q^{a_i}\Ts^{-1}\onel = (-1)^{a_{i+1}}q^{a_{i+1}}\Tp^{-1}\onel.\label{DefBraiding}
\end{gather}
It is easy to see from Proposition~\ref{PropTsTp} that both def\/initions agree,
and that this def\/inition still provides a~braiding.
We can check that we still have $C_{\Tw{1}\Tw{2}}(E_1)=E_2$ as endomorphisms of a~given representation appearing in the skew Howe context.

\subsubsection[$\mathfrak{sl}_2$ case]{$\boldsymbol{\mathfrak{sl}_2}$ case}%\label{detailedsl2case}

Let us now give a~complete description of the $\slnn{2}$ case.
In~\cite{Blan},
Blanchet introduces $\slnn{2}$ webs to be oriented trivalent graphs with two kinds of edges ($1$ and $2$-labeled,
we draw the latter doubled),
and vertices having two ingoing $1$-labeled strands and one outgoing $2$-labeled one,
or one ingoing $2$-labeled strand and two outgoing $1$-labeled ones.

The $\slnn{2}$ skein module $2{\bf Web}$ is the quotient of (linear combinations of) webs with edges labeled $1$ or $2$ by the next relations:
\begin{gather}
\xy
(0,0)*{
\begin{tikzpicture} [scale=.5]
\draw [very thick] (-1,0) arc (180:-180:1);
\end{tikzpicture}};
\endxy
= [2]
,\qquad
\xy
(0,0)*{
\begin{tikzpicture} [scale=.5]
\draw [double] (-1,0) arc (180:-180:1);
\end{tikzpicture}};
\endxy = 1
,\qquad
\xy
(0,0)*{
\begin{tikzpicture} [scale=.5]
\draw [double, directed=1] (0,0) .. controls (1,1) .. (0,2);
\draw [double, directed=1] (2,2) .. controls (1,1) .. (2,0);
\end{tikzpicture}};
\endxy
=
\xy
(0,0)*{
\begin{tikzpicture} [scale=.5]
\draw [double, directed=1] (0,0) .. controls (1,1) .. (2,0);
\draw [double, directed=1] (2,2) .. controls (1,1) .. (0,2);
\end{tikzpicture}};
\endxy \label{SkeinCircles}
\\
\xy
(0,0)*{
\begin{tikzpicture} [scale=.5,decoration={markings,
mark=at position 0.5 with {\arrow{>}}; }]
\draw [double] (0,0) -- (0,1);
\draw [very thick] (0,1) .. controls (-.5,2) .. (0,3);
\draw [very thick] (0,1) .. controls (.5,2) .. (0,3);
\draw [double] (0,3) -- (0,4);
\end{tikzpicture}};
\endxy
=
[2]
\xy
(0,0)*{
\begin{tikzpicture} [scale=.5,decoration={markings,
mark=at position 0.5 with {\arrow{>}}; }]
\draw [double] (0,0) -- (0,4);
\end{tikzpicture}};
\endxy
\; ,
\qquad
\xy
(0,0)*{
\begin{tikzpicture} [scale=.5]
\draw [very thick] (0,0) -- (0,1);
\draw [double] (0,1) .. controls (-.5,2) .. (0,3);
\draw [very thick] (0,1) .. controls (2,-1) and (2,5) .. (0,3);
\draw [very thick] (0,3) -- (0,4);
\end{tikzpicture}};
\endxy
=\xy
(0,0)*{
\begin{tikzpicture} [scale=.5]
\draw [very thick] (0,0) -- (0,4);
\end{tikzpicture}};
\endxy \label{SkeinDigons}
\end{gather}

For non-oriented webs above,
the depicted relations hold for any compatible orientation.

The def\/inition of the braidings then gives the following smoothing rules:
\begin{gather*}
\Psi_{\T}\left(
\xy
(0,0)*{
\begin{tikzpicture} [scale=.5]
\draw [very thick, directed=1] (0,0) -- (2,2);
\draw [very thick] (2,0) -- (1.2,.8);
\draw [very thick, directed=1] (.8,1.2) -- (0,2);
\end{tikzpicture}};
\endxy
\right)
\; = \; -q^{-1}
\xy
(0,0)*{
\begin{tikzpicture} [scale=.5]
\draw [very thick, directed=1] (0,0) .. controls (1,1) .. (0,2);
\draw [very thick, directed=1] (2,0) .. controls (1,1) .. (2,2);
\end{tikzpicture}};
\endxy
\; +\;
\xy
(0,0)*{
\begin{tikzpicture} [scale=.5]
\draw [very thick, directed=.55] (0,0) -- (1,.7);
\draw [very thick, directed=.55] (2,0) -- (1,.7);
\draw [double, directed=.55] (1,.7) -- (1,1.3);
\draw [very thick, directed=1] (1,1.3) -- (0,2);
\draw [very thick, directed=1] (1,1.3) -- (2,2);
\end{tikzpicture}};
\endxy
,\qquad
\Psi_{\T}\left(
\xy
(0,0)*{
\begin{tikzpicture} [scale=.5]
\draw [very thick] (0,0) -- (.8,.8);
\draw [very thick, directed=1] (1.2,1.2) -- (2,2);
\draw [very thick, directed=1] (2,0) -- (0,2);
\end{tikzpicture}};
\endxy
\right)
\; = \; -q
\xy
(0,0)*{
\begin{tikzpicture} [scale=.5]
\draw [very thick, directed=1] (0,0) .. controls (1,1) .. (0,2);
\draw [very thick, directed=1] (2,0) .. controls (1,1) .. (2,2);
\end{tikzpicture}};
\endxy
+
\xy
(0,0)*{
\begin{tikzpicture} [scale=.5]
\draw [very thick, directed=.55] (0,0) -- (1,.7);
\draw [very thick, directed=.55] (2,0) -- (1,.7);
\draw [double, directed=.55] (1,.7) -- (1,1.3);
\draw [very thick, directed=1] (1,1.3) -- (0,2);
\draw [very thick, directed=1] (1,1.3) -- (2,2);
\end{tikzpicture}};
\endxy
\\
\Psi_{\T}\left(
\xy
(0,0)*{
\begin{tikzpicture} [scale=.5]
\draw [double, directed=1] (0,0) -- (2,2);
\draw [very thick] (2,0) -- (1.2,.8);
\draw [very thick, directed=1] (.8,1.2) -- (0,2);
\end{tikzpicture}};
\endxy
\right)
=
-q^{-1} \xy
(0,0)*{
\begin{tikzpicture} [scale=.5]
\draw [double, directed=.55] (0,0) -- (0,.7);
\draw [very thick, directed=1] (0,.7) -- (0,2);
\draw [very thick, directed=.55] (2,0) -- (2,1.3);
\draw [double, directed=1] (2,1.3) -- (2,2);
\draw [very thick, directed=.55] (0,.7) -- (2,1.3);
\end{tikzpicture}};
\endxy
\; ,\qquad
\Psi_{\T}\left(
\xy
(0,0)*{
\begin{tikzpicture} [scale=.5]
\draw [double, directed=1] (2,0) -- (0,2);
\draw [very thick] (0,0) -- (.8,.8);
\draw [very thick, directed=1] (1.2,1.2) -- (2,2);
\end{tikzpicture}};
\endxy
\right)
=
-q \xy
(0,0)*{
\begin{tikzpicture} [scale=.5]
\draw [double, directed=.55] (2,0) -- (2,.7);
\draw [very thick, directed=1] (2,.7) -- (2,2);
\draw [very thick, directed=.55] (0,0) -- (0,1.3);
\draw [double, directed=1] (0,1.3) -- (0,2);
\draw [very thick, directed=.55] (2,.7) -- (0,1.3);
\end{tikzpicture}};
\endxy
\\
\Psi_{\T}\left(
\xy
(0,0)*{
\begin{tikzpicture} [scale=.5]
\draw [very thick, directed=1] (0,0) -- (2,2);
\draw [double] (2,0) -- (1.2,.8);
\draw [double, directed=1] (.8,1.2) -- (0,2);
\end{tikzpicture}};
\endxy
\right)
=
-q^{-1} \xy
(0,0)*{
\begin{tikzpicture} [scale=.5]
\draw [very thick, directed=.55] (0,0) -- (0,1.3);
\draw [double, directed=1] (0,1.3) -- (0,2);
\draw [double, directed=.55] (2,0) -- (2,.7);
\draw [very thick, directed=1] (2,.7) -- (2,2);
\draw [very thick, directed=.55] (2,.7) -- (0,1.3);
\end{tikzpicture}};
\endxy
\; ,\qquad
\Psi_{\T}\left(
\xy
(0,0)*{
\begin{tikzpicture} [scale=.5]
\draw [very thick, directed=1] (2,0) -- (0,2);
\draw [double] (0,0) -- (.8,.8);
\draw [double, directed=1] (1.2,1.2) -- (2,2);
\end{tikzpicture}};
\endxy
\right)
=
-q \xy
(0,0)*{
\begin{tikzpicture} [scale=.5]
\draw [very thick, directed=.55] (2,0) -- (2,1.3);
\draw [double, directed=1] (2,1.3) -- (2,2);
\draw [double, directed=.55] (0,0) -- (0,.7);
\draw [very thick, directed=1] (0,.7) -- (0,2);
\draw [very thick, directed=.55] (0,.7) -- (2,1.3);
\end{tikzpicture}};
\endxy
\\
\Psi_{\T}\left(
\xy
(0,0)*{
\begin{tikzpicture} [scale=.5]
\draw [double, directed=1] (0,0) -- (2,2);
\draw [double] (2,0) -- (1.2,.8);
\draw [double, directed=1] (.8,1.2) -- (0,2);
\end{tikzpicture}};
\endxy
\right)
=
q^{-2} \xy
(0,0)*{
\begin{tikzpicture} [scale=.5]
\draw [double, directed=1] (0,0) .. controls (1,1) .. (0,2);
\draw [double, directed=1] (2,0) .. controls (1,1) .. (2,2);
\end{tikzpicture}};
\endxy
,\qquad
\Psi_{\T}\left(
\xy
(0,0)*{
\begin{tikzpicture} [scale=.5]
\draw [double, directed=1] (2,0) -- (0,2);
\draw [double] (0,0) -- (.8,.8);
\draw [double, directed=1] (1.2,1.2) -- (2,2);
\end{tikzpicture}};
\endxy
\right)
=
q^{2} \xy
(0,0)*{
\begin{tikzpicture} [scale=.5]
\draw [double, directed=1] (2,0) .. controls (1,1) .. (2,2);
\draw [double, directed=1] (0,0) .. controls (1,1) .. (0,2);
\end{tikzpicture}};
\endxy
\end{gather*}

We could check,
following~\cite{Kauf},
that the previous relations \eqref{SkeinCircles},
\eqref{SkeinDigons},
and the above smoothing relations def\/ine a~framed skein module.
Checking directly all formulas is rather long and tedious,
and we note that using the description in terms of the $U_q(\slm)$-action gives us an ef\/f\/icient way to considerably simplify the proof,
in the general case.
Indeed,
most formulas we want to check are consequences of $U_q(\slm)$-relations.

The previous skein module provides invariants of framed webs.
Here are the ef\/fects of adding a~negative twist on a~$1$-strand (depicted in a~ribbon version in the two left parts of the equation below):
\begin{gather*}
\Psi_{\T}\left(
\xy
(0,0)*{
\begin{tikzpicture} [scale=1,decoration={markings,
mark=at position 0.8 with {\arrow{>}}; }]
\draw (.1,0) -- (.1,.4) .. controls (.1,.5) and (-.1,.5) .. (-.1,.6) -- (-.1,.6) .. controls (-.1,.63) .. (-.05,.66);
\draw (.05,.74) .. controls (.1,.77) .. (.1,.8) -- (.1,1.2);
\draw (-.1,0) -- (-.1,.4) .. controls (-.1,.43) .. (-.05,.46);
\draw (.05,.54)..
controls (.1,.57) .. (.1,.6) -- (.1,.6) .. controls (.1,.7) and (-.1,.7) .. (-.1,.8) -- (-.1,1.2);
\end{tikzpicture}};
\endxy
\right)
=
\Psi_{\T}\left(
\xy
(0,0)*{
\begin{tikzpicture} [scale=1,decoration={markings,
mark=at position 0.8 with {\arrow{>}}; }]
\draw (-.1,0) .. controls (-.05,.2) .. (.06,.44);
\draw (.1,0) .. controls (.13,.17) .. (.2,.32);
\draw (.22,.7) .. controls (.4,1.1) and (.75,.8) .. (.75,.5) -- (.75,.5) .. controls (.75,.1)
 and (.2,.2) .. (.1,.5) -- (.1,.5) .. controls (-.05,.8) .. (-.1,1);
\draw (.32,.52) .. controls (.47,.86) and (.6,.56) .. (.6,.5) -- (.6,.5) .. controls (.6,.34)
 and (.38,.34) .. (.29,.52) -- (.29,.52) .. controls (.15,.8) .. (.1,1);
\end{tikzpicture}};
\endxy
\right)
=
-q \xy
(0,0)*{
\begin{tikzpicture} [scale=1,decoration={markings,
mark=at position 0.8 with {\arrow{>}}; }]
\draw [very thick] (0,0) .. controls (.3,.5) .. (0,1);
\draw [very thick] (.4,.5) arc (-180:180:.2);
\end{tikzpicture}};
\endxy
 \ +\xy
(0,0)*{
\begin{tikzpicture} [scale=1,decoration={markings,
mark=at position 0.8 with {\arrow{>}}; }]
\draw [very thick] (0,0) .. controls (0,.1) and (.1,.2) .. (.2,.35);
\draw [very thick] (.2,.35) .. controls (.3,0) and (.6,.3) .. (.6,.5) -- (.6,.5) .. controls (.6,.7) and (.3,1) .. (.2,.65);
\draw [double] (.2,.35) -- (.2,.65);
\draw [very thick] (.2,.65) .. controls (.1,.8) and (0,.9) .. (0,1);
\end{tikzpicture}};
\endxy
=
-q^{2}\xy
(0,0)*{
\begin{tikzpicture} [scale=1,decoration={markings,
mark=at position 0.8 with {\arrow{>}}; }]
\draw [very thick] (0,0) .. controls (.3,.5) .. (0,1);
\end{tikzpicture}};
\endxy
\end{gather*}

The same computation for a~$2$-labeled strand gives a~$q^{2}$ coef\/f\/icient.
We may introduce twists with half-integers,
assigning to them in the negative case the multiplication by $q$,
and in the positive case the multiplication by $q^{-1}$,
up to fourth roots of the unity.
We f\/ix the value to be $(-1)^{\frac{k}{2}}q^{\frac{-2k+k(k-1)}{2}}$, where $k$ stands for the labeling of the strand,
and where we have chosen a~favorite primitive fourth root of the unity $(-1)^{\frac{1}{2}}$.
The same formula will generalize to the $\sln$ case,
assigning to a~half twist on a~$k$-strand $(-1)^{\frac{k}{2}}q^{\frac{-nk+k(k-1)}{2}}$.

\subsubsection[$\mathfrak{sl}_n$ case]{$\boldsymbol{\mathfrak{sl}_n}$ case}

The previous process applies as well to any value of $n$,
and produces a~skein module in the following sense: Cautis,
Kamnitzer and Morrison~\cite{CKM} prove that all web relations come from $\U_q(\slm)$ relations,
so we just need to extend it to crossings and prove the invariance under Reidemeister moves.
This is the purpose of the following theorem,
which is the main result of the f\/irst part of this article.

\begin{Theorem}%\label{PlaneSkein}
For all $n$, the map $\Psi_{\T}$ induced by the local smoothing rules defined by $\T$ $($from relation~\eqref{DefBraiding}$)$
for positive crossings and $\T^{-1}$ for negative crossings extends to a~Kauffman bracket on $\sln$ web-tangles.
\end{Theorem}

In a~diagrammatic version,
the smoothing rules for $(k,l)$ crossings are thus given as follows:
\begin{gather*}
\Psi_{\T}\left(
\begin{tikzpicture}[scale=.8,anchorbase]
\draw [very thick,->] (0,0) -- (0,.5) -- (1,1.5) -- (1,2);
\draw [very thick] (1,0) -- (1,.5) -- (.65,.85);
\draw [very thick,->] (.35,1.15) -- (0,1.5) -- (0,2);
\node at (0,-.3) {\tiny $k$};
\node at (0,2.3) {};
\node at (1,-.3) {\tiny $l$};
\node at (1,2.3) {};
\end{tikzpicture}
\right)
=
\sum\limits_{a-b+c=k-l}(-1)^{b-l}q^{-ac+b-l} \;
\begin{tikzpicture}[scale=.8,anchorbase,decoration={markings,
mark=at position 0.8 with {\arrow{>}};}]
\draw [very thick,->] (0,0) -- (0,2);
\node at (0,-.3) {\tiny $k$};
\node at (0,2.3) {};
\node at (1,-.3) {\tiny $l$};
\node at (1,2.3) {};
\draw [very thick,->] (1,0) -- (1,2);
\draw [very thick,
postaction={decorate}] (0,.4) -- (1,.6);
\node at (.5,.2) {\tiny $c$};
\draw [very thick,
postaction={decorate}] (1,.9) -- (0,1.1);
\node at (.5,.8) {\tiny $b$};
\draw [very thick,
postaction={decorate}] (0,1.4) -- (1,1.6);
\node at (.5,1.8) {\tiny $a$};
\end{tikzpicture}
\\
\Psi_{\T}\left(
\begin{tikzpicture}[scale=.8,anchorbase]
\draw [very thick,->] (1,0) -- (1,.5) -- (0,1.5) -- (0,2);
\draw [very thick] (0,0) -- (0,.5) -- (.35,.85);
\draw [very thick,->] (.65,1.15) -- (1,1.5) -- (1,2);
\node at (0,-.3) {\tiny $k$};
\node at (0,2.3) {};
\node at (1,-.3) {\tiny $l$};
\node at (1,2.3) {};
\end{tikzpicture}
\right)
=
\sum\limits_{a-b+c=l-k}(-1)^{b+k}q^{ac-b+k} \;
\begin{tikzpicture}[scale=.8,anchorbase,decoration={markings,
mark=at position 0.8 with {\arrow{>}};}]
\draw [very thick,->] (0,0) -- (0,2);
\node at (0,-.3) {\tiny $k$};
\node at (0,2.3) {};
\node at (1,-.3) {\tiny $l$};
\node at (1,2.3) {};
\draw [very thick,->] (1,0) -- (1,2);
\draw [very thick,
postaction={decorate}] (1,.4) -- (0,.6);
\node at (.5,.2) {\tiny $c$};
\draw [very thick,
postaction={decorate}] (0,.9) -- (1,1.1);
\node at (.5,.8) {\tiny $b$};
\draw [very thick,
postaction={decorate}] (1,1.4) -- (0,1.6);
\node at (.5,1.8) {\tiny $a$};
\end{tikzpicture}
\end{gather*}

\begin{proof}
 \textbf{Braid-like relations:} Braid-like Reidemeister II relation is direct,
and braid-like Reidemeister III relations are consequences of the braiding relation (see~\cite{CMW} for a~presentation
of all~$6$ braid-like Reidemeister III relations).

 \textbf{Framing:} The framed Reidemeister I relation is easy to check.
Furthermore,
we can deal locally with the framing as we did in the $\slnn{2}$-case.
Details will be given in Lemma~\ref{lemma_framing}.

 \textbf{Braid-like web relation~\ref{ReidWebs1}:}\footnote{Note that
similar relations are studied at the categorif\/ied level in~\cite{MSVHomfly}.}
Braid-like relations \eqref{ReidWebs1} are consequences of the next equa\-li\-ty,
or similar ones:
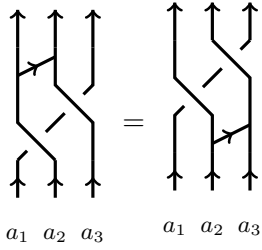
\begin{gather*}
\xy
(0,0)*{
\begin{tikzpicture} [scale=.5]
\draw [very thick, directed=.55] (0,0) -- (0,1);
\draw [very thick, directed=1] (0,4) -- (0,5);
\draw [very thick, directed=.55] (1,0) -- (1,1);
\draw [very thick, directed=1] (1,4) -- (1,5);
\draw [very thick, directed=.55] (2,0) -- (2,1);
\draw [very thick, directed=1] (2,4) -- (2,5);
\draw [very thick] (0,1) -- (.25,1.25);
\draw [very thick] (.75,1.75) -- (1,2) -- (1.25,2.25);
\draw [very thick] (1.75,2.75) -- (2,3) -- (2,4);
\draw [very thick] (1,1) -- (0,2) -- (0,4);
\draw [very thick] (2,1) -- (2,2) -- (1,3) -- (1,4);
\draw [very thick, directed=.55] (0,3.25) -- (1,3.75);
\node at (0,-1) {\text{\footnotesize{$a_1$}}};
\node at (1,-1) {\text{\footnotesize{$a_2$}}};
\node at (2,-1) {\text{\footnotesize{$a_3$}}};
\end{tikzpicture}};
\endxy
=
\xy
(0,0)*{
\begin{tikzpicture} [scale=.5]
\draw [very thick, directed=.55] (0,0) -- (0,1);
\draw [very thick, directed=1] (0,4) -- (0,5);
\draw [very thick, directed=.55] (1,0) -- (1,1);
\draw [very thick, directed=1] (1,4) -- (1,5);
\draw [very thick, directed=.55] (2,0) -- (2,1);
\draw [very thick, directed=1] (2,4) -- (2,5);
\draw [very thick] (0,1) -- (0,2) -- (.25,2.25);
\draw [very thick] (.75,2.75) -- (1,3) -- (1.25,3.25);
\draw [very thick] (1.75,3.75) -- (2,4);
\draw [very thick] (1,1) -- (1,2) -- (0,3) -- (0,4);
\draw [very thick] (2,1) -- (2,3) -- (1,4);
\draw [very thick, directed=.55] (1,1.25) -- (2,1.75);
\node at (0,-1) {\text{\footnotesize{$a_1$}}};
\node at (1,-1) {\text{\footnotesize{$a_2$}}};
\node at (2,-1) {\text{\footnotesize{$a_3$}}};
\end{tikzpicture}};
\endxy
\end{gather*}

The previous relation is a~diagrammatically depicted consequence of $C_{\Tw{1}\Tw{2}}(E_1)=E_2$,
a~relation from Propositions~\ref{propCTs} and~\ref{propCTp},
which still holds after rescalings.
For obtaining the general case,
one needs a~straightforward generalization of the previous relation: $C_{\Tw{1}\Tw{2}}(E_1^{(k)})=E_2^{(k)}$.

 \textbf{Star relations using duality:} Following~\cite{CMW},
it suf\/f\/ices to have the Reidemeister II relation with opposite orientations to deduce the following Reidemeister III star relations:
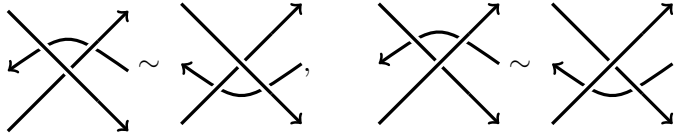
\begin{gather*}
\xy
(0,0)*{
\begin{tikzpicture} [scale=.8]
\draw [very thick,<-] (0,1) .. controls (1,1.7) .. (2,1);
\draw [very thick,
draw=white,
double=black,
double distance=1.25pt] (0,0) -- (1.5,1.5);
\draw [very thick,
->] (1.5,1.5) -- (2,2);
\draw [very thick,<-] (2,0) -- (1.5,.5);
\draw [very thick,
draw=white,
double=black,
double distance=1.25pt] (1.5,.5) -- (0,2);
\end{tikzpicture}};
\endxy
\sim \;
\xy
(0,0)*{
\begin{tikzpicture} [scale=.8]
\draw [very thick,<-] (0,1) .. controls (1,.3) .. (2,1);
\draw [very thick,
draw=white,
double=black,
double distance=1.25pt] (0,0) -- (1.5,1.5);
\draw [very thick,
->] (1.5,1.5) -- (2,2);
\draw [very thick,<-] (2,0) -- (1.5,.5);
\draw [very thick,
draw=white,
double=black,
double distance=1.25pt] (1.5,.5) -- (0,2);
\end{tikzpicture}};
\endxy
,\qquad
\xy
(0,0)*{
\begin{tikzpicture} [scale=.8,decoration={markings,
mark=at position 0.8 with {\arrow{>}}; }]
\draw [very thick,<-] (0,1) .. controls (1,1.7) .. (2,1);
\draw [very thick,
draw=white,
double=black,
double distance=1.25pt] (1.5,.5) -- (0,2);
\draw [very thick,
<-] (2,0) -- (1.5,.5);
\draw [very thick,
draw=white,
double=black,
double distance=1.25pt] (0,0) -- (1.5,1.5);
\draw [very thick,
->] (1.5,1.5) -- (2,2);
\end{tikzpicture}};
\endxy
\sim \;
\xy
(0,0)*{
\begin{tikzpicture} [scale=.8,decoration={markings,
mark=at position 0.8 with {\arrow{>}}; }]
\draw [very thick,<-] (0,1) .. controls (1,.3) .. (2,1);
\draw [very thick,
draw=white,
double=black,
double distance=1.25pt] (1.5,.5) -- (0,2);
\draw [very thick,
<-] (2,0) -- (1.5,.5);
\draw [very thick,
draw=white,
double=black,
double distance=1.25pt] (0,0) -- (1.5,1.5);
\draw [very thick,
->] (1.5,1.5) -- (2,2);
\end{tikzpicture}};
\endxy
\end{gather*}

We obtain the Reidemeister II case from the braid-like one as follows:
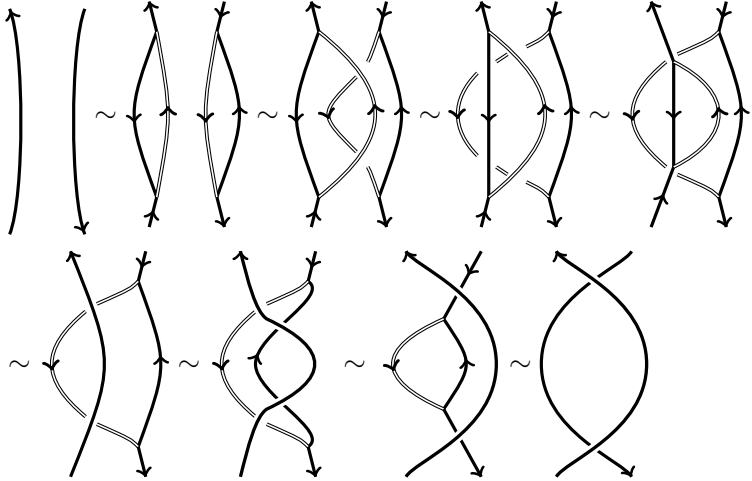
\begin{gather*}
\xy
(0,0)*{
\begin{tikzpicture} [scale=1]
\draw [very thick, directed=.99] (0,0) .. controls (.2,.5) and (.2,2.5) .. (0,3);
\draw [very thick, directed=.99] (1,3) .. controls (.8,2.5) and (.8,.5) .. (1,0);
\end{tikzpicture}};
\endxy
\sim
\xy
(0,0)*{
\begin{tikzpicture} [scale=1]
\draw [very thick, directed=.55] (0,0) -- (.1,.4);
\draw [very thick, directed=.55] (.1,2.6) .. controls (-.3,1.5) .. (.1,.4);
\draw [double, directed=.55] (.1,.4) .. controls (.3,1.5) .. (.1,2.6);
\draw [very thick, directed=1] (.1,2.6) -- (0,3);
\draw [very thick, directed=.55] (1,3) -- (.9,2.6);
\draw [very thick, directed=.55] (.9,.4) .. controls (1.3,1.5) .. (.9,2.6);
\draw [double, directed=.55] (.9,2.6) .. controls (.7,1.5) .. (.9,.4);
\draw [very thick, directed=1] (.9,.4) -- (1,0);
\end{tikzpicture}};
\endxy
\sim
\xy
(0,0)*{
\begin{tikzpicture} [scale=1]
\draw [very thick, directed=.55] (0,0) -- (.1,.4);
\draw [very thick, directed=.55] (.1,2.6) .. controls (-.3,1.5) .. (.1,.4);
\draw [double, directed=.55] (.1,.4) .. controls (1.1,1.2) and (1.1,1.8) .. (.1,2.6);
\draw [very thick, directed=1] (.1,2.6) -- (0,3);
\draw [very thick, directed=.55] (1,3) -- (.9,2.6);
\draw [very thick, directed=.55] (.9,.4) .. controls (1.3,1.5) .. (.9,2.6);
\draw [double] (.9,2.6) .. controls (.8,2.3) .. (.75,2.2);
\draw [double, directed=.55] (.6,2) .. controls (.1,1.5) .. (.6,1);
\draw [double] (.75,.8) .. controls (.8,.7) .. (.9,.4);
\draw [very thick, directed=1] (.9,.4) -- (1,0);
\end{tikzpicture}};
\endxy
\sim
\xy
(0,0)*{
\begin{tikzpicture} [scale=1]
\draw [very thick, directed=.55] (0,0) -- (.1,.4);
\draw [very thick, directed=.55] (.1,2.6) .. controls (.1,1.5) .. (.1,.4);
\draw [double, directed=.55] (.1,.4) .. controls (1.1,1.2) and (1.1,1.8) .. (.1,2.6);
\draw [very thick, directed=1] (.1,2.6) -- (0,3);
\draw [very thick, directed=.55] (1,3) -- (.9,2.6);
\draw [very thick, directed=.55] (.9,.4) .. controls (1.3,1.5) .. (.9,2.6);
\draw [double] (.9,2.6) .. controls (.8,2.5) .. (.6,2.4);
\draw [double] (.35,2.3) -- (.2,2.2);
\draw [double, directed=.55] (-.05,2.05) .. controls (-.4,1.7) and (-.4,1.3) .. (-.05,.95);
\draw [double] (.2,.8) -- (.35,.7);
\draw [double] (.6,.6) .. controls (.8,.5) .. (.9,.4);
\draw [very thick, directed=1] (.9,.4) -- (1,0);
\end{tikzpicture}};
\endxy
\sim
\xy
(0,0)*{
\begin{tikzpicture} [scale=1,]
\draw [very thick, directed=.55] (0,0) -- (.3,.8);
\draw [very thick, directed=.55] (.3,2.2) .. controls (.3,1.5) .. (.3,.8);
\draw [double, directed=.55] (.3,.8) .. controls (1.1,1.2) and (1.1,1.8) .. (.3,2.2);
\draw [very thick, directed=1] (.3,2.2) -- (0,3);
\draw [very thick, directed=.55] (1,3) -- (.9,2.6);
\draw [very thick, directed=.55] (.9,.4) .. controls (1.3,1.5) .. (.9,2.6);
\draw [double] (.9,2.6) .. controls (.8,2.5) .. (.35,2.3);
\draw [double, directed=.55] (.2,2.2) .. controls (-.4,1.7) and (-.4,1.3) .. (.2,.8);
\draw [double] (.35,.7) .. controls (.8,.5) .. (.9,.4);
\draw [very thick, directed=1] (.9,.4) -- (1,0);
\end{tikzpicture}};
\endxy
\\
\sim
\xy
(0,0)*{
\begin{tikzpicture} [scale=1]
\draw [very thick, directed=1] (0,0) .. controls (.6,1.5) .. (0,3);
\draw [very thick, directed=.55] (1,3) -- (.9,2.6);
\draw [very thick, directed=.55] (.9,.4) .. controls (1.3,1.5) .. (.9,2.6);
\draw [double] (.9,2.6) .. controls (.8,2.5) .. (.35,2.3);
\draw [double, directed=.55] (.2,2.2) .. controls (-.4,1.7) and (-.4,1.3) .. (.2,.8);
\draw [double] (.35,.7) .. controls (.8,.5) .. (.9,.4);
\draw [very thick, directed=1] (.9,.4) -- (1,0);
\end{tikzpicture}};
\endxy
\sim
\xy
(0,0)*{
\begin{tikzpicture} [scale=1]
\draw [very thick, directed=.55] (1,3) -- (.9,2.6);
\draw [very thick, directed=.55] (.9,.4) .. controls (1.2,.6) and (.2,1.1) .. (.2,1.5) -- (.2,1.5) .. controls (.2,1.9) and (1.2,2.4) .. (.9,2.6);
\draw [double] (.9,2.6) .. controls (.8,2.5) .. (.35,2.3);
\draw [double, directed=.55] (.2,2.2) .. controls (-.4,1.7) and (-.4,1.3) .. (.2,.8);
\draw [double] (.35,.7) .. controls (.8,.5) .. (.9,.4);
\draw [very thick, directed=1] (.9,.4) -- (1,0);
\draw [very thick,
draw=white,
double,
double distance =1.25 pt,
double=black] (0,0) .. controls (.1,.4) and (.2,.8) .. (.35,.9) -- (.35,.9) .. controls (1.2,1.3)
 and (1.2,1.7) .. (.35,2.1) -- (.35,2.1) .. controls (.2,2.2) and (.1,2.6) .. (0,3);
\draw [very thick, directed=1] (0,2.99)-- (0,3);
\end{tikzpicture}};
\endxy
\sim
\xy
(0,0)*{
\begin{tikzpicture} [scale=1]
\draw [very thick, directed=.35] (1,3) -- (.5,2.1);
\draw [very thick, directed=.55] (.5,.9) .. controls (.9,1.5) .. (.5,2.1);
\draw [double, directed=.55] (.5,2.1) .. controls (-.4,1.7) and (-.4,1.3) .. (.5,.9);
\draw [very thick, directed=1] (.5,.9) -- (1,0);
\draw [very thick,
draw=white,
double,
double distance =1.25 pt,
double=black] (0,0) .. controls (.1,.2) and (1.2,.6) .. (1.2,1.5) -- (1.2,1.5) .. controls (1.2,2.4) and (.1,2.8) .. (0,3);
\draw [very thick, directed=1] (0.01,2.98)-- (0,3);
\end{tikzpicture}};
\endxy
\sim
\xy
(0,0)*{
\begin{tikzpicture} [scale=1]
\draw [very thick, directed=1] (1,3) .. controls (.9,2.8) and (-.2,2.4) .. (-.2,1.5) -- (-.2,1.5) .. controls (-.2,.6) and (.9,.2) .. (1,0);
\draw [very thick,
draw=white,
double,
double distance =1.25 pt,
double=black] (0,0) .. controls (.1,.2) and (1.2,.6) .. (1.2,1.5) -- (1.2,1.5) .. controls (1.2,2.4) and (.1,2.8) .. (0,3);
\draw [very thick, directed=1] (0.01,2.98)-- (0,3);
\end{tikzpicture}};
\endxy
\end{gather*}

We used in the previous computation the star Reidemeister II relation in the $n$-$n$ case,
which is easy to prove:
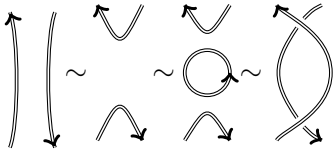
\begin{gather*}
\xy
(0,0)*{
\begin{tikzpicture} [scale=.6]
\draw [double, directed=.99] (0,0) .. controls (.2,.5) and (.2,2.5) .. (0,3);
\draw [double, directed=.99] (1,3) .. controls (.8,2.5) and (.8,.5) .. (1,0);
\end{tikzpicture}};
\endxy
\sim
\xy
(0,0)*{
\begin{tikzpicture} [scale=.6]
\draw [double, directed=1] (0,0) .. controls (.5,1) .. (1,0);
\draw [double, directed=1] (1,3) .. controls (.5,2)..
(0,3);
\end{tikzpicture}};
\endxy
\sim
\xy
(0,0)*{
\begin{tikzpicture} [scale=.6]
\draw [double, directed=1] (0,0) .. controls (.5,.8) .. (1,0);
\draw [double, directed=1] (1,3) .. controls (.5,2.2)..
(0,3);
\draw [double, directed=1] (1,1.5) arc(0:360:.5);
\end{tikzpicture}};
\endxy
\sim
\xy
(0,0)*{
\begin{tikzpicture} [scale=.6]
\draw [double] (1,3) .. controls (.8,2.9) .. (.6,2.7);
\draw [double] (.45,2.5) .. controls (-.1,1.7) and (-.1,1.3) .. (.45,.5);
\draw [double,directed=1] (.6,.3) .. controls (.8,.1) .. (1,0);
\draw [double, directed=1] (0,0) .. controls (.1,.2) and (1.2,.6) .. (1.2,1.5) -- (1.2,1.5) .. controls (1.2,2.4) and (.1,2.8) .. (0,3);
\end{tikzpicture}};
\endxy
\end{gather*}

Then,
we want to obtain the missing forms of relation \eqref{ReidWebs1}.
We proceed as follows in one case,
the other ones being similar:
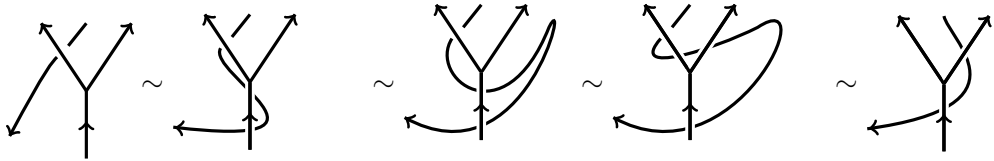
\begin{gather*}
\xy
(0,0)*{
\begin{tikzpicture} [scale=.6]
\draw [very thick, directed=.55] (1,0) -- (1,1.5);
\draw [very thick, directed=1] (1,1.5) -- (0,3);
\draw [very thick, directed=1] (1,1.5) -- (2,3);
\draw [very thick] (1,3) -- (.6,2.5);
\draw [very thick, directed=1] (.35,2.25) .. controls (.25,2.15) and (.1,2) .. (-.7,.5);
\end{tikzpicture}};
\endxy
\sim
\xy
(0,0)*{
\begin{tikzpicture} [scale=.6]
\draw [very thick, directed=1] (1,1.5) -- (0,3);
\draw [very thick, directed=1] (1,1.5) -- (2,3);
\draw [very thick] (1,3) -- (.6,2.5);
\draw [very thick, directed=.99] (.35,2.25) .. controls (0,1.7) and (3.5,0) .. (-.7,.5);
\draw [very thick,
double,
draw = white,
double= black,
double distance=1.25pt] (1,0) -- (1,1.5);
\draw [very thick, directed=.55] (1,0) -- (1,1.5);
\end{tikzpicture}};
\endxy
\sim
\xy
(0,0)*{
\begin{tikzpicture} [scale=.6]
\draw [very thick, directed=1] (1,1.5) -- (0,3);
\draw [very thick, directed=1] (1,1.5) -- (2,3);
\draw [very thick] (1,3) -- (.6,2.5);
\draw [very thick, directed=1] (.35,2.25) .. controls (-.2,1.4) and (1.5,0) .. (2.5,2.5) -- (2.5,2.5) .. controls (3,3.5) and (2,-1) ..(-.7,.5);
\draw [very thick,
double,
draw = white,
double= black,
double distance=1.25pt] (1,0) -- (1,1.5);
\draw [very thick, directed=.55] (1,0) -- (1,1.5);
\end{tikzpicture}};
\endxy
\sim
\xy
(0,0)*{
\begin{tikzpicture} [scale=.6]
\draw [very thick] (1,3) -- (.6,2.5);
\draw [very thick, directed=1] (.35,2.25) .. controls (-.4,1.4) and (1.5,2) .. (2.5,2.5) -- (2.5,2.5) .. controls (4,3.5) and (2,-1) ..(-.7,.5);
\draw [very thick,
double,
draw = white,
double= black,
double distance=1.25pt] (1,0) -- (1,1.5);
\draw [very thick, directed=.55] (1,0) -- (1,1.5);
\draw [very thick,
double,
draw = white,
double= black,
double distance=1.25pt] (1,1.5) -- (0,3);
\draw [very thick,
double,
draw = white,
double= black,
double distance=1.25pt] (1,1.5) -- (2,3);
\draw [very thick, directed=1] (1,1.5) -- (0,3);
\draw [very thick, directed=1] (1,1.5) -- (2,3);
\end{tikzpicture}};
\endxy
\sim
\xy
(0,0)*{
\begin{tikzpicture} [scale=.6]
\draw [very thick, directed=.99] (1,3) .. controls (1.1,2.5) and (3,1) ..(-.7,.5);
\draw [very thick,
double,
draw = white,
double= black,
double distance=1.25pt] (1,0) -- (1,1.5);
\draw [very thick, directed=.55] (1,0) -- (1,1.5);
\draw [very thick,
double,
draw = white,
double= black,
double distance=1.25pt] (1,1.5) -- (0,3);
\draw [very thick,
double,
draw = white,
double= black,
double distance=1.25pt] (1,1.5) -- (2,3);
\draw [very thick, directed=1] (1,1.5) -- (0,3);
\draw [very thick, directed=1] (1,1.5) -- (2,3);
\end{tikzpicture}};
\endxy
\end{gather*}

\textbf{Relation (\ref{ReidWebs2}) and framing:} The last relation,
for which we need a~better understanding of the framing,
will be proved separately in Lemma~\ref{lemma_R2webs} below.
\end{proof}

It is easy to see that a~positive framing on a~$k$-strand is equivalent to multiplication by a~polynomial in $q$ and $q^{-1}$.
Denote $t_k$ this polynomial:
\begin{gather*}
\Psi_{\T}\left(
\xy
(0,0)*{
\begin{tikzpicture} [scale=.5]
\draw [very thick, directed=1] (0,0) -- (0,2);
\node at (0,1) {$\circ$};
\node at (.3,1) {\tiny 1};
\node at (0,-.3) {\tiny k};
\node at (0,2.3) {\vphantom{\tiny k}};
\end{tikzpicture}};
\endxy
\right)
\; = \; t_k\;
\xy
(0,0)*{
\begin{tikzpicture} [scale=.5,decoration={markings,
mark=at position 1 with {\arrow{>}}; }]
\draw [very thick,
postaction={decorate}] (0,0) -- (0,2);
\node at (0,-.3) {\tiny k};
\node at (0,2.3) {\vphantom{\tiny k}};
\end{tikzpicture}};
\endxy
\end{gather*}

\begin{Lemma} \label{lemma_framing}
$t_{k}=(-1)^k q^{-kn}q^{k(k-1)}$.
\end{Lemma}

Note that this formula explains the choice for the half twists.

\begin{proof}
We claim that $t_{k+1}t_{k}^{-1}=-q^{2k-n}$.
Indeed,
the following relation holds from already proven Kauf\/fman relations:
\begin{gather*}
\xy
(0,0)*{
\begin{tikzpicture} [scale=.5]
\draw [very thick, directed=.55] (0,0) -- (0,1);
\draw [very thick] (0,1) .. controls (0,2) and (1,2) .. (1,1) -- (1,1) .. controls (1,.2) and (.4,0) .. (.2,.8);
\draw [very thick] (0,1.5) -- (0,2);
\draw [very thick] (0,2) -- (-1,3);
\draw [very thick] (0,2) -- (1,3);
\draw [very thick, directed=1] (-1,3) -- (-1,4);
\draw [very thick, directed=1] (1,3) -- (1,4);
\node at (-.8,.5) {\tiny $k+1$};
\node at (-1.3,3.5) {\tiny $1$};
\node at (1.3,3.5) {\tiny $k$};
\end{tikzpicture}};
\endxy
=
\xy
(0,0)*{
\begin{tikzpicture} [scale=.5]
\draw [very thick, directed=.55] (0,0) -- (0,1);
\draw [very thick] (0,1)..
controls (-.2,1.2) and (-.5,1.6) .. (-.5,2);
\draw [very thick] (0,1)..
controls (.2,1.2) and (.5,1.6) .. (.5,2);
\draw [very thick] (-.5,2) .. controls (-1,3.5) and (3,4.5) .. (3,3) -- (3,3) .. controls (3,1) and (.8,1.3) .. (.7,2.1);
\draw [very thick] (.5,2) .. controls (.5,3) and (2.5,3.5) .. (2.5,2.8) -- (2.5,2.8) .. controls (2.5,1.5) and (1.3,1.8) .. (1.2,2.6);
\draw [very thick] (1.1,3.1) -- (1.05,3.4);
\draw [very thick] (.55,2.6) -- (.4,3.1);
\draw [very thick, directed=1] (1,3.8) -- (1,5);
\draw [very thick, directed=1] (.3,3.5) -- (0,5);
\end{tikzpicture}};
\endxy
=
\xy
(0,0)*{
\begin{tikzpicture} [scale=.5]
\draw [very thick, directed=.55] (0,0) -- (0,1);
\draw [very thick] (0,1) -- (-.5,2) -- (.5,3) -- (.2,3.3);
\draw [very thick] (-.2,3.7) -- (-.5,4);
\draw [very thick] (0,1) -- (.5,2) -- (.2,
2.3);
\draw [very thick] (-.2,2.7) -- (-.5,3) -- (.5,4);
\draw [very thick, directed=1] (-.5,4) -- (-.5,5);
\draw [very thick, directed=1] (.5,4) -- (.5,5);
\node at (-.5,4.2) {$\circ$};
\node at (.5,4.2) {$\circ$};
\end{tikzpicture}};
\endxy
\end{gather*}
The equality of the l.h.s.\ and r.h.s\ parts implies the recurrence relation.
The general solution follows then from the computation of the value for a~$1$-strand.
An explicit computation for this gives $t_1=-q^{-n}$.
\end{proof}

We are now ready to prove the last relation by induction:

\begin{Lemma} \label{lemma_R2webs}
\begin{gather*}
\Psi_{\T}\left(
\xy
(0,0)*{
\begin{tikzpicture} [scale=.5]
\draw [very thick, directed=.8] (0,0) -- (0,1);
\draw [very thick] (0,1) -- (-.5,2) -- (.5,3);
\draw [very thick] (0,1) -- (.5,2) -- (.2,
2.3);
\draw [very thick] (-.2,2.7) -- (-.5,3);
\draw [very thick, directed=1] (-.5,3) -- (-.5,4);
\draw [very thick, directed=1] (.5,3) -- (.5,4);
\node at (-.4,.3) {$k$};
\node at (-.8,3.3) {$r$};
\node at (.8,3.3) {$l$};
\end{tikzpicture}};
\endxy
\right)
=
\Psi_{\T}\left(
\xy
(0,0)*{
\begin{tikzpicture} [scale=.5]
\draw [very thick] (0,0) -- (0,1);
\draw [very thick] (0,1) -- (-.5,2);
\draw [very thick] (0,1) -- (.5,2);
\draw [very thick, directed=1] (-.5,2) -- (-.5,3.5);
\draw [very thick, directed=1] (.5,2) -- (.5,3.5);
\node at (-.4,.3) {$k$};
\node at (-.8,2.3) {$r$};
\node at (.8,2.1) {$l$};
\node at (0,.7) {$\circ$};
\node at (-.5,2.7) {$\circ$};
\node at (.5,2.7) {$\circ$};
\node at (.4,.7) {\tiny$\frac{1}{2}$};
\node at (-.1,2.7) {\tiny$\frac{-1}{2}$};
\node at (1,2.7) {\tiny$\frac{-1}{2}$};
\end{tikzpicture}};
\endxy
\right)
\end{gather*}
\end{Lemma}

\begin{proof}
The computation is easy for $r=1$ or $l=1$.
Then we use:
\begin{gather*}
\Psi_{\T}\left(
\xy
(0,0)*{
\begin{tikzpicture} [scale=.5]
\draw [very thick, directed=.8] (0,0) -- (0,1);
\draw [very thick] (0,1) -- (-.5,2) -- (.5,3);
\draw [very thick] (0,1) -- (.5,2) -- (.2,
2.3);
\draw [very thick] (-.2,2.7) -- (-.5,3);
\draw [very thick, directed=1] (-.5,3) -- (-.5,5);
\draw [very thick, directed=1] (.5,3) -- (.5,5);
\node at (-.4,.3) {$k$};
\node at (-.8,3.3) {$r$};
\node at (.8,3.3) {$l$};
\end{tikzpicture}};
\endxy
\right)
=\frac{1}{[l]}
\Psi_{\T}\left(
\xy
(0,0)*{
\begin{tikzpicture} [scale=.5]
\draw [very thick, directed=.8] (0,0) -- (0,1);
\draw [very thick] (0,1) -- (-.5,2) -- (.5,3);
\draw [very thick] (0,1) -- (.5,2) -- (.2,
2.3);
\draw [very thick] (-.2,2.7) -- (-.5,3);
\draw [very thick, directed=1] (-.5,3) -- (-.5,5);
\draw [very thick, directed=.55] (.5,3) .. controls (0,3.5) .. (.5,4);
\node [rotate=90] at (-.2,3.5) {\tiny$l-1$};
\draw [very thick, directed=.55] (.5,3) .. controls (1,3.5) .. (.5,4);
\node at (1.2,3.5) {\tiny$1$};
\draw [very thick, directed=1] (.5,4) -- (.5,5);
\node at (-.4,.3) {$k$};
\end{tikzpicture}};
\endxy
\right)
\; = \; \frac{1}{[l]} \;
\Psi_{\T}\left(
\xy
(0,0)*{
\begin{tikzpicture} [scale=.5]
\draw [very thick, directed=.8] (0,0) -- (0,1);
\draw [very thick] (0,1) -- (-.5,2) .. controls (-.2,3) .. (.5,3);
\draw [very thick] (-.5,2) .. controls (.2,2) .. (.5,3);
\draw [very thick] (.5,3) -- (.5,4);
\draw [very thick] (0,1) -- (.5,2) -- (.3,
2.2);
\draw [very thick] (.1,2.4) -- (-.1,2.6);
\draw [very thick] (-.3,2.8) -- (-.5,3);
\draw [very thick, directed=1] (-.5,3) -- (-.5,5);
\draw [very thick, directed=1] (.5,4) -- (.5,5);
\node at (-.4,.3) {$k$};
\node at (-.8,3.5) {$r$};
\node at (.8,3.5) {$l$};
\end{tikzpicture}};
\endxy
\right)
\; = \; \frac{1}{[l]}\;
\Psi_{\T}\left(
\xy
(0,0)*{
\begin{tikzpicture} [scale=.5]
\draw [very thick, directed=.8] (0,0) -- (0,1);
\draw [very thick] (0,1) -- (-.5,2) -- (.5,4);
\draw [very thick] (0,1) -- (.5,2) .. controls (0,2.2) and (.8,2.9) .. (.5,4);
\draw [very thick] (.5,2) .. controls (1,2.4) .. (.65,2.8);
\draw [very thick] (.3,3) -- (.2,3.1);
\draw [very thick] (-.1,3.3) .. controls (-.3,3.6) and (-.5,3.8)..
(-.5,4);
\draw [very thick, directed=1] (-.5,4) -- (-.5,5);
\draw [very thick, directed=1] (.5,4) -- (.5,5);
\node at (-.4,.3) {$k$};
\node at (-.8,4.5) {$r$};
\node at (.8,4.5) {$l$};
\end{tikzpicture}};
\endxy
\right)
\\
\qquad
=\frac{a_{r+1}^{\frac{1}{2}}a_r^{\frac{-1}{2}}a_{1}^{\frac{-1}{2}}}{[l]}
\Psi_{\T}\left(
\xy
(0,0)*{
\begin{tikzpicture} [scale=.5]
\draw [very thick, directed=.8] (0,0) -- (0,1);
\draw [very thick] (0,1) -- (-.5,2) -- (.5,4);
\draw [very thick] (0,1) -- (.5,2) .. controls (1,3) .. (.5,4);
\draw [very thick] (.5,2) -- (.2,3);
\draw [very thick] (-.1,3.3) .. controls (-.3,3.6) and (-.5,3.8)..
(-.5,4);
\draw [very thick, directed=1] (-.5,4) -- (-.5,5);
\draw [very thick, directed=1] (.5,4) -- (.5,5);
\node at (-.4,.3) {$k$};
\node at (-.8,4.5) {$r$};
\node at (.8,4.5) {$l$};
\end{tikzpicture}};
\endxy
\right)
=\frac{a_{r+1}^{\frac{1}{2}}a_r^{\frac{-1}{2}}a_{1}^{\frac{-1}{2}}}{[l]} \quad
\Psi_{\T}\left(
\xy
(0,0)*{
\begin{tikzpicture} [scale=.5]
\draw [very thick, directed=.8] (0,0) -- (0,1);
\draw [very thick] (0,1) -- (-.5,2) -- (-.8,2.3) .. controls (-1,2.5) .. (.5,4);
\draw [very thick] (-.5,2) .. controls (-.1,2.5) .. (-.4,2.9);
\draw [very thick] (-.6,3.2) .. controls (-.8,3.5) .. (-.5,4);
\draw [very thick] (0,1) -- (.5,2) .. controls (1,3) .. (.5,4);
\draw [very thick, directed=1] (-.5,4) -- (-.5,5);
\draw [very thick, directed=1] (.5,4) -- (.5,5);
\node at (-.4,.3) {$k$};
\node at (-.8,4.5) {$r$};
\node at (.8,4.5) {$l$};
\end{tikzpicture}};
\endxy
\right)
\\
\qquad
=\frac{a_{r+1}^{\frac{1}{2}}a_r^{\frac{-1}{2}}a_{1}^{\frac{-1}{2}}a_{k-1}^{\frac{1}{2}}a_{l-1}^{\frac{-1}{2}}a_r^{\frac{-1}{2}}}{[l]}
\xy
(0,0)*{
\begin{tikzpicture} [scale=.5]
\draw [very thick, directed=.8] (0,0) -- (0,1);
\draw [very thick] (0,1) -- (-.5,2) -- (-.5,4);
\draw [very thick] (-.5,2) -- (.5,3);
\draw [very thick] (0,1) -- (.5,2) -- (.5,4);
\draw [very thick, directed=1] (-.5,4) -- (-.5,5);
\draw [very thick, directed=1] (.5,4) -- (.5,5);
\node at (-.4,.3) {$k$};
\node at (-.8,4.5) {$r$};
\node at (.8,4.5) {$l$};
\end{tikzpicture}};
\endxy
= a_{r+1}^{\frac{1}{2}}a_r^{\frac{-1}{2}}a_{1}^{\frac{-1}{2}}a_{k-1}^{\frac{1}{2}}a_{l-1}^{\frac{-1}{2}}a_r^{\frac{-1}{2}} \xy
(0,0)*{
\begin{tikzpicture} [scale=.5]
\draw [very thick, directed=.8] (0,0) -- (0,1);
\draw [very thick] (0,1) -- (-.5,2) -- (-.5,3);
\draw [very thick] (0,1) -- (.5,2) -- (.5,3);
\draw [very thick, directed=1] (-.5,3) -- (-.5,4);
\draw [very thick, directed=1] (.5,3) -- (.5,4);
\node at (-.4,.3) {$k$};
\node at (-.8,4.5) {$r$};
\node at (.8,4.5) {$l$};
\end{tikzpicture}};
\endxy
\end{gather*}

An explicit computation of the coef\/f\/icient shows that the previous term equals:
\begin{gather*}
\Psi_{\T}\left(
\xy
(0,0)*{
\begin{tikzpicture} [scale=.5]
\draw [very thick, directed=.8] (0,0) -- (0,1);
\draw [very thick] (0,1) -- (0,1.5);
\draw [very thick] (0,1.5) -- (-.5,2.5) -- (-.5,3);
\draw [very thick] (0,1.5) -- (.5,2.5) -- (.5,3);
\draw [very thick, directed=1] (-.5,3) -- (-.5,4);
\draw [very thick, directed=1] (.5,3) -- (.5,4);
\node at (-.4,.3) {$k$};
\node at (-.8,4.5) {$r$};
\node at (.8,4.5) {$l$};
\node at (0,1) {$\circ$};
\node at (.5,1) {$\frac{1}{2}$};
\node at (-.5,3) {$\circ$};
\node at (-1.2,3) {$-\frac{1}{2}$};
\node at (.5,3) {$\circ$};
\node at (1.2,3) {$-\frac{1}{2}$};
\end{tikzpicture}};
\endxy
\right),
\end{gather*}
which completes the proof.
\end{proof}

We therefore obtain a~well-def\/ined skein module providing an invariant of knotted webs.

Note that in the $\mathfrak{sl}_3$ case,
$2$-strands are usually translated into $1$-strands by reversing the orientation.
In this case,
smoothings of crossings would be def\/ined only up to a~power of $q$,
and understanding a~skew-Howe based way to f\/ix this power seems dif\/f\/icult.
We choose here not to apply this duality process and keep distinct $1$- and $2$-strands with their own orientations,
and more generally to keep all strands numbered $1,\dots, n$ in the $\sln$ case,
with their orientation.

So,
we have seen that to a~web-tangle in ladder position,
we can assign a~$U_q(\sln)$ morphism of tensor product of minuscule representations.
The diagrammatic form of this morphism corresponds to the image of the same web-tangle in the $\sln$ skein module.
If we start with a~non-ladder web,
we can assign to it its skein element,
but the skew-Howe process does not directly apply.
Cautis,
Kamnitzer and Morrison~\cite{CKM} explain a~process for turning upward webs to ladder form,
which we summarize in the following section.

\subsubsection{Turning a~knot to a~ladder}

Let us now consider a~tangle $T$ (possibly a~web-tangle) with only upward boundaries.
Following Cautis--Kamnitzer--Morrison,
we can present it as:
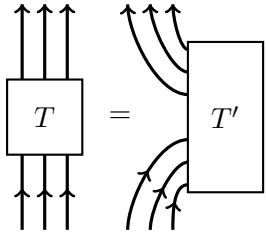
\begin{gather*}
\xy
(0,0)*{
\begin{tikzpicture}
\draw [very thick, directed=.55] (0.2,0) -- (.2,1);
\draw [very thick, directed=.55] (0.5,0) -- (.5,1);
\draw [very thick, directed=.55] (0.8,0) -- (.8,1);
\draw [thick] (0,1) rectangle (1,2);
\node at (.5,1.5) {$T$};
\draw [very thick, directed=1] (0.2,2) -- (.2,3);
\draw [very thick, directed=1] (0.5,2) -- (.5,3);
\draw [very thick, directed=1] (0.8,2) -- (.8,3);
\end{tikzpicture}};
(10,0)*{=};
(20,0)*{
\begin{tikzpicture}
\draw [very thick, directed=.55] (0.2,0) .. controls (.2,.4) and (.6,1.2) .. (1,1.2);
\draw [very thick, directed=.55] (0.5,0) .. controls (0.5,.3) and (.8,.9) .. (1,.9);
\draw [very thick, directed=.55] (0.8,0) .. controls (0.8,.2) and (.8,.6) .. (1,.6);
\draw [thick] (1,.5) rectangle (2,2.5);
\node at (1.5,1.5) {$T'$};
\draw [very thick, directed=.99] (1,1.8) .. controls (.6,1.8) and (.2,2.6) .. (.2,3);
\draw [very thick, directed=1] (1,2.1) .. controls (.8,
2.1) and (.5,2.7) .. (.5,3);
\draw [very thick, directed=1] (1,2.4) .. controls (.9,2.4) and (.8,2.8) .. (.8,3);
\end{tikzpicture}};
\endxy
\end{gather*}

The left part of the r.h.s.\ is easily presentable as a~ladder.
The tangle $T'$ is assumed to be presented as a~horizontal grid diagram generated by caps,
cups and crossings (plus 3-valent vertices for webs).
We request to have all crossings vertical,
which is possible up to some isotopy.
So, two caps or cups cannot lie one over the other one,
and we determine the number of $n$-strands we will have to add as the number of elementary pieces that contain
a~downward strand: we will then put a~strand on the right of this place.
Let this number be denoted $\alpha$.
This being done, start over,
but adjoining on the right of $T$ $\alpha$ upward $n$-strands placed at the right place
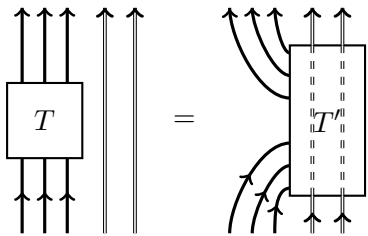
\begin{gather*}
\xy
(0,0)*{
\begin{tikzpicture}
\draw [very thick, directed=.55] (0.2,0) -- (.2,1);
\draw [very thick, directed=.55] (0.5,0) -- (.5,1);
\draw [very thick, directed=.55] (0.8,0) -- (.8,1);
\draw [thick] (0,1) rectangle (1,2);
\node at (.5,1.5) {$T$};
\draw [very thick, directed=1] (0.2,2) -- (.2,3);
\draw [very thick, directed=1] (0.5,2) -- (.5,3);
\draw [very thick, directed=1] (0.8,2) -- (.8,3);
\draw [double, directed=1] (1.3,0) -- (1.3,3);
\draw [double, directed=1] (1.7,0) -- (1.7,3);
\end{tikzpicture}};
(15,0)*{=};
(30,0)*{
\begin{tikzpicture}
\draw [very thick, directed=.55] (0.2,0) .. controls (.2,.4) and (.6,1.2) .. (1,1.2);
\draw [very thick, directed=.55] (0.5,0) .. controls (0.5,.3) and (.8,.9) .. (1,.9);
\draw [very thick, directed=.55] (0.8,0) .. controls (0.8,.2) and (.8,.6) .. (1,.6);
\draw [thick] (1,.5) rectangle (2,2.5);
\node at (1.5,1.5) {$T'$};
\draw [very thick, directed=.99] (1,1.8) .. controls (.6,1.8) and (.2,2.6) .. (.2,3);
\draw [very thick, directed=1] (1,2.1) .. controls (.8,
2.1) and (.5,2.7) .. (.5,3);
\draw [very thick, directed=1] (1,2.4) .. controls (.9,2.4) and (.8,2.8) .. (.8,3);
\draw [double, directed=.55] (1.3,0) -- (1.3,.5);
\draw [double,
dashed] (1.3,.5) -- (1.3,2.5);
\draw [double, directed=1] (1.3,2.5) -- (1.3,3);
\draw [double, directed=.55] (1.7,0) -- (1.7,.5);
\draw [double,
dashed] (1.7,.5) -- (1.7,2.5);
\draw [double, directed=1] (1.7,2.5) -- (1.7,3);
\end{tikzpicture}};
\endxy
\end{gather*}

By performing some moves and simplif\/ications near the downward strands,
we get a~ladder $L$.
These local changes Cautis--Kamnitzer--Morrison perform are smoothings and simplif\/ications of some Reidemeister moves,
and so the image of the tangle is equal in the previous skein module to $T$ with $\alpha$ disjoint $n$-strands added to it.

For example,
if we start from the elementary web we considered in the introduction:
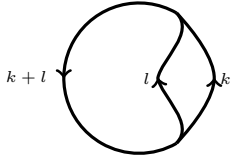
\begin{gather*}
\xy
(0,0)*{
\begin{tikzpicture} [scale=.5,decoration={markings,
mark=at position .5 with {\arrow{>}}; }]
\draw[postaction={decorate}, very thick] (1,1.73) arc (60:300:2);
\draw[postaction={decorate}, very thick] (1,-1.73) .. controls (1.5,-1.2) and (.6,-.5) .. (.5,0) -- (.5,0) .. controls (.6,.5) and (1.5,1.2) .. (1,1.73);
\draw[postaction={decorate}, very thick] (1,-1.73) .. controls (1.5,-1.2) and (2,-.5) .. (2,0) -- (2,0) .. controls (2,.5) and (1.5,1.2) .. (1,1.73);
\node at (-3,0) {\tiny $k+l$};
\node at (2.3,0) {\tiny $k$};
\node at (.2,0) {\tiny $l$};
\end{tikzpicture}};
\endxy
\end{gather*}
we can add on the right one $n$-labeled strand and perform Reidemeister--Kauf\/fman moves:
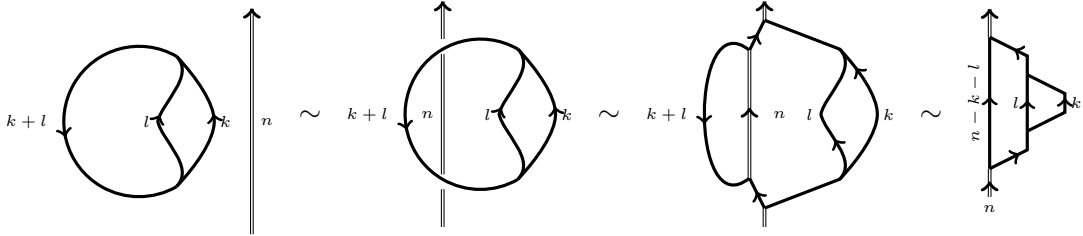
\begin{gather*}
\xy
(0,0)*{
\begin{tikzpicture} [scale=.5]
\draw[directed=.55,
very thick] (1,1.73) arc (60:300:2);
\draw[directed=.55,
very thick] (1,-1.73) .. controls (1.5,-1.2) and (.6,-.5) .. (.5,0) -- (.5,0) .. controls (.6,.5) and (1.5,1.2) .. (1,1.73);
\draw[directed=.55,
very thick] (1,-1.73) .. controls (1.5,-1.2) and (2,-.5) .. (2,0) -- (2,0) .. controls (2,.5) and (1.5,1.2) .. (1,1.73);
\draw[directed=1,
double] (3,-3) -- (3,3);
\node at (-3,0) {\tiny $k+l$};
\node at (2.3,0) {\tiny $k$};
\node at (.2,0) {\tiny $l$};
\node at (3.4,0) {\tiny $n$};
\end{tikzpicture}};
\endxy
\; \sim \;
\xy
(0,0)*{
\begin{tikzpicture} [scale=.5]
\draw[directed=.55,
very thick] (1,1.73) arc (60:300:2);
\draw[directed=.55,
very thick] (1,-1.73) .. controls (1.5,-1.2) and (.6,-.5) .. (.5,0) -- (.5,0) .. controls (.6,.5) and (1.5,1.2) .. (1,1.73);
\draw[directed=.55,
very thick] (1,-1.73) .. controls (1.5,-1.2) and (2,-.5) .. (2,0) -- (2,0) .. controls (2,.5) and (1.5,1.2) .. (1,1.73);
\draw[double] (-1,-3) -- (-1,-2);
\draw[double] (-1,-1.6) -- (-1,1.6);
\draw[directed=1,
double] (-1,2) -- (-1,3);
\node at (-3,0) {\tiny $k+l$};
\node at (2.3,0) {\tiny $k$};
\node at (.2,0) {\tiny $l$};
\node at (-1.4,0) {\tiny $n$};
\end{tikzpicture}};
\endxy
\; \sim \;
\xy
(0,0)*{
\begin{tikzpicture} [scale=.5]
\draw[double] (-1,-3) -- (-1,-2.5);
\draw[directed=1,double] (-1,2.5) -- (-1,3);
\draw[very thick, directed=.55] (-1.4,1.7) .. controls (-3,3) and (-3,-3) .. (-1.4,-1.7);
\draw [very thick, directed=.55] (-1.4,1.7) -- (-1,2.5);
\draw [very thick, directed=.55] (-1,-2.5) -- (-1.4,-1.7);
\draw [double, directed=.55] (-1.4,-1.7) -- (-1.4,1.7);
\draw[directed=.55,
very thick] (-1,-2.5) -- (1,-1.73) .. controls (1.5,-1.2) and (.6,-.5) .. (.5,0) -- (.5,0) .. controls (.6,.5) and (1.5,1.2) .. (1,1.73);
\draw[directed=.55,
very thick] (1,-1.73) .. controls (1.5,-1.2) and (2,-.5) .. (2,0) -- (2,0) .. controls (2,.5) and (1.5,1.2) .. (1,1.73) -- (-1,2.5);
\node at (-3.6,0) {\tiny $k+l$};
\node at (2.3,0) {\tiny $k$};
\node at (.2,0) {\tiny $l$};
\node at (-.6,0) {\tiny $n$};
\end{tikzpicture}};
\endxy
\; \sim \;
\xy
(0,0)*{
\begin{tikzpicture} [scale=.5]
\draw [double, directed=.55] (0,0) -- (0,.75);
\draw [very thick, directed=.55] (0,.75) -- (0,4.25);
\draw [double, directed=1] (0,4.25) -- (0,5);
\draw [very thick, directed=.55] (0,.75) -- (1,1.25) -- (1,1.75);
\draw [very thick, directed=.55] (1,1.75) -- (1,3.25);
\draw [very thick, directed=.55] (1,1.75) -- (2,2.25) -- (2,2.75) -- (1,3.25);
\draw [very thick, directed=.55] (1,3.25) -- (1,3.75) -- (0,4.25);
\node at (0,-.3) {\tiny $n$};
\node at (.7,2.5) {\tiny $l$};
\node at (2.3,2.5) {\tiny $k$};
\node [rotate=90] at (-.4,2.5) {\tiny $n-k-l$};
\end{tikzpicture}};
\endxy
\end{gather*}

The last isomorphism above is a~digon removal,
which can be found in~\cite{CKM} for example.

So,
from any upward web-tangle union $n$-strands,
we can obtain by a~succession of Reidemeister moves and equivalences a~ladder diagram.
The morphism of representation we compute by the skew Howe process has then a~diagrammatic depiction equivalent
to the skein element associated to the web-tangle we started from union the $n$-strands.

Note now that instead of pulling $n$-strands from far away,
we could have performed a~Jones--Kauf\/fman product.
Recall that the skein module may be endowed with an algebra structure by def\/ining $\alpha \ast \beta$ to be
the smoothing of the superposition of diagrams of $\alpha$ over $\beta$.
This superposition is usually assumed to be a~knotted web diagram,
meaning that the only singularities are crossings.
However,
we can allow a~singular case:
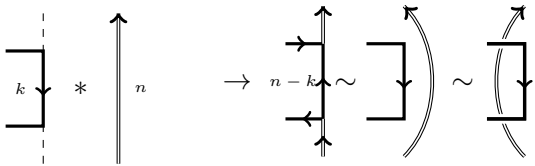
\begin{gather*}
\xy
(0,0)*{
\begin{tikzpicture}
\draw[very thick, directed=.55] (-.5,1) -- (0,1) -- (0,0) -- (-.5,0);
\node at (-.3,.5) {\tiny $k$};
\draw [dashed] (0,-.5) -- (0,1.5);
\node at (.5,.5) {$\ast$};
\draw [double, directed=1] (1,-.5) -- (1,1.5);
\node at (1.3,.5) {\tiny $n$};
\end{tikzpicture}};
\endxy
\qquad \rightarrow
\xy
(0,0)*{
\begin{tikzpicture}
\draw[very thick, directed=.55] (-.5,1) -- (0,1);
\draw[very thick, directed=.55] (0,0) -- (0,1);
\node at (-.4,.5) {\tiny $n-k$};
\draw [very thick, directed=.55] (0,0) -- (-.5,0);
\draw [double, directed=.55] (0,-.5) -- (0,0);
\draw [double, directed=1] (0,1) -- (0,1.5);
\end{tikzpicture}};
\endxy
\sim
\xy
(0,0)*{
\begin{tikzpicture}
\draw [double, directed=.99] (0,-.5) .. controls (.5,0) and (.5,1) .. (0,1.5);
\draw[very thick, directed=.55] (-.5,1) -- (0,1) -- (0,0) -- (-.5,0);
\end{tikzpicture}};
\endxy
\sim
\xy
(0,0)*{
\begin{tikzpicture}
\draw [double, directed=.99] (0,-.5) .. controls (-.5,0) and (-.5,1) .. (0,1.5);
\draw[very thick,
double,
draw= white,
double= black,
double distance= 1.25pt] (-.5,1) -- (0,1) -- (0,0) -- (-.5,0);
\draw[very thick, directed=.55] (-.5,1) -- (0,1) -- (0,0) -- (-.5,0);
\end{tikzpicture}};
\endxy
\end{gather*}
where the dashed line on the left above indicates the place
we want to put the $n$ strand: this allows not to perform any simplif\/ication on the diagram.
This re-interpretation of the process will show useful when we turn to the annular case,
where we have no free space where to put the $n$-strands before pulling them on the place they are needed.

We have seen here only the case where all the boundary of the tangle is upward.
First,
notice that this is enough for dealing with knots.
However,
as explained in~\cite{CKM},
any tangle is actually isomorphic to such an upward tangle.

\section{Af\/f\/ine extensions}

We have seen how the skew-Howe duality process,
that involves two commuting actions of $U_q(\sln)$ and $U_q(\slm)$,
helps redef\/ine Reshetikhin--Turaev $\sln$ invariants for knots and links,
that extend to invariants of knotted webs.
The f\/irst quantum group controls the invariant we are looking at,
and we therefore want to keep it unchanged.
But the second one plays the role of a~parameter related to the topology of the space we are working in.
We can thus try to modify the topology of this space.

One of the easiest extensions we can perform starting from $U_q(\slm)$ is to pass to its af\/f\/ine version $U_q(\hslm)$,
and we will show that the topological analogue of this is to close the square the knots were drawn in into an annulus.

We begin by def\/ining dif\/ferent versions of the quantum af\/f\/ine algebra $U_q(\hslm)$ that we will use here,
before turning toward easy representations of it.
We then relate this extension to knots,
and study the invariants we can deduce from it.

\subsection[Af\/f\/ine $\slm$]{Af\/f\/ine $\boldsymbol{\slm}$}

$U_q(\hslm)$ is the quantum af\/f\/ine algebra corresponding to $U_q(\slm)$,
that is the Kac--Moody algebra described by the following Dynkin diagram:
\begin{gather}
\xy
(0,0)*{
\begin{tikzpicture} [decoration={markings,
mark=at position 0.5 with {\arrow{>}}; }]
\node at (0,0) {$\circ$};
\node at (0,.3) {\tiny 1};
\node at (1,0) {$\circ$};
\node at (1,.3) {\tiny 2};
\node at (1.5,0) {$\cdots$};
\node at (2,0) {$\circ$};
\node at (2,.3) {\tiny $m-2$};
\node at (3,0) {$\circ$};
\node at (3,.3) {\tiny $m-1$};
\node at (1.5,-2) {$\circ$};
\node at (1.5,-2.3) {\tiny 0};
\draw (.1,0) -- (.9,0);
\draw (2.1,0) -- (2.9,0);
\draw (.1,-.1) -- (1.4,-1.9);
\draw (1.6,-1.9) -- (2.9,-.1);
\end{tikzpicture}};
\endxy \label{DynkinDiagAffine}
\end{gather}

Following~\cite{HongKang},
we consider the algebra $U_q(\hslm)$ as generated by Chevalley generators $E_i$,
$F_i$ and $K_i^{\pm 1}$ for $0\leq i\leq m-1$,
and extra generators $K^{\pm 1}_d$ corresponding to the null root.
The elements $E_i$,
$F_i$ and $K_i^{\pm 1}$ are subject to $\slm$ relations,
where we identify $m$ and $0$,
so that the quantum Serre relations hold for pairs $(E_0,E_1)$,
$(F_0,F_1)$,
$(E_0,E_{m-1})$ and $(F_0,F_{m-1})$:
\begin{gather*}
K_iK_i^{-1} = K_i^{-1}K_i = 1, \qquad K_iK_j = K_jK_i,
\\
K_iE_jK_i^{-1} = q^{a_{ij}} E_j, \qquad K_iF_jK_i^{-1} = q^{-a_{ij}} F_j,
\qquad
E_iF_j - F_jE_i = \delta_{ij} \frac{K_i-K_{i}^{-1}}{q-q^{-1}},
\\
E_i^2E_j-\big(q+q^{-1}\big)E_iE_jE_i+E_jE_i^2 =0\qquad \text{if} \ \ j=i\pm 1,
\\
F_i^2F_j-\big(q+q^{-1}\big)F_iF_jF_i+F_jF_i^2=0
\qquad
\text{if} \ \ j=i\pm 1,
\\
E_iE_j=E_jE_i,\qquad F_iF_j = F_jF_i
\qquad
\text{if} \ \ |i-j|>1.
\end{gather*}

Furthermore,
$K_d$ and $K_d^{-1}$ are subject to the following relations:
\begin{gather*}
K_dK_i=K_iK_d\qquad \forall\, i\in \{0,\dots,m-1\},
\qquad
K_d E_iK_d^{-1}=q^{\delta_{0,i}}E_i
\qquad
\forall\, i\in \{0,\dots,m-1\},
\\
K_d F_iK_d^{-1}=q^{-\delta_{0,i}}F_i
\qquad
\forall\, i\in \{0,\dots,m-1\}.
\end{gather*}

If we restrict to the sub-algebra generated by $E_i$,
$F_i$ and $K_i^{\pm 1}$,
we produce a~quantum group usually denoted $U'_q(\hslm)$.
A~key dif\/ference between the two versions is that the second one has f\/inite dimensional irreducible modules,
while the f\/irst one admits no non-trivial f\/inite dimensional representations.

We will also use an idempotented version of $U'_q(\hslm)$,
that we denote $\U'_q(\hslm)$,
generated by~$\onel$,
$E_i\onel$ and $F_i\onel$ with the obvious generalization of the relations of the $\slm$ case.
Weights here are $m$-tuples $\lambda=(\lambda_0,\dots,\lambda_{m-1})$,
which in our case,
with $N$ f\/ixed,
will be related to the sequences $(a_1,\dots,a_m)$ by $\lambda_i=a_{i+1}-a_i$ for $i\neq 0$ and $\lambda_0=a_1-a_m$.
Note that we have $\sum\limits \lambda_i=0$.

\subsection{Evaluation representations}

$U_q(\slm)$-representations may be extended to representations of $U'_q(\hslm)$.
The complete formulas (that require a~step through $U_q(\mathfrak{gl}_m)$) can be found in~\cite[p.~400]{ChariPressley}.
These formulas seem at f\/irst sight a~bit mysterious,
and rather than directly using them,
we choose to def\/ine dif\/ferently extensions based on the use of the braidings $\Ts$ and $\T$.
We will denote by $\rho_a$ and $\tilde{\rho}_a$ the two morphisms used in these def\/initions.
We will then investigate \textit{a posteriori} the relation with usual formulas as written in~\cite{ChariPressley}
in Proposition~\ref{prop_CP_evalrep}.

First,
recall that $C_{\Ts}$ is the automorphism of $U_q(\slm)$ (or its idempotented version) given by conjugation by $\Ts$ (that is,
$X\in U_q(\slm)\;\mapsto \Ts X \Ts^{-1}$).
We assume here that whenever we refer to a~weight $\lambda$,
a number $N$ is f\/ixed,
and we freely refer to the associated $m$-tuple $(a_1,\dots,a_m)$.

For $a$ complex number,
def\/ine $\rho_a\colon \U_q(\hslm)\mapsto \U_q(\slm)$ by:
\begin{gather*}
\rho_a(E_0\onel) = aq^{-(a_1+a_m)}C_{T''_{m-1}}\cdots C_{T''_{2}}(F_1)\onel,
\\
\rho_a(F_0\onel) = a^{-1}q^{a_1+a_m}C_{T''_{m-1}}\cdots C_{T''_{2}}(E_1)\onel.
\end{gather*}

$\rho_a$ is extended to other generators $E_i$ ($i\neq 0$) by sending $E_i\in U_q(\hslm)$ to $E_i\in U_q(\slm)$,
and similarly for $F_i$.

In terms of the graphical calculus previously def\/ined,
this can be drawn as:
\begin{gather}
E_0\onel
\mapsto
a q^{-(a_1+a_m)}
\Psi_{\Ts}\left(
\xy
(0,0)*{
\begin{tikzpicture} [scale=.5]
\draw [very thick, directed=.55] (0,-1) -- (0,2.75);
\draw [very thick, directed=1] (0,2.75) -- (0,6);
\draw [very thick, directed=.55] (1,-1) -- (1,1);
\draw [very thick] (1,1) -- (1.4,1.4);
\draw [very thick] (1.6,1.6) -- (2,2);
\draw [very thick, directed=.55] (2,2) -- (2,3);
\draw [very thick] (2,3) -- (1.6,3.4);
\draw [very thick] (1.4,3.6) -- (1,4);
\draw [very thick, directed=1] (1,4) -- (1,6);
\draw [very thick] (3,-1) -- (3.4,-.6);
\draw [very thick] (3.6,-.4)-- (4,0);
\draw [very thick, directed=.55] (4,0) -- (4,5);
\draw [very thick] (4,5) -- (3.6,5.4);
\draw [very thick,directed=1] (3.4,5.6) -- (3,6);
\draw [very thick] (4,-1) -- (2.5,.5);
\draw [very thick] (2.5,.5) -- (1,2);
\draw [very thick] (1,2) -- (1,3);
\draw [very thick, directed=1] (1,3) -- (4,6);
\draw [very thick, directed=.55] (1,2.25) -- (0,2.75);
\node at (2,-1) {$\cdots$};
\node at (3,2.5) {$\cdots$};
\node at (2,6) {$\cdots$};
\node at (0,-1.3) {\tiny $a_1$};
\node at (1,-1.3) {\tiny $a_2$};
\node at (2.9,-1.3) {\tiny $a_{m-1}$};
\node at (4.1,-1.3) {\tiny $a_m$};
\end{tikzpicture}
};
\endxy
\right)
\nonumber\\
F_0\onel
\mapsto
a^{-1} q^{a_1+a_m}
\Psi_{\Ts}\left(
\xy
(0,0)*{
\begin{tikzpicture} [scale=.5]
\draw [very thick, directed=.55] (0,-1) -- (0,2.25);
\draw [very thick, directed=1] (0,2.25) -- (0,6);
\draw [very thick, directed=.55] (1,-1) -- (1,1);
\draw [very thick] (1,1) -- (1.4,1.4);
\draw [very thick] (1.6,1.6) -- (2,2);
\draw [very thick, directed=.55] (2,2) -- (2,3);
\draw [very thick] (2,3) -- (1.6,3.4);
\draw [very thick] (1.4,3.6) -- (1,4);
\draw [very thick, directed=1] (1,4) -- (1,6);
\draw [very thick] (3,-1) -- (3.4,-.6);
\draw [very thick] (3.6,-.4)-- (4,0);
\draw [very thick, directed=.55] (4,0) -- (4,5);
\draw [very thick] (4,5) -- (3.6,5.4);
\draw [very thick,directed=1] (3.4,5.6) -- (3,6);
\draw [very thick] (4,-1) -- (2.5,.5);
\draw [very thick] (2.5,.5) -- (1,2);
\draw [very thick] (1,2) -- (1,3);
\draw [very thick, directed=1] (1,3) -- (4,6);
\draw [very thick, directed=.55] (0,2.25) -- (1,2.75);
\node at (2,-1) {$\cdots$};
\node at (3,2.5) {$\cdots$};
\node at (2,6) {$\cdots$};
\node at (0,-1.3) {\tiny $a_1$};
\node at (1,-1.3) {\tiny $a_2$};
\node at (2.9,-1.3) {\tiny $a_{m-1}$};
\node at (4.1,-1.3) {\tiny $a_m$};
\end{tikzpicture}
};
\endxy
\right)\label{evalRep}
\end{gather}

The above depiction corresponds to letting $\U_q(\hslm)$ act on a~tensor product of minuscule representations of $U_q(\slm)$ via the map $\rho_a$.
This diagrammatic def\/inition gives us an easy way to reprove (in Proposition~\ref{PropEvalRep})
in a~diagrammatic fashion the results of~\cite[Proposition~12.2.10]{ChariPressley} applied
to these particular representations ${\bigwedge}_q(\C^n\otimes \C^m)$.
However,
since all fundamental $U_q(\slm)$ representations ${\bigwedge}^k_q(\C^m)$ appear as summands of a~${\bigwedge}^N(\C^n\otimes\C^m)$,
using the coproduct,
we can recover from this construction all evaluation representations associated to f\/inite-dimensional representations.

The relations in Propositions~\ref{propCTs} and~\ref{propCTp} together with the fact that the $\Ts$'s provide a~brai\-ding
give the next useful diagrammatic relations.
In order to simplify the notation,
we omit here to write $\Psi_{\Ts}$ in all relations below:
\begin{gather}\label{diagRel1}
\xy
(0,0)*{
\begin{tikzpicture} [scale=.5]
\draw [very thick] (0,-1) -- (0,1);
\draw [very thick] (0,1) -- (2,3);
\draw [very thick, directed=1] (2,3) -- (2,4);
\draw [very thick, directed=.55] (1,-1) -- (1,1);
\draw [very thick] (1,1) -- (.6,1.4);
\draw [very thick] (.4,1.6) -- (0,2);
\draw [very thick, directed=1] (0,2) -- (0,4);
\draw [very thick] (2,-1) -- (2,0);
\draw [very thick, directed=.55] (2,0) -- (2,2);
\draw [very thick] (2,2) -- (1.6,2.4);
\draw [very thick] (1.4,2.6) -- (1,3);
\draw [very thick, directed=1] (1,3) -- (1,4);
\draw [very thick, directed=.55] (2,.25) -- (1,.75);
\end{tikzpicture}
};
\endxy
=
\xy
(0,0)*{
\begin{tikzpicture} [scale=.5]
\draw [very thick] (0,-1) -- (0,0);
\draw [very thick] (0,0) -- (2,2);
\draw [very thick,directed=1] (2,2) -- (2,4);
\draw [very thick, directed=.55] (1,-1) -- (1,0);
\draw [very thick] (1,0) -- (.6,.4);
\draw [very thick] (.4,.6) -- (0,1);
\draw [very thick,directed=1] (0,1) -- (0,4);
\draw [very thick, directed=.55] (2,-1) -- (2,1);
\draw [very thick] (2,1) -- (1.6,1.4);
\draw [very thick] (1.4,1.6) -- (1,2);
\draw [very thick,directed=1] (1,2) -- (1,4);
\draw [very thick, directed=.55] (1,2.25) -- (0,2.75);
\end{tikzpicture}
};
\endxy
\\
\label{diagRel2}
\xy
(0,0)*{
\begin{tikzpicture} [scale=.5]
\draw [very thick, directed=.55] (0,0) -- (0,2);
\draw [very thick] (0,2) -- (1,3);
\draw [very thick, directed=1] (1,3) -- (1,4);
\draw [very thick, directed=.55] (1,0) -- (1,2);
\draw [very thick] (1,2) -- (.6,2.4);
\draw [very thick] (.4,2.6) -- (0,3);
\draw [very thick, directed=1] (0,3) -- (0,4);
\draw [very thick, directed=.55] (0,1.25) -- (1,1.75);
\node at (0,-.3) {\text{\tiny $a_i$}};
\node at (1,-.3) {\text{\tiny $a_{i+1}$}};
\node at (0,4.3) {};
\node at (1,4.3) {};
\end{tikzpicture}
};
\endxy
=
-q^{a_i-a_{i+1}}\xy
(0,0)*{
\begin{tikzpicture} [scale=.5]
\draw [very thick, directed=.55] (0,0) -- (0,1);
\draw [very thick] (0,1) -- (1,2);
\draw [very thick, directed=1] (1,2) -- (1,4);
\draw [very thick, directed=.55] (1,0) -- (1,1);
\draw [very thick] (1,1) -- (.6,1.4);
\draw [very thick] (.4,1.6) -- (0,2);
\draw [very thick, directed=1] (0,2) -- (0,4);
\draw [very thick, directed=.55] (1,2.25) -- (0,2.75);
\node at (0,-.3) {\text{\tiny $a_i$}};
\node at (1,-.3) {\text{\tiny $a_{i+1}$}};
\node at (0,4.3) {};
\node at (1,4.3) {};
\end{tikzpicture}
};
\endxy
, \qquad \xy
(0,0)*{
\begin{tikzpicture} [scale=.5]
\draw [very thick, directed=.55] (0,0) -- (0,2);
\draw [very thick] (0,2) -- (1,3);
\draw [very thick, directed=1] (1,3) -- (1,4);
\draw [very thick, directed=.55] (1,0) -- (1,2);
\draw [very thick] (1,2) -- (.6,2.4);
\draw [very thick] (.4,2.6) -- (0,3);
\draw [very thick, directed=1] (0,3) -- (0,4);
\draw [very thick, directed=.55] (1,1.25) -- (0,1.75);
\node at (0,4.3) {\text{\tiny $a_i$}};
\node at (1,4.3) {\text{\tiny $a_{i+1}$}};
\node at (0,-.3) {};
\node at (1,-.3) {};
\end{tikzpicture}
};
\endxy
=
-q^{a_i-a_{i+1}}\xy
(0,0)*{
\begin{tikzpicture} [scale=.5]
\draw [very thick, directed=.55] (0,0) -- (0,1);
\draw [very thick] (0,1) -- (1,2);
\draw [very thick, directed=1] (1,2) -- (1,4);
\draw [very thick, directed=.55] (1,0) -- (1,1);
\draw [very thick] (1,1) -- (.6,1.4);
\draw [very thick] (.4,1.6) -- (0,2);
\draw [very thick, directed=1] (0,2) -- (0,4);
\draw [very thick, directed=.55] (0,2.25) -- (1,2.75);
\node at (0,4.3) {\text{\tiny $a_i$}};
\node at (1,4.3) {\text{\tiny $a_{i+1}$}};
\node at (0,-.3) {};
\node at (1,-.3) {};
\end{tikzpicture}
};
\endxy
\\
\label{diagRel3}
\xy
(0,0)*{
\begin{tikzpicture} [scale=.5]
\draw [very thick, directed=.55] (0,0) -- (0,2);
\draw [very thick] (0,2) -- (.4,2.4);
\draw [very thick] (.6,2.6) -- (1,3);
\draw [very thick, directed=1] (1,3) -- (1,4);
\draw [very thick, directed=.55] (1,0) -- (1,2);
\draw [very thick] (1,2) -- (0,3);
\draw [very thick, directed=1] (0,3) -- (0,4);
\draw [very thick, directed=.55] (0,1.25) -- (1,1.75);
\node at (0,-.3) {};
\node at (1,-.3) {};
\node at (0,4.3) {\text{\tiny $a_i$}};
\node at (1,4.3) {\text{\tiny $a_{i+1}$}};
\end{tikzpicture}
};
\endxy
=
-q^{a_i-a_{i+1}}
\xy
(0,0)*{
\begin{tikzpicture} [scale=.5]
\draw [very thick, directed=.55] (0,0) -- (0,1);
\draw [very thick] (0,1) -- (.4,1.4);
\draw [very thick] (.6,1.6) -- (1,2);
\draw [very thick, directed=1] (1,2) -- (1,4);
\draw [very thick, directed=.55] (1,0) -- (1,1);
\draw [very thick] (1,1) -- (0,2);
\draw [very thick, directed=1] (0,2) -- (0,4);
\draw [very thick, directed=.55] (1,2.25) -- (0,2.75);
\node at (0,-.3) {};
\node at (1,-.3) {};
\node at (0,4.3) {\text{\tiny $a_i$}};
\node at (1,4.3) {\text{\tiny $a_{i+1}$}};
\end{tikzpicture}
};
\endxy
, \qquad \xy
(0,0)*{
\begin{tikzpicture} [scale=.5]
\draw [very thick, directed=.55] (0,0) -- (0,2);
\draw [very thick] (0,2) -- (.4,2.4);
\draw [very thick] (.6,2.6) -- (1,3);
\draw [very thick, directed=1] (1,3) -- (1,4);
\draw [very thick, directed=.55] (1,0) -- (1,2);
\draw [very thick] (1,2) -- (0,3);
\draw [very thick, directed=1] (0,3) -- (0,4);
\draw [very thick, directed=.55] (1,1.25) -- (0,1.75);
\node at (0,4.3) {};
\node at (1,4.3) {};
\node at (0,-.3) {\text{\tiny $a_i$}};
\node at (1,-.3) {\text{\tiny $a_{i+1}$}};
\end{tikzpicture}
};
\endxy
=
-q^{a_i-a_{i+1}}
\xy
(0,0)*{
\begin{tikzpicture} [scale=.5]
\draw [very thick, directed=.55] (0,0) -- (0,1);
\draw [very thick] (0,1) -- (.4,1.4);
\draw [very thick] (.6,1.6) -- (1,2);
\draw [very thick, directed=1] (1,2) -- (1,4);
\draw [very thick, directed=.55] (1,0) -- (1,1);
\draw [very thick] (1,1) -- (0,2);
\draw [very thick, directed=1] (0,2) -- (0,4);
\draw [very thick, directed=.55] (0,2.25) -- (1,2.75);
\node at (0,4.3) {};
\node at (1,4.3) {};
\node at (0,-.3) {\text{\tiny $a_i$}};
\node at (1,-.3) {\text{\tiny $a_{i+1}$}};
\end{tikzpicture}
};
\endxy
\end{gather}

\begin{Proposition} \label{PropEvalRep}
The action of $\U_q(\slm)$ on ${\bigwedge}_q(\C^n\otimes \C^m)$ extends via $\rho_a$ to an action of $\U'_q(\hat{\slm})$.
\end{Proposition}
\begin{proof}
The proof consists in checking the def\/ining relations of $\U'_q(\hat{\slm})$ that involve $E_0$ or $F_0$,
that is
\begin{gather*}
[E_0,F_0] \onel = [\lambda_0]\onel,
\qquad
[E_0,F_j] \onel=0,
\qquad
[F_0,E_j]\onel=0,
\\
E_0^2E_j\onel - \big(q+q^{-1}\big)E_0E_jE_0\onel +E_jE_0^2\onel=0 \qquad \text{if} \ \ j=1\quad \text{or}\quad j=m-1,
\\
E_j^2E_0\onel - \big(q+q^{-1}\big)E_jE_0E_j\onel +E_0E_j^2\onel=0 \qquad \text{if} \ \ j=1\quad \text{or}\quad j=m-1,
\\
F_0^2F_j\onel - \big(q+q^{-1}\big)F_0F_jF_0\onel +F_jF_0^2\onel=0 \qquad \text{if} \ \ j=1\quad \text{or}\quad j=m-1,
\\
F_j^2F_0\onel - \big(q+q^{-1}\big)F_jF_0F_j\onel +F_0F_j^2\onel=0 \qquad \text{if} \ \ j=1\quad \text{or}\quad j=m-1,
\\
E_0E_j=E_jE_0,
\qquad F_0F_j=F_jF_0 \qquad \text{if} \ \  j\neq1\quad \text{and} \quad j\neq m-1.
\end{gather*}

By symmetry between the $E's$ and the $F$'s,
it is enough to check only relations involving $E_0$ and $E_j$ or $F_j$,
plus the f\/irst one.
The proof mostly relies on a~straightforward use of relations \eqref{diagRel1},
\eqref{diagRel2},
\eqref{diagRel3},
and the braiding relation.
We present some of them below,
where we identify an element of $\U'_q(\slm)$ with its image under the representation
\begin{gather*}
[E_0,F_0]\onel = E_0F_0\onel - F_0E_0\onel = q^{0}[C_{T''_{w}}(F_1),C_{T''_{w}}(E_1)]\onel=C_{T''_{w}}([F_1,E_1]\onell{w\cdot \lambda})
\\
\phantom{[E_0,F_0]\onel}
= C_{T''_{w}}\big({-}\big[\big(w^{-1}\cdot \lambda\big)_1\big]\big)\onel = -(-[\lambda_0])\onel =[\lambda_0]\onel.
\end{gather*}

Under the braiding relation,
we can reduce $[E_0,F_1]$ to (we again omit to write $\Psi_{\Ts}$):
\begin{gather*}
[E_0,F_1]\onel =
aq^{-a-1-c}
 \xy
(0,0)*{
\begin{tikzpicture} [scale=.5]
\draw [very thick, directed=.55] (0,0) -- (0,1);
\draw [very thick] (0,1) -- (0,5);
\draw [very thick, directed=1] (0,5) -- (0,6);
\draw [very thick, directed=.55] (1,0) -- (1,1);
\draw [very thick] (1,1) -- (1,2) -- (1.4,2.4);
\draw [very thick] (1.6,2.6) -- (2,3) -- (2,4) -- (1.6,4.4);
\draw [very thick] (1.4,4.6) -- (1,5);
\draw [very thick, directed=1] (1,5) -- (1,6);
\draw [very thick, directed=.55] (2,0) -- (2,1);
\draw [very thick] (2,1) -- (2,2) -- (1,3) -- (1,4) -- (2,5);
\draw [very thick, directed=1] (2,5) -- (2,6);
\draw [very thick, directed=.55] (1,1.25) -- (0,1.75);
\draw [very thick, directed=.55] (1,3.25) -- (0,3.75);
\node at (0,7) {};
\node at (1,7) {};
\node at (2,7) {};
\node at (0,-1) {\text{\footnotesize{$a$}}};
\node at (1,-1) {\text{\footnotesize{$b$}}};
\node at (2,-1) {\text{\footnotesize{$c$}}};
\end{tikzpicture}
};
\endxy
-
aq^{-a-c}\xy
(0,0)*{
\begin{tikzpicture} [scale=.5]
\draw [very thick, directed=.55] (0,0) -- (0,1);
\draw [very thick] (0,1) -- (0,5);
\draw [very thick, directed=1] (0,5) -- (0,6);
\draw [very thick, directed=.55] (1,0) -- (1,1);
\draw [very thick] (1,1) -- (1.4,1.4);
\draw [very thick] (1.6,1.6) -- (2,2) -- (2,3) -- (1.6,3.4);
\draw [very thick] (1.4,3.6) -- (1,4) -- (1,5);
\draw [very thick, directed=1] (1,5) -- (1,6);
\draw [very thick, directed=.55] (2,0) -- (2,1);
\draw [very thick] (2,1) -- (1,2) -- (1,3) -- (2,4) -- (2,5);
\draw [very thick, directed=1] (2,5) -- (2,6);
\draw [very thick, directed=.55] (1,2.25) -- (0,2.75);
\draw [very thick, directed=.55] (1,4.25) -- (0,4.75);
\node at (0,7) {};
\node at (1,7) {};
\node at (2,7) {};
\node at (0,-1) {\text{\footnotesize{$a$}}};
\node at (1,-1) {\text{\footnotesize{$b$}}};
\node at (2,-1) {\text{\footnotesize{$c$}}};
\end{tikzpicture}
};
\endxy
\end{gather*}
which equals
\begin{gather*}
=
aq^{-a-1-c}
\xy
(0,0)*{
\begin{tikzpicture} [scale=.5]
\draw [very thick, directed=.55] (0,0) -- (0,1);
\draw [very thick] (0,1) -- (0,2) -- (1,3) -- (1,4) -- (0,5);
\draw [very thick, directed=1] (0,5) -- (0,6);
\draw [very thick, directed=.55] (1,0) -- (1,1);
\draw [very thick] (1,1) -- (1,2) -- (.6,2.4);
\draw [very thick] (.4,2.6) -- (0,3) -- (0,4) -- (.4,4.4);
\draw [very thick] (.6,4.6) -- (1,5);
\draw [very thick, directed=1] (1,5) -- (1,6);
\draw [very thick, directed=.55] (2,0) -- (2,1);
\draw [very thick] (2,1) -- (2,5);
\draw [very thick, directed=1] (2,5) -- (2,6);
\draw [very thick, directed=.55] (1,1.25) -- (0,1.75);
\draw [very thick, directed=.55] (2,3.25) -- (1,3.75);
\node at (0,7) {};
\node at (1,7) {};
\node at (2,7) {};
\node at (0,-1) {\text{\footnotesize{$a$}}};
\node at (1,-1) {\text{\footnotesize{$b$}}};
\node at (2,-1) {\text{\footnotesize{$c$}}};
\end{tikzpicture}
};
\endxy
-
aq^{-a-c}\xy
(0,0)*{
\begin{tikzpicture} [scale=.5]
\draw [very thick, directed=.55] (0,0) -- (0,1);
\draw [very thick] (0,1) -- (1,2) -- (1,3) -- (0,4) -- (0,5);
\draw [very thick, directed=1] (0,5) -- (0,6);
\draw [very thick, directed=.55] (1,0) -- (1,1);
\draw [very thick] (1,1) -- (.6,1.4);
\draw [very thick] (.4,1.6) -- (0,2) -- (0,3) -- (.4,3.4);
\draw [very thick] (.6,3.6) -- (1,4) -- (1,5);
\draw [very thick, directed=1] (1,5) -- (1,6);
\draw [very thick, directed=.55] (2,0) -- (2,1);
\draw [very thick] (2,1) -- (2,5);
\draw [very thick, directed=1] (2,5) -- (2,6);
\draw [very thick, directed=.55] (2,2.25) -- (1,2.75);
\draw [very thick, directed=.55] (1,4.25) -- (0,4.75);
\node at (0,7) {};
\node at (1,7) {};
\node at (2,7) {};
\node at (0,-1) {\text{\footnotesize{$a$}}};
\node at (1,-1) {\text{\footnotesize{$b$}}};
\node at (2,-1) {\text{\footnotesize{$c$}}};
\end{tikzpicture}
};
\endxy
\\
= aq^{-a-1-c} \big(-q^{-a-1+b-1}\big)\xy
(0,0)*{
\begin{tikzpicture} [scale=.5]
\draw [very thick, directed=.55] (0,0) -- (0,1);
\draw [very thick] (0,1) -- (1,2) -- (1,4) -- (0,5);
\draw [very thick, directed=1] (0,5) -- (0,6);
\draw [very thick, directed=.55] (1,0) -- (1,1);
\draw [very thick] (1,1) -- (.6,1.4);
\draw [very thick] (.4,1.6) -- (0,2) -- (0,4) -- (.4,4.4);
\draw [very thick] (.6,4.6) -- (1,5);
\draw [very thick, directed=1] (1,5) -- (1,6);
\draw [very thick, directed=.55] (2,0) -- (2,1);
\draw [very thick] (2,1) -- (2,5);
\draw [very thick, directed=1] (2,5) -- (2,6);
\draw [very thick, directed=.55] (2,3.25) -- (1,3.75);
\draw [very thick, directed=.55] (0,2.25) -- (1,2.75);
\node at (0,7) {};
\node at (1,7) {};
\node at (2,7) {};
\node at (0,-1) {\text{\footnotesize{$a$}}};
\node at (1,-1) {\text{\footnotesize{$b$}}};
\node at (2,-1) {\text{\footnotesize{$c$}}};
\end{tikzpicture}
};
\endxy
-
aq^{-a-c} \big({-}q^{-a-2+b-1}\big)\xy
(0,0)*{
\begin{tikzpicture} [scale=.5]
\draw [very thick, directed=.55] (0,0) -- (0,1);
\draw [very thick] (0,1) -- (1,2) -- (1,4) -- (0,5);
\draw [very thick, directed=1] (0,5) -- (0,6);
\draw [very thick, directed=.55] (1,0) -- (1,1);
\draw [very thick] (1,1) -- (.6,1.4);
\draw [very thick] (.4,1.6) -- (0,2) -- (0,4) -- (.4,4.4);
\draw [very thick] (.6,4.6) -- (1,5);
\draw [very thick, directed=1] (1,5) -- (1,6);
\draw [very thick, directed=.55] (2,0) -- (2,1);
\draw [very thick] (2,1) -- (2,5);
\draw [very thick, directed=1] (2,5) -- (2,6);
\draw [very thick, directed=.55] (2,2.25) -- (1,2.75);
\draw [very thick, directed=.55] (0,3.25) -- (1,3.75);
\node at (0,7) {};
\node at (1,7) {};
\node at (2,7) {};
\node at (0,-1) {\text{\footnotesize{$a$}}};
\node at (1,-1) {\text{\footnotesize{$b$}}};
\node at (2,-1) {\text{\footnotesize{$c$}}};
\end{tikzpicture}
};
\endxy
\end{gather*}
We then obtain:
\begin{gather*}
-aq^{-2a+b-c-3}
\left(
\xy
(0,0)*{
\begin{tikzpicture} [scale=.5]
\draw [very thick, directed=.55] (0,0) -- (0,1);
\draw [very thick] (0,1) -- (1,2) -- (1,4) -- (0,5);
\draw [very thick, directed=1] (0,5) -- (0,6);
\draw [very thick, directed=.55] (1,0) -- (1,1);
\draw [very thick] (1,1) -- (.6,1.4);
\draw [very thick] (.4,1.6) -- (0,2) -- (0,4) -- (.4,4.4);
\draw [very thick] (.6,4.6) -- (1,5);
\draw [very thick, directed=1] (1,5) -- (1,6);
\draw [very thick, directed=.55] (2,0) -- (2,1);
\draw [very thick] (2,1) -- (2,5);
\draw [very thick, directed=1] (2,5) -- (2,6);
\draw [very thick, directed=.55] (2,2.25) -- (1,2.75);
\draw [very thick, directed=.55] (0,3.25) -- (1,3.75);
\node at (0,7) {};
\node at (1,7) {};
\node at (2,7) {};
\node at (0,-1) {\text{\footnotesize{$a$}}};
\node at (1,-1) {\text{\footnotesize{$b$}}};
\node at (2,-1) {\text{\footnotesize{$c$}}};
\end{tikzpicture}
};
\endxy
- \xy
(0,0)*{
\begin{tikzpicture} [scale=.5]
\draw [very thick, directed=.55] (0,0) -- (0,1);
\draw [very thick] (0,1) -- (1,2) -- (1,4) -- (0,5);
\draw [very thick, directed=1] (0,5) -- (0,6);
\draw [very thick, directed=.55] (1,0) -- (1,1);
\draw [very thick] (1,1) -- (.6,1.4);
\draw [very thick] (.4,1.6) -- (0,2) -- (0,4) -- (.4,4.4);
\draw [very thick] (.6,4.6) -- (1,5);
\draw [very thick, directed=1] (1,5) -- (1,6);
\draw [very thick, directed=.55] (2,0) -- (2,1);
\draw [very thick] (2,1) -- (2,5);
\draw [very thick, directed=1] (2,5) -- (2,6);
\draw [very thick, directed=.55] (2,3.25) -- (1,3.75);
\draw [very thick, directed=.55] (0,2.25) -- (1,2.75);
\node at (0,7) {};
\node at (1,7) {};
\node at (2,7) {};
\node at (0,-1) {\text{\footnotesize{$a$}}};
\node at (1,-1) {\text{\footnotesize{$b$}}};
\node at (2,-1) {\text{\footnotesize{$c$}}};
\end{tikzpicture}
};
\endxy
\right) =0
\end{gather*}
using the $U_q(\slm)$ relation $[E_1,F_2]=0$.
We now turn to
\begin{gather*}
E_0^2E_1-\big(q+q^{-1}\big)E_0E_1E_0+E_1E_0^2=a^2 q^{-2a-2c+2}
\xy
(0,0)*{
\begin{tikzpicture} [scale=.5]
\draw [very thick, directed=.55] (0,0) -- (0,1);
\draw [very thick] (0,1) -- (0,8);
\draw [very thick, directed=1] (0,8) -- (0,9);
\draw [very thick, directed=.55] (1,0) -- (1,1);
\draw [very thick] (1,1) -- (1,2) -- (1.4,2.4);
\draw [very thick] (1.6,2.6) -- (2,3) -- (2,4) -- (1.6,4.4);
\draw [very thick] (1.4,4.6) -- (1,5) -- (1.4,5.4);
\draw [very thick] (1.6,5.6) -- (2,6) -- (2,7) -- (1.6,7.4);
\draw [very thick] (1.4,7.6) -- (1,8);
\draw [very thick, directed=1] (1,8) -- (1,9);
\draw [very thick, directed=.55] (2,0) -- (2,1);
\draw [very thick] (2,1) -- (2,2) -- (1,3) -- (1,4) -- (2,5) -- (1,6) -- (1,7) -- (2,8);
\draw [very thick, directed=1] (2,8) -- (2,9);
\draw [very thick, directed=.55] (0,1.25) -- (1,1.75);
\draw [very thick, directed=.55] (1,3.25) -- (0,3.75);
\draw [very thick, directed=.55] (1,6.25) -- (0,6.75);
\node at (0,10) {};
\node at (1,10) {};
\node at (2,10) {};
\node at (0,-1) {\text{\footnotesize{$a$}}};
\node at (1,-1) {\text{\footnotesize{$b$}}};
\node at (2,-1) {\text{\footnotesize{$c$}}};
\end{tikzpicture}
};
\endxy
\\
{}-a^2 \big(q+q^{-1}\big) q^{-2a-2c+1}
\xy
(0,0)*{
\begin{tikzpicture} [scale=.5]
\draw [very thick, directed=.55] (0,0) -- (0,1);
\draw [very thick] (0,1) -- (0,8);
\draw [very thick, directed=1] (0,8) -- (0,9);
\draw [very thick, directed=.55] (1,0) -- (1,1);
\draw [very thick] (1,1) -- (1.4,1.4);
\draw [very thick] (1.6,1.6) -- (2,2) -- (2,3) -- (1.6,3.4);
\draw [very thick] (1.4,3.6) -- (1,4) -- (1,5) -- (1.4,5.4);
\draw [very thick] (1.6,5.6) -- (2,6) -- (2,7) -- (1.6,7.4);
\draw [very thick] (1.4,7.6) -- (1,8);
\draw [very thick, directed=1] (1,8) -- (1,9);
\draw [very thick, directed=.55] (2,0) -- (2,1);
\draw [very thick] (2,1) -- (1,2) -- (1,3) -- (2,4) -- (2,5) -- (1,6) -- (1,7) -- (2,8);
\draw [very thick, directed=1] (2,8) -- (2,9);
\draw [very thick, directed=.55] (1,2.25) -- (0,2.75);
\draw [very thick, directed=.55] (0,4.25) -- (1,4.75);
\draw [very thick, directed=.55] (1,6.25) -- (0,6.75);
\node at (0,10) {};
\node at (1,10) {};
\node at (2,10) {};
\node at (0,-1) {\text{\footnotesize{$a$}}};
\node at (1,-1) {\text{\footnotesize{$b$}}};
\node at (2,-1) {\text{\footnotesize{$c$}}};
\end{tikzpicture}
};
\endxy
\ +a^2 q^{-2a-2c} \xy
(0,0)*{
\begin{tikzpicture} [scale=.5]
\draw [very thick, directed=.55] (0,0) -- (0,1);
\draw [very thick] (0,1) -- (0,8);
\draw [very thick, directed=1] (0,8) -- (0,9);
\draw [very thick, directed=.55] (1,0) -- (1,1);
\draw [very thick] (1,1) -- (1.4,1.4);
\draw [very thick] (1.6,1.6) -- (2,2) -- (2,3) -- (1.6,3.4);
\draw [very thick] (1.4,3.6) -- (1,4) -- (1.4,4.4);
\draw [very thick] (1.6,4.6) -- (2,5) -- (2,6) -- (1.6,6.4);
\draw [very thick] (1.4,6.6) -- (1,7) -- (1,8);
\draw [very thick, directed=1] (1,8) -- (1,9);
\draw [very thick, directed=.55] (2,0) -- (2,1);
\draw [very thick] (2,1) -- (1,2) -- (1,3) -- (2,4) -- (1,5) -- (1,6) -- (2,7) -- (2,8);
\draw [very thick, directed=1] (2,8) -- (2,9);
\draw [very thick, directed=.55] (1,2.25) -- (0,2.75);
\draw [very thick, directed=.55] (1,5.25) -- (0,5.75);
\draw [very thick, directed=.55] (0,7.25) -- (1,7.75);
\node at (0,10) {};
\node at (1,10) {};
\node at (2,10) {};
\node at (0,-1) {\text{\footnotesize{$a$}}};
\node at (1,-1) {\text{\footnotesize{$b$}}};
\node at (2,-1) {\text{\footnotesize{$c$}}};
\end{tikzpicture}
};
\endxy
\\
=   a^2 q^{-2a-2c+2}
\xy
(0,0)*{
\begin{tikzpicture} [scale=.5]
\draw [very thick, directed=.55] (0,-1) -- (0,1);
\draw [very thick] (0,1) -- (0,2) -- (1,3) -- (1,5) -- (0,6);
\draw [very thick, directed=1] (0,6) -- (0,8);
\draw [very thick, directed=.55] (1,-1) -- (1,1);
\draw [very thick] (1,1) -- (1,2) -- (.6,2.4);
\draw [very thick] (.4,2.6) -- (0,3) -- (0,5) -- (.4,5.4);
\draw [very thick] (.6,5.6) -- (1,6);
\draw [very thick, directed=1] (1,6) -- (1,8);
\draw [very thick, directed=.55] (2,-1) -- (2,1);
\draw [very thick] (2,1) -- (2,6);
\draw [very thick, directed=1] (2,6) -- (2,8);
\draw [very thick, directed=.55] (0,1.25) -- (1,1.75);
\draw [very thick, directed=.55] (2,3.25) -- (1,3.75);
\draw [very thick, directed=.55] (2,4.25) -- (1,4.75);
\node at (0,9) {};
\node at (1,9) {};
\node at (2,9) {};
\node at (0,-2) {\text{\footnotesize{$a$}}};
\node at (1,-2) {\text{\footnotesize{$b$}}};
\node at (2,-2) {\text{\footnotesize{$c$}}};
\end{tikzpicture}
};
\endxy \; -\big(q+q^{-1}\big) a^2 q^{-2a-2c+1}
\xy
(0,0)*{
\begin{tikzpicture} [scale=.5]
\draw [very thick, directed=.55] (0,0) -- (0,1);
\draw [very thick] (0,1) -- (1,2) -- (1,3) -- (0,4) -- (0,5) -- (1,6) -- (1,7) -- (0,8);
\draw [very thick, directed=1] (0,8) -- (0,9);
\draw [very thick, directed=.55] (1,0) -- (1,1);
\draw [very thick] (1,1) -- (.6,1.4);
\draw [very thick] (.4,1.6) -- (0,2) -- (0,3) -- (.4,3.4);
\draw [very thick] (.6,3.6) -- (1,4) -- (1,5) -- (.6,5.4);
\draw [very thick] (.4,5.6) -- (0,6) -- (0,7) -- (.4,7.4);
\draw [very thick] (.6,7.6) -- (1,8);
\draw [very thick, directed=1] (1,8) -- (1,9);
\draw [very thick, directed=.55] (2,0) -- (2,1);
\draw [very thick] (2,1) -- (2,8);
\draw [very thick, directed=1] (2,8) -- (2,9);
\draw [very thick, directed=.55] (2,2.25) -- (1,2.75);
\draw [very thick, directed=.55] (0,4.25) -- (1,4.75);
\draw [very thick, directed=.55] (2,6.25) -- (1,6.75);
\node at (0,10) {};
\node at (1,10) {};
\node at (2,10) {};
\node at (0,-1) {\text{\footnotesize{$a$}}};
\node at (1,-1) {\text{\footnotesize{$b$}}};
\node at (2,-1) {\text{\footnotesize{$c$}}};
\end{tikzpicture}
};
\endxy
\  +a^2 q^{-2a-2c}
\xy
(0,0)*{
\begin{tikzpicture} [scale=.5]
\draw [very thick, directed=.55] (0,-1) -- (0,1);
\draw [very thick] (0,1) -- (1,2) -- (1,4) -- (0,5) -- (0,6);
\draw [very thick, directed=1] (0,6) -- (0,8);
\draw [very thick, directed=.55] (1,-1) -- (1,1);
\draw [very thick] (1,1) -- (.6,1.4);
\draw [very thick] (.4,1.6) -- (0,2) -- (0,4) -- (.4,4.4);
\draw [very thick] (.6,4.6) -- (1,5) -- (1,6);
\draw [very thick, directed=1] (1,6) -- (1,8);
\draw [very thick, directed=.55] (2,-1) -- (2,1);
\draw [very thick] (2,1)-- (2,6);
\draw [very thick, directed=1] (2,6) -- (2,8);
\draw [very thick, directed=.55] (2,2.25) -- (1,2.75);
\draw [very thick, directed=.55] (2,3.25) -- (1,3.75);
\draw [very thick, directed=.55] (0,5.25) -- (1,5.75);
\node at (0,9) {};
\node at (1,9) {};
\node at (2,9) {};
\node at (0,-2) {\text{\footnotesize{$a$}}};
\node at (1,-2) {\text{\footnotesize{$b$}}};
\node at (2,-2) {\text{\footnotesize{$c$}}};
\end{tikzpicture}
};
\endxy
\\
= -a^2q^{-a-b-2c+2}
\xy
(0,0)*{
\begin{tikzpicture} [scale=.5]
\draw [very thick, directed=.55] (0,0) -- (0,1);
\draw [very thick] (0,1) -- (1,2) -- (1,5) -- (0,6);
\draw [very thick, directed=1] (0,6) -- (0,7);
\draw [very thick, directed=.55] (1,0) -- (1,1);
\draw [very thick] (1,1) -- (.6,1.4);
\draw [very thick] (.4,1.6) -- (0,2) -- (0,5) -- (.4,5.4);
\draw [very thick] (.6,5.6) -- (1,6);
\draw [very thick, directed=1] (1,6) -- (1,7);
\draw [very thick, directed=.55] (2,0) -- (2,1);
\draw [very thick] (2,1) -- (2,6);
\draw [very thick, directed=1] (2,6) -- (2,7);
\draw [very thick, directed=.55] (1,2.25) -- (0,2.75);
\draw [very thick, directed=.55] (2,3.25) -- (1,3.75);
\draw [very thick, directed=.55] (2,4.25) -- (1,4.75);
\node at (0,8) {};
\node at (1,8) {};
\node at (2,8) {};
\node at (0,-1) {\text{\footnotesize{$a$}}};
\node at (1,-1) {\text{\footnotesize{$b$}}};
\node at (2,-1) {\text{\footnotesize{$c$}}};
\end{tikzpicture}
};
\endxy
\  + a^2 \big(q+q^{-1}\big) q^{-a-b-2c+2}
\xy
(0,0)*{
\begin{tikzpicture} [scale=.5]
\draw [very thick, directed=.55] (0,0) -- (0,1);
\draw [very thick] (0,1) -- (1,2) -- (1,5) -- (0,6);
\draw [very thick, directed=1] (0,6) -- (0,7);
\draw [very thick, directed=.55] (1,0) -- (1,1);
\draw [very thick] (1,1) -- (.6,1.4);
\draw [very thick] (.4,1.6) -- (0,2) -- (0,5) -- (.4,5.4);
\draw [very thick] (.6,5.6) -- (1,6);
\draw [very thick, directed=1] (1,6) -- (1,7);
\draw [very thick, directed=.55] (2,0) -- (2,1);
\draw [very thick] (2,1) -- (2,6);
\draw [very thick, directed=1] (2,6) -- (2,7);
\draw [very thick, directed=.55] (2,2.25) -- (1,2.75);
\draw [very thick, directed=.55] (1,3.25) -- (0,3.75);
\draw [very thick, directed=.55] (2,4.25) -- (1,4.75);
\node at (0,8) {};
\node at (1,8) {};
\node at (2,8) {};
\node at (0,-1) {\text{\footnotesize{$a$}}};
\node at (1,-1) {\text{\footnotesize{$b$}}};
\node at (2,-1) {\text{\footnotesize{$c$}}};
\end{tikzpicture}
};
\endxy
\  - a^2 q^{-a-b-2c+2}
\xy
(0,0)*{
\begin{tikzpicture} [scale=.5]
\draw [very thick, directed=.55] (0,0) -- (0,1);
\draw [very thick] (0,1) -- (1,2) -- (1,5) -- (0,6);
\draw [very thick, directed=1] (0,6) -- (0,7);
\draw [very thick, directed=.55] (1,0) -- (1,1);
\draw [very thick] (1,1) -- (.6,1.4);
\draw [very thick] (.4,1.6) -- (0,2) -- (0,5) -- (.4,5.4);
\draw [very thick] (.6,5.6) -- (1,6);
\draw [very thick, directed=1] (1,6) -- (1,7);
\draw [very thick, directed=.55] (2,0) -- (2,1);
\draw [very thick] (2,1) -- (2,6);
\draw [very thick, directed=1] (2,6) -- (2,7);
\draw [very thick, directed=.55] (2,2.25) -- (1,2.75);
\draw [very thick, directed=.55] (2,3.25) -- (1,3.75);
\draw [very thick, directed=.55] (1,4.25) -- (0,4.75);
\node at (0,8) {};
\node at (1,8) {};
\node at (2,8) {};
\node at (0,-1) {\text{\footnotesize{$a$}}};
\node at (1,-1) {\text{\footnotesize{$b$}}};
\node at (2,-1) {\text{\footnotesize{$c$}}};
\end{tikzpicture}
};
\endxy
\end{gather*}

The latter equals $0$ since $F_2^2F_1-(q+q^{-1}) F_2F_1F_2 +F_1F_2^2=0$.
\end{proof}

Note that the usual def\/inition of the evaluation representations requires to have a~$U_q(\glm)$ action.
In our case,
the $U_q(\slm)$ action,
with a~choice of $N$ made,
is actually a~disguised $U_q(\glm)$-representation,
and the factor $q^{\pm(a_1+a_m)}$ utilizes this $\glm$-information.

It will be useful when we turn to skein modules to rescale the braidings and replace the $\Ts$'s by the $\T$'s.
Doing so,
we have analogues of relations \eqref{diagRel1},
\eqref{diagRel2} and \eqref{diagRel3}.
This time,
we omit $\Psi_{\T}$ everywhere
\begin{gather*}
\xy
(0,0)*{
\begin{tikzpicture} [scale=.5]
\draw [very thick, directed=.55] (0,-1) -- (0,1);
\draw [very thick] (0,1) -- (2,3);
\draw [very thick, directed=1] (2,3) -- (2,4);
\draw [very thick, directed=.55] (1,-1) -- (1,1);
\draw [very thick] (1,1) -- (.6,1.4);
\draw [very thick] (.4,1.6) -- (0,2);
\draw [very thick, directed=1] (0,2) -- (0,4);
\draw [very thick] (2,-1) -- (2,0);
\draw [very thick, directed=.55] (2,0) -- (2,2);
\draw [very thick] (2,2) -- (1.6,2.4);
\draw [very thick] (1.4,2.6) -- (1,3);
\draw [very thick,directed=1] (1,3) -- (1,4);
\draw [very thick, directed=.55] (2,.25) -- (1,.75);
\end{tikzpicture}
};
\endxy
=
\xy
(0,0)*{
\begin{tikzpicture} [scale=.5]
\draw [very thick, directed=.55] (0,-1) -- (0,0);
\draw [very thick, directed=.55] (0,0) -- (2,2);
\draw [very thick, directed=1] (2,2) -- (2,4);
\draw [very thick, directed=.55] (1,-1) -- (1,0);
\draw [very thick] (1,0) -- (.6,.4);
\draw [very thick] (.4,.6) -- (0,1);
\draw [very thick, directed=1] (0,1) -- (0,4);
\draw [very thick, directed=.55] (2,-1) -- (2,1);
\draw [very thick] (2,1) -- (1.6,1.4);
\draw [very thick] (1.4,1.6) -- (1,2);
\draw [very thick, directed=1] (1,2) -- (1,4);
\draw [very thick, directed=.55] (1,2.25) -- (0,2.75);
\end{tikzpicture}
};
\endxy
\\
%\label{diagRel2_T}
\xy
(0,0)*{
\begin{tikzpicture} [scale=.5]
\draw [very thick, directed=.55] (0,0) -- (0,2);
\draw [very thick] (0,2) -- (1,3);
\draw [very thick, directed=1] (1,3) -- (1,4);
\draw [very thick, directed=.55] (1,0) -- (1,2);
\draw [very thick] (1,2) -- (.6,2.4);
\draw [very thick] (.4,2.6) -- (0,3);
\draw [very thick, directed=1] (0,3) -- (0,4);
\draw [very thick, directed=.55] (0,1.25) -- (1,1.75);
\node at (0,-.3) {\text{\tiny $a_i$}};
\node at (1,-.3) {\text{\tiny $a_{i+1}$}};
\node at (0,4.3) {};
\node at (1,4.3) {};
\end{tikzpicture}
};
\endxy
=
-q^{a_i-a_{i+1}+1}\xy
(0,0)*{
\begin{tikzpicture} [scale=.5]
\draw [very thick, directed=.55] (0,0) -- (0,1);
\draw [very thick] (0,1) -- (1,2);
\draw [very thick, directed=1] (1,2) -- (1,4);
\draw [very thick, directed=.55] (1,0) -- (1,1);
\draw [very thick] (1,1) -- (.6,1.4);
\draw [very thick] (.4,1.6) -- (0,2);
\draw [very thick, directed=1] (0,2) -- (0,4);
\draw [very thick, directed=.55] (1,2.25) -- (0,2.75);
\node at (0,-.3) {\text{\tiny $a_i$}};
\node at (1,-.3) {\text{\tiny $a_{i+1}$}};
\node at (0,4.3) {};
\node at (1,4.3) {};
\end{tikzpicture}
};
\endxy
, \qquad \xy
(0,0)*{
\begin{tikzpicture} [scale=.5]
\draw [very thick, directed=.55] (0,0) -- (0,2);
\draw [very thick] (0,2) -- (1,3);
\draw [very thick, directed=1] (1,3) -- (1,4);
\draw [very thick, directed=.55] (1,0) -- (1,2);
\draw [very thick] (1,2) -- (.6,2.4);
\draw [very thick] (.4,2.6) -- (0,3);
\draw [very thick, directed=1] (0,3) -- (0,4);
\draw [very thick, directed=.55] (1,1.25) -- (0,1.75);
\node at (0,4.3) {\text{\tiny $a_i$}};
\node at (1,4.3) {\text{\tiny $a_{i+1}$}};
\node at (0,-.3) {};
\node at (1,-.3) {};
\end{tikzpicture}
};
\endxy
=
-q^{a_i-a_{i+1}-1}\xy
(0,0)*{
\begin{tikzpicture} [scale=.5]
\draw [very thick, directed=.55] (0,0) -- (0,1);
\draw [very thick] (0,1) -- (1,2);
\draw [very thick, directed=1] (1,2) -- (1,4);
\draw [very thick, directed=.55] (1,0) -- (1,1);
\draw [very thick] (1,1) -- (.6,1.4);
\draw [very thick] (.4,1.6) -- (0,2);
\draw [very thick, directed=1] (0,2) -- (0,4);
\draw [very thick, directed=.55] (0,2.25) -- (1,2.75);
\node at (0,4.3) {\text{\tiny $a_i$}};
\node at (1,4.3) {\text{\tiny $a_{i+1}$}};
\node at (0,-.3) {};
\node at (1,-.3) {};
\end{tikzpicture}
};
\endxy
\\
%\label{diagRel3_T}
\xy
(0,0)*{
\begin{tikzpicture} [scale=.5]
\draw [very thick, directed=.55] (0,0) -- (0,2);
\draw [very thick] (0,2) -- (.4,2.4);
\draw [very thick] (.6,2.6) -- (1,3);
\draw [very thick, directed=1] (1,3) -- (1,4);
\draw [very thick, directed=.55] (1,0) -- (1,2);
\draw [very thick] (1,2) -- (0,3);
\draw [very thick, directed=1] (0,3) -- (0,4);
\draw [very thick, directed=.55] (0,1.25) -- (1,1.75);
\node at (0,-.3) {};
\node at (1,-.3) {};
\node at (0,4.3) {\text{\tiny $a_i$}};
\node at (1,4.3) {\text{\tiny $a_{i+1}$}};
\end{tikzpicture}
};
\endxy
=
-q^{a_i-a_{i+1}+1}\xy
(0,0)*{
\begin{tikzpicture} [scale=.5]
\draw [very thick, directed=.55] (0,0) -- (0,1);
\draw [very thick] (0,1) -- (.4,1.4);
\draw [very thick] (.6,1.6) -- (1,2);
\draw [very thick, directed=1] (1,2) -- (1,4);
\draw [very thick, directed=.55] (1,0) -- (1,1);
\draw [very thick] (1,1) -- (0,2);
\draw [very thick, directed=1] (0,2) -- (0,4);
\draw [very thick, directed=.55] (1,2.25) -- (0,2.75);
\node at (0,-.3) {};
\node at (1,-.3) {};
\node at (0,4.3) {\text{\tiny $a_i$}};
\node at (1,4.3) {\text{\tiny $a_{i+1}$}};
\end{tikzpicture}
};
\endxy
, \qquad \xy
(0,0)*{
\begin{tikzpicture} [scale=.5]
\draw [very thick, directed=.55] (0,0) -- (0,2);
\draw [very thick] (0,2) -- (.4,2.4);
\draw [very thick] (.6,2.6) -- (1,3);
\draw [very thick, directed=1] (1,3) -- (1,4);
\draw [very thick, directed=.55] (1,0) -- (1,2);
\draw [very thick] (1,2) -- (0,3);
\draw [very thick, directed=1] (0,3) -- (0,4);
\draw [very thick, directed=.55] (1,1.25) -- (0,1.75);
\node at (0,4.3) {};
\node at (1,4.3) {};
\node at (0,-.3) {\text{\tiny $a_i$}};
\node at (1,-.3) {\text{\tiny $a_{i+1}$}};
\end{tikzpicture}
};
\endxy
=
-q^{a_i-a_{i+1}-1}\xy
(0,0)*{
\begin{tikzpicture} [scale=.5]
\draw [very thick, directed=.55] (0,0) -- (0,1);
\draw [very thick] (0,1) -- (.4,1.4);
\draw [very thick] (.6,1.6) -- (1,2);
\draw [very thick, directed=1] (1,2) -- (1,4);
\draw [very thick, directed=.55] (1,0) -- (1,1);
\draw [very thick] (1,1) -- (0,2);
\draw [very thick, directed=1] (0,2) -- (0,4);
\draw [very thick, directed=.55] (0,2.25) -- (1,2.75);
\node at (0,4.3) {};
\node at (1,4.3) {};
\node at (0,-.3) {\text{\tiny $a_i$}};
\node at (1,-.3) {\text{\tiny $a_{i+1}$}};
\end{tikzpicture}
};
\endxy
\end{gather*}

We denote $\tilde{\rho}_a$,
for $a$ complex number,
the analogue of $\rho_a$:
\begin{gather*}
\tilde{\rho}_a(E_0\onel) = aq^{-(a_1+a_m)}C_{T_{m-1}}\cdots C_{T_{2}}(F_1)\onel,
\\
\tilde{\rho}_a(F_0\onel) = a^{-1}q^{a_1+a_m}C_{T_{m-1}}\cdots C_{T_{2}}(E_1)\onel.
\end{gather*}

$\rho_a$ and $\tilde{\rho}_a$ are def\/ined similarly for other generators.

The proof of Proposition~\ref{PropEvalRep} remains valid if we replace the braidings $\Ts$'s by the ($\glm$) rescaled $\T$'s.

\begin{Proposition}%\label{PropEvalRepResc}
The action of $\U_q(\slm)$ on ${\bigwedge}_q(\C^n\otimes \C^m)$ extends via $\tilde{\rho_a}$ to an action of $\U'_q(\hat{\slm})$.
\end{Proposition}

Now,
using the fact that the rescaled braidings $\T$ allow some liberty,
we can give \emph{a posteriori} a~formula very close to Chari-Pressley original evaluation representations~\cite{ChariPressley}.

Assume that we have chosen a~square root $q^{\frac{1}{2}}$ of $q$,
and def\/ine the bracket $[\,\cdot\,,\,\cdot\,]_{q^{\frac{1}{2}}}$ by $[u,v]_{q^{\frac{1}{2}}}=q^{\frac{1}{2}}uv-q^{-\frac{1}{2}}vu$.

\begin{Proposition} \label{prop_CP_evalrep}
We have
\begin{gather*}
\tilde{\rho}_a(E_0\onel)
=(-1)^{m-2}aq^{-a_1-a_m+\frac{m-2}{2}}\big[F_{m-1},\big[F_{m-2},\dots,[F_{2},F_{1}]_{q^{\frac{1}{2}}}\cdots\big]_{q^{\frac{1}{2}}}\big]_{q^{\frac{1}{2}}}
\\
\phantom{\tilde{\rho}_a(E_0\onel)}
=(-1)^{m-2}aq^{-a_1-a_m+\frac{m-2}{2}}\big[\big[\cdots[F_{m-1},F_{m-2}]_{q^{\frac{1}{2}}}\cdots,F_{2}\big]_{q^{\frac{1}{2}}},F_{1}\big]_{q^{\frac{1}{2}}}.
\end{gather*}
\end{Proposition}

\begin{proof}
The core of the identif\/ication is to understand the bracket process in the diagrammatic def\/inition.
By Reidemeister or Kauf\/fman-type moves,
we can rewrite $\tilde{\rho}_a(E_0\onel)$ as (which doesn't make sense in $U_q(\slm)$ anymore):
\begin{gather*}
E_0\onel\quad \mapsto
\quad
a q^{-(a_1+a_m)}
\xy
(0,0)*{
\begin{tikzpicture} [scale=.5]
\draw [very thick, directed=.55] (0,0) -- (0,1);
\draw [very thick] (0,1) -- (0,5);
\draw [very thick, directed=1] (0,5) -- (0,6);
\draw [very thick, directed=.55] (1,0) -- (1,1);
\draw [very thick] (1,1) -- (1,3.6);
\draw [very thick] (1,4.2) -- (1,5);
\draw [very thick, directed=1] (1,5) -- (1,6);
\draw [very thick, directed=.55] (3,0) -- (3,1);
\draw [very thick] (3,1) -- (3,1.8);
\draw [very thick] (3,2.4)-- (3,5);
\draw [very thick, directed=1] (3,5) -- (3,6);
\draw [very thick, directed=.55] (4,0) -- (4,1);
\draw [very thick] (4,1) -- (4,5);
\draw [very thick, directed=1] (4,5) -- (4,6);
\draw [very thick, directed=.55] (4,1.25) -- (0,4.75);
\node at (2,.5) {$\cdots$};
\node at (2,4.5) {$\cdots$};
\node at (0,-.3) {\tiny $a_1$};
\node at (1,-.3) {\tiny $a_2$};
\node at (2.9,-.3) {\tiny $a_{m-1}$};
\node at (4.1,-.3) {\tiny $a_m$};
\end{tikzpicture}
};
\endxy
\end{gather*}
Then,
we can smooth the rightmost crossing using:
\begin{gather*}
\Psi_{\T}\left(
\xy
(0,0)*{
\begin{tikzpicture} [scale=.5,decoration={markings,
mark=at position 1 with {\arrow{>}}; }]
\draw [very thick,
postaction={decorate}] (1,0) -- (0,1);
\draw [very thick] (0,0) -- (.3,.3);
\draw [very thick,
postaction={decorate}] (.7,.7) -- (1,1);
\node at (0,-.5) {\tiny $k$};
\node at (1,-.5) {\tiny $1$};
\node at (0,1.5) {\tiny $1$};
\node at (1,1.5) {\tiny $k$};
\end{tikzpicture}
};
\endxy
\right)
=
-q \;
\xy
(0,0)*{
\begin{tikzpicture} [scale=.5]
\draw [very thick, directed=1] (0,0) -- (0,2);
\draw [very thick, directed=1] (1,0) -- (1,2);
\draw [very thick, directed=.55] (0,.75) -- (1,1.25);
\node at (0,-.5) {\tiny $k$};
\node at (1,-.5) {\tiny $1$};
\node at (.5,.55) {\tiny $k-1$};
\node at (0,2.5) {\tiny $1$};
\node at (1,2.5) {\tiny $k$};
\end{tikzpicture}
};
\endxy
+
\xy
(0,0)*{
\begin{tikzpicture} [scale=.5]
\draw [very thick, directed=.55] (0,0) -- (0,.75);
\draw [very thick] (1,0) -- (1,.5);
\draw [very thick, directed=.55] (1,.5) -- (0,.75);
\draw [very thick] (0,.75) -- (0,1.25);
\draw [very thick,
postaction={decorate}] (0,1.25) -- (1,1.5);
\draw [very thick, directed=1] (1,1.5) -- (1,2);
\draw [very thick, directed=1] (0,1.25) -- (0,2);
\node at (0,-.5) {\tiny $k$};
\node at (1,-.5) {\tiny $1$};
\node at (0,2.5) {\tiny $1$};
\node at (1,2.5) {\tiny $k$};
\end{tikzpicture}
};
\endxy.
\end{gather*}

This gives:
\begin{gather*}
E_0\onel
\mapsto
-q^ea q^{-(a_1+a_m)}
\xy
(0,0)*{
\begin{tikzpicture} [scale=.5]
\draw [very thick, directed=.55] (0,0) -- (0,1);
\draw [very thick] (0,1) -- (0,5);
\draw [very thick, directed=1] (0,5) -- (0,6);
\draw [very thick, directed=.55] (1,0) -- (1,1);
\draw [very thick] (1,1) -- (1,3.6);
\draw [very thick] (1,4.2) -- (1,5);
\draw [very thick, directed=1] (1,5) -- (1,6);
\draw [very thick, directed=.55] (3,0) -- (3,1);
\draw [very thick, directed=.55] (3,1) -- (3,5);
\draw [very thick, directed=1] (3,5) -- (3,6);
\draw [very thick, directed=.55] (4,0) -- (4,1);
\draw [very thick] (4,1) -- (4,5);
\draw [very thick, directed=1] (4,5) -- (4,6);
\draw [very thick, directed=.55] (3,2.25) -- (0,4.75);
\draw [very thick, directed=.55] (4,4.25) -- (3,4.75);
\node at (2,.5) {$\cdots$};
\node at (2,4.5) {$\cdots$};
\node at (0,-.3) {\tiny $a_1$};
\node at (1,-.3) {\tiny $a_2$};
\node at (2.9,-.3) {\tiny $a_{m-1}$};
\node at (4.1,-.3) {\tiny $a_m$};
\end{tikzpicture}
};
\endxy
+a q^{-(a_1+a_m)}
\xy
(0,0)*{
\begin{tikzpicture} [scale=.5]
\draw [very thick, directed=.55] (0,0) -- (0,1);
\draw [very thick] (0,1) -- (0,5);
\draw [very thick, directed=1] (0,5) -- (0,6);
\draw [very thick, directed=.55] (1,0) -- (1,1);
\draw [very thick] (1,1) -- (1,3.6);
\draw [very thick] (1,4.2) -- (1,5);
\draw [very thick, directed=1] (1,5) -- (1,6);
\draw [very thick, directed=.55] (3,0) -- (3,1);
\draw [very thick] (3,1) -- (3,5);
\draw [very thick, directed=1] (3,5) -- (3,6);
\draw [very thick, directed=.55] (4,0) -- (4,1);
\draw [very thick] (4,1) -- (4,5);
\draw [very thick, directed=1] (4,5) -- (4,6);
\draw [very thick, directed=.55] (3,2.25) -- (0,4.75);
\draw [very thick, directed=.55] (4,1.25) -- (3,1.75);
\node at (2,.5) {$\cdots$};
\node at (2,4.5) {$\cdots$};
\node at (0,-.3) {\tiny $a_1$};
\node at (1,-.3) {\tiny $a_2$};
\node at (2.9,-.3) {\tiny $a_{m-1}$};
\node at (4.1,-.3) {\tiny $a_m$};
\end{tikzpicture}
};
\endxy
\\
\phantom{E_0\onel}
=-qa q^{-(a_1+a_m)} F_{m-1}\xy
(0,0)*{
\begin{tikzpicture} [scale=.5]
\draw [very thick, directed=.55] (0,0) -- (0,1);
\draw [very thick] (0,1) -- (0,5);
\draw [very thick, directed=1] (0,5) -- (0,6);
\draw [very thick, directed=.55] (1,0) -- (1,1);
\draw [very thick] (1,1) -- (1,3.6);
\draw [very thick] (1,4.2) -- (1,5);
\draw [very thick, directed=1] (1,5) -- (1,6);
\draw [very thick, directed=.55] (3,0) -- (3,1);
\draw [very thick] (3,1) -- (3,5);
\draw [very thick, directed=1] (3,5) -- (3,6);
\draw [very thick, directed=.55] (4,0) -- (4,1);
\draw [very thick] (4,1) -- (4,5);
\draw [very thick, directed=1] (4,5) -- (4,6);
\draw [very thick, directed=.55] (3,2.25) -- (0,4.75);
\node at (2,.5) {$\cdots$};
\node at (2,4.5) {$\cdots$};
\node at (0,-.3) {\tiny $a_1$};
\node at (1,-.3) {\tiny $a_2$};
\node at (2.9,-.3) {\tiny $a_{m-1}$};
\node at (4.1,-.3) {\tiny $a_m$};
\end{tikzpicture}
};
\endxy
+a q^{-(a_1+a_m)}
\xy
(0,0)*{
\begin{tikzpicture} [scale=.5]
\draw [very thick, directed=.55] (0,0) -- (0,1);
\draw [very thick] (0,1) -- (0,5);
\draw [very thick, directed=1] (0,5) -- (0,6);
\draw [very thick, directed=.55] (1,0) -- (1,1);
\draw [very thick] (1,1) -- (1,3.6);
\draw [very thick] (1,4.2) -- (1,5);
\draw [very thick, directed=1] (1,5) -- (1,6);
\draw [very thick, directed=.55] (3,0) -- (3,1);
\draw [very thick] (3,1) -- (3,5);
\draw [very thick, directed=1] (3,5) -- (3,6);
\draw [very thick, directed=.55] (4,0) -- (4,1);
\draw [very thick] (4,1) -- (4,5);
\draw [very thick, directed=1] (4,5) -- (4,6);
\draw [very thick, directed=.55] (3,2.25) -- (0,4.75);
\node at (2,.5) {$\cdots$};
\node at (2,4.5) {$\cdots$};
\node at (0,-.3) {\tiny $a_1$};
\node at (1,-.3) {\tiny $a_2$};
\node at (2.9,-.3) {\tiny $a_{m-1}$};
\node at (4.1,-.3) {\tiny $a_m$};
\end{tikzpicture}
};
\endxy F_{m-1}
\\
\phantom{E_0\onel}
=-aq^{-a_1-a_m+\frac{1}{2}}\left[F_{m-1},
\xy
(0,0)*{
\begin{tikzpicture} [scale=.5]
\draw [very thick, directed=.55] (0,0) -- (0,1);
\draw [very thick] (0,1) -- (0,5);
\draw [very thick, directed=1] (0,5) -- (0,6);
\draw [very thick, directed=.55] (1,0) -- (1,1);
\draw [very thick] (1,1) -- (1,3.6);
\draw [very thick] (1,4.2) -- (1,5);
\draw [very thick, directed=1] (1,5) -- (1,6);
\draw [very thick, directed=.55] (3,0) -- (3,1);
\draw [very thick] (3,1) -- (3,5);
\draw [very thick, directed=1] (3,5) -- (3,6);
\draw [very thick, directed=.55] (4,0) -- (4,1);
\draw [very thick] (4,1) -- (4,5);
\draw [very thick, directed=1] (4,5) -- (4,6);
\draw [very thick, directed=.55] (3,2.25) -- (0,4.75);
\node at (2,.5) {$\cdots$};
\node at (2,4.5) {$\cdots$};
\node at (0,-.3) {\tiny $a_1$};
\node at (1,-.3) {\tiny $a_2$};
\node at (2.9,-.3) {\tiny $a_{m-1}$};
\node at (4.1,-.3) {\tiny $a_m$};
\end{tikzpicture}
};
\endxy
\right]_{q^{\frac{1}{2}}}.
\end{gather*}

We can then iterate by smoothing successive crossings.
This proves the f\/irst part of the equality.
For the second part,
we perform the same process,
but we start from the left-most crossing.
\end{proof}

Because of the similarity of this def\/inition to the usual one,
we will from now on refer to the representations induced by $\tilde{\rho}_a$ as \emph{evaluation representations}.

This process helps us extend ${\bigwedge}_q(\C^n\otimes \C^m)$ to a~$U'_q(\hslm)$ representation.
The $U'_q(\hslm)$ and $U_q(\sln)$ actions still commute,
but this certainly cannot provide new information,
since the new representation is entirely built on the old one.
An extension of this representation will be studied later.

Since the usual skew Howe duality context is closely related to skein modules (or spider categories),
it is natural to wonder if we can understand a~skein analogue for the evaluation representations.
It appears that closing the Dynkin diagram of $\sln$ corresponds to gluing two opposite sides of the box
on which one usually looks at tangles: this produces an annulus.

\subsection{Annular knots}

Represent ${\bigwedge}_q^N(\C^n\otimes \C^m) \simeq \bigoplus_{a_1+\cdots+a_m=N}{\bigwedge}_q^{a_1}(\C^n)
\otimes \cdots \otimes {\bigwedge}_q^{a_m}(\C^n)$ by a~sequence $(a_1,\dots,a_m)$,
as before, but now drawn on a~circle instead of drawing it on a~segment,
so that $a_1$ and $a_m$ lie next to each other.
$E_0$ is a~map: ${\bigwedge}_q^{a_1}(\C^n)\otimes \cdots \otimes {\bigwedge}_q^{a_m}(\C^n)
\rightarrow {\bigwedge}_q^{a_1+1}(\C^n)\otimes \cdots \otimes {\bigwedge}_q^{a_m-1}(\C^n)$.
This idea gives a~diagrammatic presentation of the af\/f\/ine extension of the skew Howe duality phenomenon.
See below the diagram corresponding to $E_0E_1$ acting on a~sequence $(2,0)$,
with $m=2$, $n=2$, $N=2$:
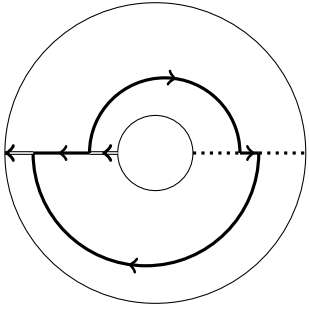
\begin{gather*}
\xy
(0,0)*{
\begin{tikzpicture} [scale=.5]
\draw (-1,0) arc (180:-180:1);
\draw[double, directed=.55] (-1,0) -- (-1.75,0);
\draw[very thick, directed=.55] (-1.75,0) arc (180:0:2);
\draw [very thick,
dotted] (1,0) -- (2.25,0);
\draw [very thick,
dotted] (2.75,0) -- (4,0);
\draw[very thick,directed=.8] (2.25,0) -- (2.75,0);
\draw[very thick, directed=.55] (2.75,0) arc(0:-180:3);
\draw[very thick, directed=.55] (-1.75,0) -- (-3.25,0);
\draw[double, directed=1] (-3.25,0) -- (-4,0);
\draw (-4,0) arc (180:-180:4);
\end{tikzpicture}};
\endxy
\end{gather*}

Annular webs can be def\/ined in a~very similar fashion as in the plane case.
One just embeds trivalent graphs in the annulus instead of the plane,
and all relations being local,
they look the same for webs considered in any surface.

\begin{Definition} \label{def:AWeb}
Let $n{\bf AWeb}$,
the $\sln$ annular web skein module,
be the $\Z[q,q^{-1}]$-module ge\-ne\-rated by annular webs (oriented trivalent graphs with preserved f\/low embedded in an annulus)
possibly with boundary (embedded in the two boundary circles of the annulus) up to isotopy and the local $\sln$ web relations \eqref{webrel1},
\eqref{webrel2} and \eqref{webrel3}.
\end{Definition}

All web relations and all $\hslm$ relations being ``local relations'',
it is easy to observe that web relations are implied by $\hslm$.
However,
if we consider the evaluation representation to extend the $\slm$ action,
there are other relations than web relations: for example,
the above diagram corresponds to a~scalar action,
while in the skein module,
this would rather correspond to a~generator.
We shall f\/irst identify here a~skein module that corresponds to these representations,
before seeking for a~situation closer to the skein module of the annulus.
For this purpose,
let us introduce the annular extension of the category $n{\bf Web}^+_m$.

\begin{Definition}
Def\/ine $n{\bf AWeb}^+_m$ to be the category with objects,
sequences $(a_1,\dots,a_m)$ ($0\leq a_i\leq n$) labeling points regularly drawn on a~circle,
together with a~zero object.
Points labeled by zero can be erased,
but,
as in~\cite{Blan} and~\cite{LQR1},
we will keep the $n$-strands as well.
Morphisms are formal sums over $\Z[q^{\pm1}]$ of upward $\sln$-webs (in the sense of Def\/inition~\ref{def:AWeb}) drawn on a~cylinder,
generated by\footnote{Note that this a~priori doesn't recover the complete braid group of the annulus.
We miss in the braid group associated to $U_q(\hslm)$ the element given by sending each point on one boundary of the annulus
to the one corresponding to its right (for example) neighbor in the other boundary of the annulus.
In case of ladders,
this can be artif\/icially solved by adding a~$0$-strand.}:
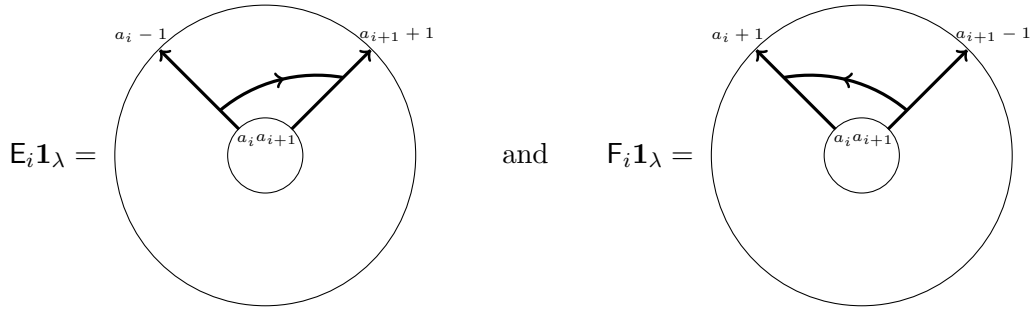
\begin{gather*}
\E_i\onel
=\xy
(0,0)*{
\begin{tikzpicture} [scale=.5]
\draw (-1,0) arc (180:-180:1);
\draw[very thick, directed=1] (-.7,.7) -- (-2.8,2.8);
\draw[very thick, directed=1] (.7,.7) -- (2.8,2.8);
\draw [very thick, directed=.55] (-1.2,1.2) arc (130:80:4);
\draw (-4,0) arc (180:-180:4);
\node at (.3,.4) {\text{\tiny{$a_{i+1}$}}};
\node at (-.5,.4) {\text{\tiny{$a_i$}}};
\node at (3.5,3.2) {\text{\tiny{$a_{i+1}+1$}}};
\node at (-3.3,3.2) {\text{\tiny{$a_i-1$}}};
\end{tikzpicture}};
\endxy
\qquad \text{and}
\qquad
\F_i\onel
=\xy
(0,0)*{
\begin{tikzpicture} [scale=.5]
\draw (-1,0) arc (180:-180:1);
\draw[very thick, directed=1] (-.7,.7) -- (-2.8,2.8);
\draw[very thick, directed=1] (.7,.7) -- (2.8,2.8);
\draw [very thick, directed=.55] (1.2,1.2) arc (50:100:4);
\draw (-4,0) arc (180:-180:4);
\node at (.3,.4) {\text{\tiny{$a_{i+1}$}}};
\node at (-.5,.4) {\text{\tiny{$a_i$}}};
\node at (3.5,3.2) {\text{\tiny{$a_{i+1}-1$}}};
\node at (-3.3,3.2) {\text{\tiny{$a_i+1$}}};
\end{tikzpicture}};
\endxy
\end{gather*}
In the above pictures,
we have drawn only strands $i$ and $i+1$.
\end{Definition}

$n{\bf{AWeb}}^+_m(N)$ will be the full subcategory with objects such that $\sum\limits a_i=N$.

Given $m$,
$n$,
$N$,
we can def\/ine on weights a~map $\Phi$ that sends $\lambda=(\lambda_0,\dots,
\lambda_{m-1})$ to a~sequence $(a_1,\dots,a_m)$ such that $a_{i+1}-a_i=\lambda_i$ and $a_1-a_{m-1}=\lambda_0$,
with $\sum\limits a_i=N$.
If such a~solution doesn't exist,
the weight is sent to the zero object.

\begin{Proposition}%\label{prop:affineFunctor}
For all $m$,
$n$,
$N$,
there is a~functor $\Phi\;\colon\;\U'_q(\hslm)\;\mapsto\; n{\bf{AWeb}}^+_m(N)$ defined on weights as above and on morphisms
by $E_i\onel\mapsto \E_i\onel$ and $F_i\onel\mapsto \F_i\onel$.
\end{Proposition}

$\hslm$ relations are locally $\slm$ relations: the proof is straightforward.
\cite[Proposition~7.4]{Cautis} will receive a~direct translation.
Let us consider a~knotted annular tangle as a~composition of $\E_i\onel$,
$\F_i\onel$ and crossings between the $i$-th and $(i+1)$-th strands.
From such a~presentation,
we can read of\/f an element of $\tilde{\U'_q(\hslm)}$ as the corresponding product of $E_i\onel$,
$F_i\onel$ and $\T^{\pm 1}$.
Let us call $X_w$ this element corresponding to a~presentation of a~web-tangle $w$.

\begin{Proposition}%\label{prop:affineInvEvalrep}
The evaluation representations produce annular web invariants for ladder-type webs.
In other words,
if $X_w \in \tilde{\U'_q(\hslm)}$ corresponds to a~presentation $w$ of a~web-tangle,
then the morphism of $U_q(\sln)$ representations given by $X_w$ is an invariant of the web-tangle.
\end{Proposition}

$\tilde{\U'_q(\hslm)}$ above denotes the completion of $\U'_q(\hslm)$ given by the quotient of the ring of series
of elements of $\U'_q(\hslm)$ acting by evaluation representation on each irreducible $\U'_q(\slm)$ representation $V_{\lambda}$
by zero but for f\/initely many terms,
mod out by the two-sided ideal of elements acting by zero on all $V_{\lambda}$.

As explained previously,
the above process produces an invariant of annular web-tangles.
However,
it appears to come with only little information about the topology of the annulus,
and we can indeed explicitly identify this invariant.
It may be useful for this purpose to see the annulus as a~cylinder rather than f\/lattened on a~plane: this way,
we can f\/ill it.

\begin{Proposition}%\label{prop:identifyingFilledCylinder}
The evaluation representation with $a=-q^{n+1}$ recovers the $\sln$ skein module of the filled cylinder.
In other words, for each fixed value of $N$ we have the following commutative diagram:
\begin{gather*}
\xymatrix{
\U'_q(\hslm) \ar[r]^{\Phi} \ar[d]^{ev} & n{\bf{AWeb}}^+_m \ar[d]^{\rm Filling}
\\
\U_q(\slm) \ar[r]^{\Phi} & n{\bf{Web}}^+_m
}
\end{gather*}

\end{Proposition}

\begin{proof}
We just have to check that the action of $E_0\onel$ and $F_0\onel$ corresponds to $\E_0\onel$ and $\F_0\onel$ seen
in the skein module of the f\/illed cylinder.

For $E_0$ for example,
we have the following situation:
\begin{gather*}
{\rm Filling}(\Phi(E_0\onel))
=
\Psi_{\T}\left(
\xy
(0,0)*{
\begin{tikzpicture} [scale=.7,decoration={markings,
mark=at position 0.5 with {\arrow{>}}; }]
\draw [very thick,
postaction={decorate}] (2,0) -- (2,1);
\node at (2,-.4) {\footnotesize $a_m$};
\draw [very thick] (2,1) .. controls (2,1.2) and (3,1.8) .. (3,2) -- (3,2) .. controls (3,2.5)
and (-1,3.5) .. (-1,4) -- (-1,4) .. controls (-1,4.2) and (0,4.8)..
(0,5);
\draw [very thick] (2,1) -- (2,2.4);
\draw [very thick] (2,2.8) -- (2,5);
\draw [very thick,
->] (2,5) -- (2,6);
\node at (2,6.4) {\footnotesize $a_m-1$};
\draw [very thick,
postaction={decorate}] (0,0) -- (0,1);
\node at (0,-.4) {\footnotesize $a_1$};
\draw [very thick] (0,1) -- (0,3.2);
\draw [very thick] (0,3.6) -- (0,5);
\draw [very thick,
->] (0,5) -- (0,6);
\node at (0,6.4) {\footnotesize $a_1+1$};
\node at (1,3) {$\circ$};
\node at (1.4,3.3) {\text{\footnotesize{$-1$}}};
\end{tikzpicture}};
\endxy
\right)
\; =\;
\Psi_{\T}\left(
\xy
(0,0)*{
\begin{tikzpicture} [scale=.7,decoration={markings,
mark=at position 0.5 with {\arrow{>}}; }]
\draw [very thick,
postaction={decorate}] (2,0) -- (2,1);
\node at (2,-.4) {\footnotesize $a_m$};
\draw [very thick] (2,1) .. controls (2,1.2) and (0,4.8) .. (0,5);
\draw [very thick,
->] (0,5) -- (0,6);
\node at (2,6.4) {\footnotesize $a_{m}-1$};
\draw [very thick,
postaction={decorate}] (0,0) -- (0,1);
\draw [very thick] (0,1) -- (0,5);
\node at (0,-.4) {\footnotesize $a_1$};
\draw [very thick] (2,1) -- (2,5);
\draw [very thick,
->] (2,5) -- (2,6);
\node at (0,6.4) {\footnotesize $a_1$+1};
\node at (.15,3.8) {\text{\footnotesize{$\circ \frac{1}{2}$}}};
\node at (.45,5.4) {\text{\footnotesize{$\circ -\frac{1}{2}$}}};
\node at (2.45,.6) {\text{\footnotesize{$\circ -\frac{1}{2}$}}};
\node at (2.15,2.2) {\text{\footnotesize{$\circ \frac{1}{2}$}}};
\end{tikzpicture}};
\endxy
\right)
\end{gather*}

The negative twist on the l.h.s.\ of the above equation comes from the fact that the strand goes along the cylinder with framing parallel
to the cylinder.
When f\/illing the cylinder and projecting it on the back side,
this produces a~twist.

We have depicted here only the leftmost and the rightmost strands.
This whole process is to be understood in front of the other strands.
Then,
a succession of Reidemeister II moves presents this piece of tangle as the elements def\/ining the evaluation representation,
presented in~\eqref{evalRep}.

In the previous computation,
the twists produce the following coef\/f\/icient:
\begin{gather*}
t_{a_1}^{\frac{1}{2}}t_{a_1+1}^{-\frac{1}{2}}t_{a_m}^{-\frac{1}{2}}t_{a_{m-1}}^{\frac{1}{2}}=-q^{-(a_1+a_m)+n+1},
\end{gather*}
while the evaluation representation provides $aq^{-(a_1+a_m)}$.
Choosing $a=-q^{n+1}$ adjusts the coef\/f\/icients.
Checking the results for $F_0$ is similar.
\end{proof}

So,
it appears that extending the skew Howe duality phenomenon to the af\/f\/ine case by evaluation representation gives a~coherent process,
but is too weak to recover the skein module of the annulus.
We miss the fact that acting by $E_1E_2\cdots E_{m-1}E_0$,
although it does not change anything on the weight,
has no reason to be something trivial in the skein module.
This is a~well-known phenomenon in the study of $\hslm$: if we want to understand the non-triviality of this action,
we have to keep the whole data coming from the Dynkin diagram and not only generators $E_i$,
$F_i$,
$K_i^{\pm}$.
We should therefore work in the whole $U_q(\hslm)$ and not only with $U'_q(\hslm)$.

Nonetheless,
we want to keep working with an analogue of Howe duality,
and it would be convenient to have a~process built on these particular representations.
It turns out that there is an easy way to do it,
called \emph{af\/f\/inization},
as explained for example in~\cite[p.233]{HongKang}.

\subsection{Af\/f\/inization}

Following~\cite{HongKang},
we now consider $\C(q)[z,z^{-1}]\otimes_{\C(q)}{\bigwedge}_q^N(\C^n\otimes \C^m)$.
The $U_q(\sln)$ action can be extended by acting by an identity on the $z$-part,
and the previous $U'_q(\hslm)$ action may be extended to an $U_q(\hslm)$ action by the following rules:
\begin{gather*}
E_0\big(z^m\otimes v\big)=z^{m+1}\otimes (E_0v),
\qquad
E_i\big(z^m\otimes v\big) = z^m\otimes (E_iv) \quad \text{for}\quad  i\neq0,
\\
F_0\big(z^m\otimes v\big) = z^{m-1}\otimes (F_0v),
\qquad
F_i\big(z^m\otimes v\big) = z^m \otimes (F_iv) \quad \text{for}\quad  i\neq 0,
\\
K_i\big(z^m\otimes v\big) = z^m \otimes (K_iv),
\qquad
K_d\big(z^m\otimes v\big) = q^m z^m\otimes v.
\end{gather*}
Here,
$K_d$ is the derivation element,
that was neglected in the previous subsection.

We have the same decomposition as before:
\begin{gather*}
\bigoplus_{a_1+\cdots+a_m=N}\C(q)\big[z,z^{-1}\big]\otimes_{\C(q)}{\bigwedge}_q^{a_1}(\C^n)\otimes \cdots \otimes {\bigwedge}_q^{a_m}(\C^n),
\end{gather*}
and more precisely:
\begin{gather*}
\bigoplus_{a_1+\cdots+a_m=N,k}\C(q)\cdot z^k\otimes{\bigwedge}_q^{a_1}(\C^n)\otimes \cdots \otimes {\bigwedge}_q^{a_m}(\C^n).
\end{gather*}
The $U_q(\hslm)$-weight of one summand is $(a_1-a_m,
a_2-a_1,\dots,
a_m-a_{m-1})+k\delta$, where $\delta$ is the null root (that was neglected in the previous subsections).
Note that these representations are of level $0$.

We obtain here a~new process: given a~knotted ladder,
we can turn it into an element of $\U_q(\hslm)$ acting on ${\bigwedge}_q^N(\C^n\otimes \C^m)$,
and extend this into an action on $\C(q)[z^{\pm 1}]\otimes_{\C(q)}{\bigwedge}_q^N(\C^n\otimes \C^m)$.
If we restrict to knots drawn on a~cylinder,
that is ladders with a~boundary sequence with only~$0$'s and~$n$'s,
we obtain an element of $\End ( \C(q)[z^{\pm 1}]\otimes_{\C(q)}\C[q^{\pm 1}])=\C(q)[z^{\pm 1}]$.

Since all relations are local,
\cite[Proposition~7.4]{Cautis} receives a~direct translation:

\begin{Proposition}
The previous process defines a~web tangle invariant.
In other words,
if $\onell{\lambda'}X_w\onel\in \tilde{\U'_q(\hslm)}$ corresponds to a~presentation $w$ of an annular web-tangle mapping
sequences $(a_1^{\lambda},\dots,
a_m^{\lambda})$ to $(a_1^{\lambda'},\dots,
a_m^{\lambda'})$,
then the morphism of $U_q(\sln)$ representations mapping
$\C(q)[z,z^{-1}]\otimes_{\C(q)}{\bigwedge}_q^{a_1^{\lambda}}(\C^n)\otimes \cdots \otimes {\bigwedge}_q^{a_m^{\lambda}}(\C^n)$
to $\C(q)[z,z^{-1}]\otimes_{\C(q)}{\bigwedge}_q^{a_1^{\lambda'}}(\C^n)\otimes \cdots \otimes {\bigwedge}_q^{a_m^{\lambda'}}(\C^n)$
given by $X_w$ is an invariant of the web-tangle.
\end{Proposition}

So, to an annular web-tangle presented in a~ladder form,
we can assign a~morphism of $U_q(\sln)$ representations which is an invariant of the web-tangle.
This morphism may be expressed in a~diagrammatic way,
producing a~skein element.
A~natural question would be to know whether the web relations alone are suf\/f\/icient to identify if two diagrams yield the same
skein element (see for example~\cite{GrantModuli} for the planar case).
Unfortunately,
this process is still not faithful enough: $E_1F_0F_1E_0\onel$ acting for example on a~2-strands sequence $(0,2)$ acts as $[2]^2\cdot \onel$,
while the skein element corresponding to both these elements are not equal:
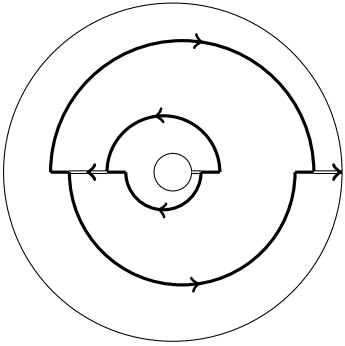
\begin{gather*}
\xy
(0,0)*{
\begin{tikzpicture} [scale=.5]
\draw (-.5,0) arc (180:-180:.5);
\draw[double] (.5,0) -- (.75,0);
\draw[very thick, directed=.55] (.75,0) arc (0:-180:1);
\draw [very thick] (.75,0) -- (1.25,0);
\draw [very thick] (-1.25,0) -- (-1.75,0);
\draw [very thick, directed=.55] (1.25,0) arc (0:180:1.5);
\draw [double, directed=.55] (-1.75,0) -- (-2.75,0);
\draw[very thick, directed=.55] (-2.75,0) arc (-180:0:3);
\draw [very thick] (-2.75,0) -- (-3.25,0);
\draw [very thick] (3.25,0) -- (3.75,0);
\draw [very thick, directed=.55] (-3.25,0) arc (180:0:3.5);
\draw [double, directed=1] (3.75,0) -- (4.5,0);
\draw (-4.5,0) arc (180:-180:4.5);
\end{tikzpicture}};
\endxy
\end{gather*}

It would be interesting to compute the kernel of the map,
but we do not know how to address this question.
This also suggests to look for richer $U_q(\hslm)$ representations extending the skew-Howe duality phenomenon.

The previous invariant contains two pieces of information: the same as the evaluation representation,
that is,
the skein element associated to the web tangle in the skein module of the f\/illed annulus,
and an information given by the action on the $z$-part.
This traces a~kind of algebraic linking number with the core of the annulus.
The problem is that this algebraic number doesn't detect the possible non-triviality of an algebraically non-linked web,
as explained above: the representations we have been working with are still too weak to give a~full representation-theoretic counterpart
of the annular webs.

In the $\slnn{2}$-case,
the unoriented skein module of the annulus is well understood,
isomorphic to $\Z[q^{\pm 1}][z]$,
with $z$ the generator given by a~circle around the hole.
Note that,
if we try to compare the obtained invariant with the usual unoriented $\slnn{2}$ skein module of the annulus,
the f\/irst dif\/ference comes from the fact that the $2$-labeled strands don't play any role in the unoriented version,
while when wrapped around the hole,
they produce a~coef\/f\/icient $z^{\pm 2}$ in the oriented version.
A~second dif\/ference comes from the issue explained above.

An easy way to compare the (usual,
unoriented) Kauf\/fman bracket skein module of the annulus and the invariant obtained using enhanced $\slnn{2}$ webs
and the af\/f\/ine representation would thus be to mod both out by the ideal generated by $(z^2-1)$,
after some renormalizations (since $z$ in the Kauf\/fman bracket version corresponds in spirit to $[2]z$ in the oriented case).
More ref\/ined comparisons seem to non-trivially involve the orientation and the behavior of the $2$-labeled lines,
which makes it very hard to control the power of $z$.

\subsection[Forgetting about $\sln$ ...]{Forgetting about $\boldsymbol{\sln}$ \dots}

Note that the issue that prevents us to obtain an algebraic object that would mimic the behavior of the skein module comes
from the fact that we consider particular $U_q(\sln)$ representations that are not powerful enough to detect all the topological data.
A~kind of virtual analogue would be to only keep the $U_q(\slm)$-part in the duality,
mod out by information extracted from the usual case,
and extend only this to the annular case.

Recall from~\cite{CKM} that we can understand the quotient of $U_q(\slm)$ which corresponds to classes of $\sln$-webs.
Fix $N$ (from the ladders we are looking at) and a~dominant weight $\lambda$ (corresponding to a~sequence $(a^{\lambda}_1,\dots,a^{\lambda}_m)$),
and denote $I_{\lambda}$ the ideal of $\U_q(\slm)$ generated by all weights which do not lie in the Weyl orbit
of any weight $\mu$ so that $\lambda$ dominates $\mu$.
Furthermore,
denote $\U_q(\slm)^n$ the quotient of $\U_q(\slm)$ by the set of weights whose associated sequence $(a_1,\dots,a_m)$
either does not exist or has at least one coef\/f\/icient $a_i< 0$ or $a_i>n$.

Then,
Cautis,
Kamnitzer and Morrison~\cite[Theorem 4.4.1,
Lemma 4.4.2]{CKM} tell us that the morphism $\U_q(\slm)^n/I_{\lambda}\mapsto n{\bf Web}^+_m(a^{\lambda}_1,\dots,a^{\lambda}_m)$
is an isomorphism\footnote{Again,
note that the result they state uses $\glm$ and not $\slm$.
The dependency in the $\glm$ weights in our case is hidden in the $n$-bounded quotient: we need $N$ to know the value of the $a_i$'s.},
where $n{\bf Web}^+_m(a^{\lambda}_1,\dots,a^{\lambda}_m)$ is the algebra of ladders that can be reached starting
from the sequence $(a^{\lambda}_1,\dots,a^{\lambda}_m)$.
In other words,
these are the ladders $W\colon (a^1_1,\dots,a^1_m)\mapsto (a^2_1,\dots,a^2_m)$,
so that the set of ladders between $(a^{\lambda}_1,\dots,a^{\lambda}_m)$ and $(a^1_1,\dots,a^1_m)$ is non-empty.

Let us now denote $n{\bf{AWeb}}^+_m(a_1,\dots,a_m)$ the algebra of annular ladder webs built from a~sequence $(a_1,\dots,a_m)$.
Each $I_{\lambda}$ extends in the af\/f\/ine case to a~module $I_{\hat{\lambda}}$ simply given by assigning to any $\mu\in I_{\lambda}$
the af\/f\/ine weight $\hat{\mu}$ so that if $\mu=(a^{\mu}_2-a^{\mu}_1,\dots,a^{\mu}_m-a^{\mu}_{m-1})$,
$\hat{\mu}=(a^{\mu}_1-a^{\mu}_{m},a^{\mu}_2-a^{\mu}_1,\dots,a^{\mu}_m-a^{\mu}_{m-1})$.
We can then consider the quotient $\U'_q(\hslm)^n$ of $\U'_q(\hslm)$ by weights whose associated sequence has indices lower than $0$
or bigger than $n$, and we have the quotient $\U'_q(\hslm)^n/I_{\hat{\lambda}}$.
The next result then gives us somehow the result we hoped to f\/ind with an explicit $U_q(\sln)$-representation
but just by looking on the dual side of the picture!

\begin{Proposition}%\label{prop:CKMAffine}
For $\lambda$ a~dominant $\slm$ weight,
the map:
\begin{gather*}
\U'_q\big(\hslm\big)^n/I_{\hat{\lambda}} \mapsto n{\bf{AWeb}}^+_m(a^{\lambda}_1,\dots,a^{\lambda}_m)
\end{gather*}
is an isomorphism.
The pre-image of a~knotted ladder in $\U'_q(\hslm)^n/I_{\hat{\lambda}}$ is therefore an invariant of the knot.
\end{Proposition}

\begin{proof}
First,
note that since $E_0$ acts on weights as $F_1\cdots F_n$,
the weights $\hat{\mu}$ are those that cannot be reached from $\hat{\lambda}$.
On objects,
the statement is therefore obvious.

The surjectivity on morphisms comes from the def\/inition.

For morphisms,
the injectivity argument is the same as before: since all relations are local (and elementary relations involve at most three strands),
either the generators $E_0$ and $F_0$ are not involved and we can assume we work with $\slm$,
or they are but (for $m\geq 3$) there exists $i$ so that $E_i$ and $F_i$ are not involved.
There is then an inclusion of $U_q(\slm)$ in $U_q(\hslm)$ that does not involve $E_i$ and $F_i$ and we can assume we work in this one.
\end{proof}

\subsection{\dots\ to better recover it?}

The objects that appear in the previous paragraph are closely related to $q$-Schur algebras and af\/f\/ine versions of them.
We refer to~\cite{DotyGreen} for a~clear presentation of the context in which they appear.

Following~\cite{DotyGreen},
we will consider in what follows an extension\footnote{This extension has the nice property that it gives us the missing
generator of the braid group of the annulus that was previously discussed.} of $U'_q(\hslm)$ with two extra generators $R$ and $R^{-1}$
subject to relations (indices are to be understood modulo $m$):
\begin{gather*}
\begin{split}
& RR^{-1}=R^{-1} R= 0,
\qquad
R^{-1}K_{i+1}R=K_i,
\qquad
R^{-1}K_{i+1}^{-1}R=K_i^{-1},
\\
& R^{-1} E_{i+1}R=E_i,
\qquad
R^{-1}F_{i+1}R=F_i.
\end{split}
\end{gather*}

This algebra will be denoted $\hat{U}_q(\hslm)$,
and its idempotented version $\hat{\U}_q(\hslm)$.
Note that the previous relations in particular give us that $E_0=R^{-1}E_1R$ and $F_0=R^{-1}F_1R$.

Doty and Green suggest us to replace the fundamental representation $\C^m$ of $U_q(\slm)$ by an inf\/inite-dimensional
version $V_{\infty}=\C^{\infty}=\langle X_i,i\in \Z\rangle$ with action:
\begin{gather*}
E_iX_j=X_{j+1}\quad \text{if} \quad  i=j\; \bmod(m),
\qquad
E_iX_j=0\quad \text{otherwise,}
\\
F_iK_{j+1}=X_j\quad \text{if}\quad i=j\; \bmod(m),
\qquad
F_iX_j=0\quad \text{otherwise,}
\\
K_i X_j=q^{-1}X_j\quad \text{and}\quad K_iX_{j+1}=qX_{j+1}\quad \text{if}\quad i=j\; \bmod(m),
\\
K_iX_j=X_j\quad \text{otherwise},
\qquad
RX_j=X_{j+1}.
\end{gather*}

As in the linear case,
we can endow $\C^n\otimes V_{\infty}$ with two commutative actions of $U_q(\sln)$ and $\hat{U}_q(\hslm)$.
We now wish to consider the quantum exterior power of this tensor product of representations and perform the same type of process as before.
However,
this quantum exterior power is more complicated to def\/ine in the af\/f\/ine case than in the linear one.

We refer to~\cite{Stern,TakemuraUglov,Uglov99} and references therein for details about these representations and tools one could use to study them.
We intend here to sketch a~process allowing us to relate $m$-uprights annular ladders to $\sln$-representation theory,
but the question of completely understanding this relation,
and also relating higher-level analogous phenomenon to knot theory remains open.

Denote $V_m=\langle X_1,\dots,X_m\rangle$ that we see as a~subspace of $V_{\infty}$.
$V_m$ is the vector representation of $U_q(\slm)$ (but is not a~$U_q(\hslm)$ module),
and its af\/f\/inization $V_m\otimes \C[z^{\pm 1}]$ is isomorphic to $V_{\infty}$.
The precise def\/inition of the quantum exterior power of $V_{\infty}$ can be found in~\cite[Section 2.1]{Uglov99}.
We have in particular: $X_{i+km}\wedge X_{i+lm}+X_{i+lm}\wedge X_{i+km}=0$.

Uglov's process consists in making ${\bigwedge}_q\left(\C^n\otimes \C[z^{\pm 1}] \otimes V_m \right)$ into a~$U'_q(\hsln)\otimes U'_q(\hslm)$ module.
In particular,
this implies that ${\bigwedge}_q\left(\C^n\otimes V_{\infty} \right)$ has two commuting actions of $U_q(\sln)$ and $U'_q(\hslm)$,
and it is not hard to see that the latter can be extended into a~$\hat{U}_q(\hslm)$ action.

We now want to relate annular webs,
the algebra $\hat{\U}_q(\hslm)$ and the morphisms of the $U_q(\sln)$ representation ${\bigwedge}_q\left(\C^n\otimes V_{\infty} \right)$.

Let us slightly extend the def\/inition of $n{\bf{AWeb}}^+_m(N)$ into $n{\bf{\hat{A}Web}}^+_m(N)$ by adding to it the image of $R$
as the next elementary annular ladder web:
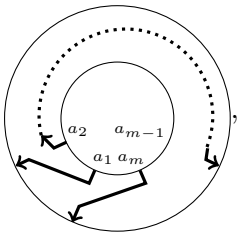
\begin{gather*}
\xy
(0,0)*{
\begin{tikzpicture} [scale=.75,decoration={markings,
mark=at position 0.5 with {\arrow{>}}; }]
\draw (1,0) arc (0:360:1);
\draw (2,0) arc (0:360:2);
\draw [very thick,
->] (1.6,-.53) -- (1.57,-.73) -- (1.8,-.84);
\draw [very thick,
->] (.4,-.92) -- (.5,-1.15) -- (-.68,-1.56) -- (-.8,-1.84);
\draw [very thick,
->] (-.4,-.92) -- (-.5,-1.15) -- (-1.57,-.73) -- (-1.8,-.84);
\draw [very thick,->] (-.9,-.42) -- (-1.125,-.525) -- (-1.35,-.3);
\draw [dotted,
very thick] (-1.35,-.3) arc (193:-23:1.54);
\node at (-.25,-.72) {\tiny $a_1$};
\node at (-.7,-.22) {\tiny $a_2$};
\node at (.4,-.25) {\tiny $a_{m-1}$};
\node at (.25,-.72) {\tiny $a_m$};
\end{tikzpicture}
};
\endxy,
\end{gather*}
and similarly for $R^{-1}$.

It is a~direct extension of the usual case that we have an isomorphism $\oplus_N \hat{\U}_q(\hslm)^{n}\mapsto n{\bf{\hat{A}Web}}^+_m$.

Usually,
when we relate $U_q(\slm)$ and $U_q(\sln)$ endomorphisms,
the fact that we can kill in~$U_q(\slm)$ all weights corresponding to sequences $(a_1,\dots,a_m)$ with an $a_i>n$ comes from the fact
that exterior powers of $\C^n$ of degree more than $n$ are zero.
Here,
we can decompose:
\begin{gather*}
{\bigwedge}^N_q\big(\C^n\otimes V_{\infty} \big)
\simeq
{\bigwedge}^N_q\big(\C^n\otimes \C\big[z^{\pm 1}\big]
\otimes V_m \big)\simeq {\bigwedge}^N_q\big(\C^n\otimes \C\big[z^{\pm 1}\big] \otimes (\C\oplus \C \oplus \cdots \oplus \C) \big)
\\
\phantom{{\bigwedge}^N_q\big(\C^n\otimes V_{\infty} \big)}{}
\simeq
\bigoplus_{a_1+\cdots+a_m=N}{\bigwedge}^{a_1}_q\big(\C^n\otimes \C\big[z^{\pm 1}\big]\big)\otimes\cdots \otimes {\bigwedge}^{a_m}_q\big(\C^n\otimes \C\big[z^{\pm 1}\big]\big).
\end{gather*}
However,
the relation $X_{i+km}\wedge X_{i+lm}+X_{i+lm}\wedge X_{i+km}=0$ does not ensure that ${\bigwedge}^{n}_q(\C^n\otimes \C[z^{\pm 1}])\simeq \C$
nor that it is zero after.

There is nonetheless a~case where this issue does not arise: it is when $N\leq n$ (this corresponds to looking at a~``generic'' version,
where the choice of $n$ actually doesn't matter).
Let us assume that we are in that situation.
We then have:
\begin{gather*}
n{\bf{\hat{A}Web}}^+_m(N) \simeq \hat{\U}_q(\hslm)^{n}(N) \mapsto \End_{U_q(\sln)}\left({\bigwedge}_q^N\big(\C^n\otimes V_{\infty}\big)\right),
\end{gather*}
where $\hat{\U}_q(\hslm)^{n}(N)$ denotes the full subcategory of $\hat{\U}_q(\hslm)^{n}$ whose lifts $(a_1,\dots,a_m)$
of the weights $(\lambda_1,\dots,\lambda_{m-1})$ are such that $\sum\limits a_i=N$ (note again that here,
taking the $n$-bounded quotient only kills weight whose lifts have negative entries,
as all entries will be less than or equal to $n$).
We want to show that this map is an inclusion.
A~very useful tool for this is provided by~\cite[Proposition 3.14]{McKTh}.

Assume that $N<m$ (we can restrict to that case by adding $0$-labeled strands),
and deno\-te~$\1_r$ the idempotent corresponding to the sequence $(a_1,\dots,a_m)=(1,\dots,1,0,\dots,0)$ contai\-ning~$r$ times the number~$1$.
Then, $\1_r\hat{\U}_q(\hslm)\1_r$ is isomorphic to the af\/f\/ine Hecke algebra~$\hat{\mathcal{H}}_{\hat{A}_{r-1}}$ (see~\cite{DotyGreen, McKTh}).
We can present it as generated by $b_1,\dots,b_r$ and $T_{\rho}$,
$T_{\rho}^{-1}$ subject to the following relations:
\begin{gather*}
b_i^2= \big(q+q^{-1}\big)b_i
\quad
\text{for}
\quad
i=1,\dots,r,
\\
b_ib_j=b_jb_i
\quad
\text{for distant}
\quad
i,j=1,\dots,r,
\\
b_ib_{i+1}b_i+b_{i+1}=b_{i+1}b_ib_{i+1}+b_i
\quad
\text{for}
\quad
i=1,\dots,r,
\\
T_{\rho}b_iT_{\rho}^{-1}=b_{i+1}
\quad
\text{for}
\quad
i=1,\dots,r.
\end{gather*}
The last equation is to be understood with indices modulo $r$.
Note that we can obtain the generator $b_r$ as $b_r=T_{\rho}b_{r-1}T_{\rho}^{-1}$.

Mackaay and Thiel then state, where $\hat{S}(m,N)$ is the af\/f\/ine Schur algebra:

\begin{Proposition}[\protect{\cite[Proposition~3.14]{McKTh}}]\label{prop:McKThHecke}
Let $N<m$. Suppose that $A$ is a~$\Q(q)$ algebra and
\begin{gather*}
f\colon \hat{S}(m,N)\mapsto A
\end{gather*}
is a~surjective $\Q(q)$-algebra homomorphism which is an embedding when restricted to
\begin{gather*}
\1_r\hat{S}(m,N)\1_r\simeq \hat{\mathcal{H}}_{\hat{A}_{r-1}}.
\end{gather*}
Then $f$ is a~$\Q(q)$-algebra isomorphism $A \simeq \hat{S}(m,N)$.
\end{Proposition}

The idea is that it is enough to check the injectivity on the Hecke algebra for deducing it for the whole Schur algebra.
Now, it is proven in~\cite[Theorem~3.4.8]{Green} that the action of~$\hat{S}(m,N)$ on~$V_{\infty}^{\otimes N}$ is faithful.
This is the central ingredient in the following proposition.

\begin{Proposition}%\label{prop:injectivity_gen}
If $N<m$ and $N\leq n$,
the map $n{\bf{\hat{A}Web}}^+_m(N)\mapsto \End_{U_q(\sln)}\big({\bigwedge}_q^N(\C^n\otimes V_{\infty})\big)$ is injective.
\end{Proposition}
\begin{proof}
From~\cite{Green},
the action of $\hat{S}(m,N)\simeq n{\bf \hat{A}Web}^+_m$ on $V_{\infty}^{\otimes N}$ is faithful.
In particular,
for all $x\in \1_N\hat{S}(m,N)\1_N$,
there exists a~vector $\theta$ which is acted on non-trivially.
We have a~map $V_{\infty}^{\otimes N}\mapsto \bigwedge^N(\C^n\otimes V_{\infty})$ sending $X_{j_1}\otimes \cdots \otimes X_{j_N}$
onto $(v_1\otimes X_{j_1})\wedge \cdots \wedge (v_N\otimes X_{j_N})$ (where $v_1,\dots,v_n$ are the generators of $\C^n$).
Then,
$x$ acts (on $\bigwedge_q(\C^n\otimes V_{\infty})$) on the image of $x$ under this map and the result is the image of $x(\theta)$
under the same map.
We would like to ensure that we can f\/ind a~$\theta$ so that this image is non-zero.

To prove this,
assign to $v_i\otimes (X_j\otimes z^r)$,
with $X_j\in V_m$,
the integer $k=i+nj-nmr$.
This def\/ines a~bijective correspondence between $\Z$ and vectors of $\C^n\otimes V_{\infty}$.
Denote $u_k=v_i\otimes (X_j\otimes z^r)$.
Then Uglov~\cite{Uglov99} proves that ordered wedges,
that is wedges of the kind $u_{k_1}\wedge \cdots \wedge u_{k_N}$ with $k_1>k_2>\cdots>k_N$,
form a~basis of the wedge product $\bigwedge^N(\C^n\otimes V_{\infty})$.

Note now that the action of $\hat{S}(m,N)$ on $V_{\infty}^{\otimes N}$ commutes with the multiplication by $z$ on any of the factors.
Hence, up to multiplication of the $\theta$ by $(z^{p_1},\dots, z^{p_N})$,
we can assume that all wedges that appear in $x(\theta)$ are ordered: they all survive (and are linearly independent)
in $\bigwedge^N(\C^n\otimes V_{\infty})$.

We thus obtain that $\1_N\hat{S}(m,N)\1_N$ injects into $\End_{U_q(\sln)}\big({\bigwedge}_q(\C^n\otimes V_{\infty})\big)$,
and using Proposition~\ref{prop:McKThHecke},
we obtain the desired claim.
\end{proof}

As explained earlier,
when trying to extend this result for any $n$ and $N$,
one goes into trouble because the $n$-bounded quotient does not naturally translate in $\bigwedge_q(\C^n\otimes V_{\infty})$.
A natural way to solve this issue is to consider the quotient of ${\bigwedge}_q^N(\C^n\otimes V_{\infty})$ by the following subspace:
\begin{gather*}%\label{wedgeA}
\bigoplus_{\lambda|\1_{\lambda}=0\; \in\hat{\U_q}(\hslm)^{n}}\hat{\U_q}\big(\hslm\big)\1_{\lambda}\cdot{\bigwedge}_q^N\big(\C^n\otimes V_{\infty}\big).
\end{gather*}
This space is clearly a~sub $\hat{\U_q}(\hslm)$ and $U_q(\sln)$ representation,
so the quotient is acted on by both quantum groups.
Furthermore,
the action of $\hat{\U_q}(\hslm)$ descends to an action of $\hat{\U_q}(\hslm)^{n}$,
precisely because the ideal by which we quotient $\hat{\U_q}(\hslm)$ acts by zero.
We denote ${}^a{\bigwedge}_q^N(\C^n\otimes V_{\infty})$ this quotient,
and let ${}^a{\bigwedge}_q(\C^n\otimes V_{\infty})=\bigoplus_N{}^a{\bigwedge}_q^N(\C^n\otimes V_{\infty})$.

It is interesting to note that since the proof of Proposition~\ref{prop:McKThHecke} consists in factorizing elements of the Hecke algebra
through other weights,
it can be extended to the $n$-bounded case as well.
In our case,
we have to consider that if the weight is killed,
then the associated element will be killed.
We can thus state:

\begin{Proposition}
Let $N<m$.
Suppose that $A$ is a~$\Q(q)$ algebra and
\begin{gather*}
f\colon \  \hat{S}(m,N)^{n}\mapsto A
\end{gather*}
is a~surjective $\Q(q)$-algebra homomorphism which is an embedding when restricted to
\begin{gather*}
1_r\hat{S}(m,N)^{n}\1_r\simeq \hat{\mathcal{H}}^{2}_{\hat{A}_{r-1}}.
\end{gather*}
Then $f$ is a~$\Q(q)$-algebra isomorphism $A \simeq \hat{S}(m,N)^{n}$.
\end{Proposition}

However, it is now not clear anymore that the Hecke algebra will inject into:
\begin{gather*}
\End_{U_q(\sln)}\left({}^a{\bigwedge}_q\big(\C^n\otimes V_{\infty}\big)\right).
\end{gather*}
Such a~result would be very interesting,
as it would give to annular webs a~representation-theory-based interpretation in the same f\/lavor as the original def\/inition
of webs as an algebra of intertwiners for minuscule $U_q(\sln)$ representations.
More generally,
it would also be very interesting to better understand the relations between this representation and knot theory,
and to see whether there is any translation in the knot theory side of the more general phenomenons studied in~\cite{Uglov99}.

\subsection{Turning an annular knot to a~ladder}

So far, we saw that to a~ladder annular web-tangle,
we can assign a~$U_q(\sln)$ morphism of tensor product of minuscule representations (possibly tensorized with $\C(q)[z,z^{-1}]$),
whose diagrammatic depiction equals the skein element associated to the web-tangle.
This holds for any annular web-tangle isotopic to a~ladder,
but we can extend the process used in the case of usual webs for turning webs to ladder webs.

We can present any upward annular web-tangle in a~similar form as in~\cite{CKM}:
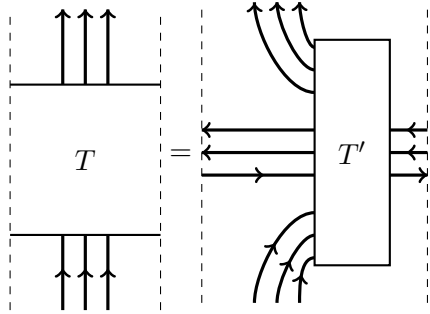
\begin{gather*}
\xy
(0,0)*{
\begin{tikzpicture}
\draw [very thick, directed=.55] (0.2,0) -- (.2,1);
\draw [very thick, directed=.55] (0.5,0) -- (.5,1);
\draw [very thick, directed=.55] (0.8,0) -- (.8,1);
\draw [thick] (-.5,1) -- (1.5,1);
\draw [thick] (-.5,3) -- (1.5,3);
\node at (.5,2) {$T$};
\draw [very thick, directed=1] (0.2,3) -- (.2,4);
\draw [very thick, directed=1] (0.5,3) -- (.5,4);
\draw [very thick, directed=1] (0.8,3) -- (.8,4);
\draw [dashed] (-.5,0) -- (-.5,4);
\draw [dashed] (1.5,0) -- (1.5,4);
\end{tikzpicture}};
\endxy
=\xy
(0,0)*{
\begin{tikzpicture}
\draw [very thick, directed=.55] (0.2,0) .. controls (.2,.4) and (.6,1.2) .. (1,1.2);
\draw [very thick, directed=.55] (0.5,0) .. controls (0.5,.3) and (.8,.9) .. (1,.9);
\draw [very thick, directed=.55] (0.8,0) .. controls (0.8,.2) and (.8,.6) .. (1,.6);
\draw [thick] (1,.5) rectangle (2,3.5);
\node at (1.5,2) {$T'$};
\draw [very thick, directed=.55] (-.5,1.7) -- (1,1.7);
\draw [very thick, directed=1] (2,1.7) -- (2.5,1.7);
\draw [very thick, directed=.55] (2.5,2) -- (2,2);
\draw [very thick, directed=1] (1,2) -- (-.5,2);
\draw [very thick, directed=.55] (2.5,2.3) -- (2,2.3);
\draw [very thick, directed=1] (1,2.3) -- (-.5,2.3);
\draw [very thick, directed=.99] (1,2.8) .. controls (.6,2.8) and (.2,3.6) .. (.2,4);
\draw [very thick, directed=1] (1,3.1) .. controls (.8,
3.1) and (.5,3.7) .. (.5,4);
\draw [very thick, directed=1] (1,3.4) .. controls (.9,3.4) and (.8,3.8) .. (.8,4);
\draw [dashed] (-.5,0) -- (-.5,4);
\draw [dashed] (2.5,0) -- (2.5,4);
\end{tikzpicture}};
\endxy
\end{gather*}

The above pictures are to be understood on an annulus.

Then, taking the Jones--Kauf\/fman product with a~set of well-placed $n$-strands,
we can apply to $T'$ the same process as usual and turn $T$ to a~ladder form.
Note that taking the Jones--Kauf\/fman product with $n$-strands is an invertible process,
by taking the product with the same number of $n$-strands oriented downward and pairing the obtained couples of oppositely oriented $n$-strands.

\section{Categorif\/ication}

From the Dynkin diagram \eqref{DynkinDiagAffine},
one can build a~categorif\/ied quantum group $\Ucat_Q(\hslm)$ fol\-lo\-wing~\cite{CLau},
which generalizes works by Khovanov and Lauda~\cite{KL,KL3,KL2,Lau1}.

Similarly,
it is a~straightforward generalization of $\Bfoam{n}{m}$ to consider foams on an annulus $n{\bf{ABFoam}}_m(N)$:
local generators are the same ones as in the disk case,
but instead of embedding the foams in the thickening of a~disk,
we embed them in the thickening of an annulus.
Then, the main result from~\cite{LQR1} generalizes at no cost:

\begin{Proposition}%\label{prop_LQR1}
For $n=2,3$,
for each $N>0$ there is a~$2$-representation $\Phi_n \colon \Ucat_Q(\hslm) \to n{\bf{ABFoam}}_m(N)$ defined on single strand $2$-morphisms by:
\begin{gather*}
\Phi_n \left(
\xy 0;/r.17pc/:
(0,7);(0,-7); **\dir{-} ?(1)*\dir{>};
(7,3)*{ \scriptstyle \lambda};
(-2.5,-6)*{\scriptstyle i};
(-10,0)*{};(10,0)*{};
\endxy
\right) =
\Efoam[.5]
,
\qquad
\Phi_n \left(
\xy 0;/r.17pc/:
(0,7);(0,-7); **\dir{-} ?(1)*\dir{>};
(0,0)*{\bullet};
(7,3)*{ \scriptstyle \lambda};
(-2.5,-6)*{\scriptstyle i};
(-10,0)*{};(10,0)*{};
\endxy
\right) =
\dotEfoam[.5]
\end{gather*}
on crossings by:
\begin{gather*}
\Phi_n \left(
\xy 0;/r.20pc/:
(0,0)*{\xybox{
(-4,-4)*{};(4,4)*{} **\crv{(-4,-1) & (4,1)}?(1)*\dir{>} ;
(4,-4)*{};(-4,4)*{} **\crv{(4,-1) & (-4,1)}?(1)*\dir{>};
(-5.5,-3)*{\scriptstyle i};
(5.5,-3)*{\scriptstyle i};
(9,1)*{\scriptstyle \lambda};
(-10,0)*{};(10,0)*{};}};
\endxy
\right) =
\crossingEEfoam[.5]
\\
\Phi_n \left(
\xy 0;/r.20pc/:
(0,0)*{\xybox{
(-4,-4)*{};(4,4)*{} **\crv{(-4,-1) & (4,1)}?(1)*\dir{>} ;
(4,-4)*{};(-4,4)*{} **\crv{(4,-1) & (-4,1)}?(1)*\dir{>};
(-5.5,-3)*{\scriptstyle i};
(5.5,-3)*{\scriptstyle \ \ \ i+1};
(9,1)*{\scriptstyle \lambda};
(-10,0)*{};(10,0)*{};}};
\endxy
\right) =
\crossingEoneEtwofoam[.5]
,
\qquad
\Phi_n \left(
\xy 0;/r.20pc/:
(0,0)*{\xybox{
(-4,-4)*{};(4,4)*{} **\crv{(-4,-1) & (4,1)}?(1)*\dir{>} ;
(4,-4)*{};(-4,4)*{} **\crv{(4,-1) & (-4,1)}?(1)*\dir{>};
(-5.5,-3)*{\scriptstyle i+1 \ \ \ };
(5.5,-3)*{\scriptstyle i};
(9,1)*{\scriptstyle \lambda};
(-10,0)*{};(10,0)*{};}};
\endxy
\right) =
\crossingEtwoEonefoam[.5]
\\
\Phi_n \left(
\xy 0;/r.20pc/:
(0,0)*{\xybox{
(-4,-4)*{};(4,4)*{} **\crv{(-4,-1) & (4,1)}?(1)*\dir{>} ;
(4,-4)*{};(-4,4)*{} **\crv{(4,-1) & (-4,1)}?(1)*\dir{>};
(-5.5,-3)*{\scriptstyle j};
(5.5,-3)*{\scriptstyle i};
(9,1)*{\scriptstyle \lambda};
(-10,0)*{};(10,0)*{};}};
\endxy
\right) =
\crossingEthreeEonefoam[.5]
,
\qquad
\Phi_n \left(
\xy 0;/r.20pc/:
(0,0)*{\xybox{
(-4,-4)*{};(4,4)*{} **\crv{(-4,-1) & (4,1)}?(1)*\dir{>} ;
(4,-4)*{};(-4,4)*{} **\crv{(4,-1) & (-4,1)}?(1)*\dir{>};
(-5.5,-3)*{\scriptstyle i};
(5.5,-3)*{\scriptstyle j};
(9,1)*{\scriptstyle \lambda};
(-10,0)*{};(10,0)*{};}};
\endxy
\right) =
\crossingEoneEthreefoam[.5]
\end{gather*}
where $j-i > 1$,
and on caps and cups by:
\begin{gather*}
\Phi_n \left(
\xy 0;/r.20pc/:
(0,0)*{\bbcef{i}};
(8,4)*{\scriptstyle \lambda};
(-10,0)*{};(10,0)*{};
\endxy
\right) =
\capFEfoam[.5]
,
\qquad
\Phi_n \left(
\xy 0;/r.20pc/:
(0,0)*{\bbcfe{i}};
(8,4)*{\scriptstyle \lambda};
(-10,0)*{};(10,0)*{};
\endxy
\right) = (-1)^{a_i}
\capEFfoam[.5]
\\
\Phi_n \left(
\xy 0;/r.20pc/:
(0,0)*{\bbpef{i}};
(8,-4)*{\scriptstyle \lambda};
(-10,0)*{};(10,0)*{};
\endxy
\right) = (-1)^{a_i + 1}
\cupFEfoam[.5]
,
\qquad
\Phi_n \left(
\xy 0;/r.20pc/:
(0,0)*{\bbpfe{i}};
(8,-4)*{\scriptstyle \lambda};
(-10,0)*{};(10,0)*{};
\endxy
\right) =
\cupEFfoam[.5]
\end{gather*}
where in the above diagrams the $i^{th}$ sheet is always in the front.
\end{Proposition}

Then, following~\cite{LQR1}, we can build from any annular knot,
turned into an annular entangled ladder, a complex over the category $\Ucat_Q(\slm)$.
Applying to it $\Phi_n$ (for $n=2$ or $n=3$),
we obtain extensions to the annulus of Khovanov's homology~\cite{Kh1,Kh2,Kh5} built in the spirit of Bar-Natan~\cite{BN2}.
Again, the proof of the invariance relies on checking Reidemeister moves,
which are local and therefore directly extend from the usual case to the af\/f\/ine one.

\subsection*{Acknowledgements}

I would like to thank my advisors Christian Blanchet and Catharina Stroppel for their constant support,
and Aaron Lauda for his great help.
Many thanks also to David Rose,
Peng Shan, Pedro Vaz and Emmanuel Wagner for all useful discussions we had,
Marco Mackaay for pointing out the interest of the af\/f\/ine Hecke algebra,
and especially to Mathieu Mansuy for teaching me everything I know about af\/f\/ine algebras.
I~also wish to acknowledge the great help of the anonymous referees.

\pdfbookmark[1]{References}{ref}
\LastPageEnding

\end{document}